\documentclass[reqno,11pt]{amsart}
\usepackage{amsmath, latexsym, amsfonts, amssymb, amsthm, amscd}
\usepackage[rightcaption]{sidecap}
\usepackage{multirow,stmaryrd}
\usepackage{color}
\usepackage{graphics,epsf,psfrag}
\numberwithin{equation}{section}

\newtheorem{theorem}{Theorem}[section]
\newtheorem{lemma}[theorem]{Lemma}
\newtheorem{proposition}[theorem]{Proposition}

\newtheorem{rem}[theorem]{Remark}

\DeclareMathOperator{\sign}{\mathrm{sign}}

\newcommand{\ind}{\mathbf{1}}

\renewcommand{\tilde}{\widetilde}

\newcommand{\cF}{{\ensuremath{\mathcal F}} }

\newcommand{\cE}{{\ensuremath{\mathcal E}} }
\newcommand{\cH}{{\ensuremath{\mathcal H}} }
\newcommand{\cC}{{\ensuremath{\mathcal C}} }

\newcommand{\cL}{{\ensuremath{\mathcal L}} }

\newcommand{\cW}{{\ensuremath{\mathcal W}} }
\newcommand{\cB}{{\ensuremath{\mathcal B}} }
\newcommand{\cI}{{\ensuremath{\mathcal I}} }
\newcommand{\cZ}{{\ensuremath{\mathcal Z}} }

\newcommand{\bP}{{\ensuremath{\mathbf P}} }
\newcommand{\bE}{{\ensuremath{\mathbf E}} }

\DeclareMathSymbol{\leqslant}{\mathalpha}{AMSa}{"36} % nicer `smaller or equal'
\DeclareMathSymbol{\geqslant}{\mathalpha}{AMSa}{"3E} % nicer `larger or equal'
\DeclareMathSymbol{\eset}{\mathalpha}{AMSb}{"3F}     % nicer `emptyset'
\newcommand{\dd}{\,\text{\rm d}}             % a straight d for differentials
\newcommand{\suptwo}[2]{\sup_{\substack{#1 \\ #2}}} % sup with 2 lines
\newcommand{\sumtwo}[2]{\sum_{\substack{#1 \\ #2}}} % sum with 2 lines
\newcommand{\limtwo}[2]{\lim_{\substack{#1 \\ #2}}}     % \lim with 2 lines
\newcommand{\bbN}{{\ensuremath{\mathbb N}} }

\newcommand{\bbR}{{\ensuremath{\mathbb R}} }

\newcommand{\bbZ}{{\ensuremath{\mathbb Z}} }

\newcommand{\ga}{\alpha}
\newcommand{\gb}{\beta}
\newcommand{\gd}{\delta}
\newcommand{\gep}{\varepsilon}       % \ge already exists...
\newcommand{\gp}{\varphi}

\newcommand{\gl}{\lambda}
\newcommand{\gL}{\Lambda}

\makeatletter
\def\captionfont@{\footnotesize}
\def\captionheadfont@{\scshape}

\long\def\@makecaption#1#2{%
  \vspace{2mm}
  \setbox\@tempboxa\vbox{\color@setgroup
    \advance\hsize-6pc\noindent
    \captionfont@\captionheadfont@#1\@xp\@ifnotempty\@xp
        {\@cdr#2\@nil}{.\captionfont@\upshape\enspace#2}%
    \unskip\kern-6pc\par
    \global\setbox\@ne\lastbox\color@endgroup}%
  \ifhbox\@ne % the normal case
    \setbox\@ne\hbox{\unhbox\@ne\unskip\unskip\unpenalty\unkern}%
  \fi
  \ifdim\wd\@tempboxa=\z@ % this means caption will fit on one line
    \setbox\@ne\hbox to\columnwidth{\hss\kern-6pc\box\@ne\hss}%
  \else % tempboxa contained more than one line
    \setbox\@ne\vbox{\unvbox\@tempboxa\parskip\z@skip
        \noindent\unhbox\@ne\advance\hsize-6pc\par}%
\fi
  \ifnum\@tempcnta<64 % if the float IS a figure...
    \addvspace\abovecaptionskip
    \moveright 3pc\box\@ne
  \else % if the float IS NOT a figure...
    \moveright 3pc\box\@ne
    \nobreak
    \vskip\belowcaptionskip
  \fi
\relax
}
\makeatother
\def\writefig#1 #2 #3 {\rlap{\kern #1 truecm
\raise #2 truecm \hbox{#3}}}

\newcommand{\tf}{\textsc{f}}
\newcommand{\tg}{\textsc{g}}

\begin{document}

\title[Generalized Poland-Scheraga model and two-dimensional renewal processes]{Generalized Poland-Scheraga denaturation model and two-dimensional renewal processes}

\author{Giambattista Giacomin}
\address{
  Universit\'e Paris Diderot, Sorbonne Paris Cit\'e,  Laboratoire de Probabilit{\'e}s et Mod\`eles Al\'eatoires, UMR 7599,
            F- 75205 Paris, France
}

\author{Maha Khatib}
\address{
Universit\'e Paris Diderot, Sorbonne Paris Cit\'e,  Laboratoire de Probabilit{\'e}s et Mod\`eles Al\'eatoires, UMR 7599,
            F- 75205 Paris, France
            }

\begin{abstract}
The Poland-Scheraga  model describes the denaturation transition of two complementary -- in particular, equally long --  strands of DNA, and it has enjoyed a remarkable success both for quantitative modeling purposes and at a more theoretical level. The solvable character of the homogeneous version of the model is 
one of  features to which its success is due.
In the bio-physical literature a generalization of the model, allowing different length and non complementarity of the strands,
  has been considered and the solvable character extends to this substantial generalization.  
  We present a mathematical analysis of the homogeneous generalized Poland-Scheraga model. Our approach is based on the fact that such a model is  a homogeneous pinning model based on a bivariate renewal process, much like the basic Poland-Scheraga model is a pinning model based on a univariate, i.e. standard, renewal. We present a complete analysis of the free energy singularities, which include  the localization-delocalization critical point and (in general)  other critical points that have been only partially captured in the physical literature.  We  obtain also precise estimates on 
  the path properties of the model. 
   \\[10pt]
  \\[10pt]
\end{abstract}

\maketitle

%\tableofcontents

\section{Introduction and main results}
\label{sec:intro}

\subsection{General overview}
The localization-delocalization phenomenon in various polymer models has been the object of much attention
  in the physics, biophysics and mathematics literature \cite{cf:Fisher,cf:GB,cf:dH,cf:KMP,cf:PSbook}. 
  One of the main biological and physical phenomenon that motivates this work  is
    DNA denaturation, that is the separation of the two DNA strands at high temperature and, more generally, the fluctuation phenomena observed below denaturation, when the two strands are tied together. 
The most basic and studied model in this field 
is  the Poland-Sheraga (PS) model \cite{cf:PSbook} which is limited to the case of sharp
 complementarity of  two equal length strands: only bases with the same index can form pairs. 
   From the theoretical physics and mathematical viewpoint what is most remarkable in the homogeneous version of the  model is its solvable character and the fact that at the delocalization (or denaturation) transition the behavior -- i.e. the critical behavior -- can be fully captured. 
A mathematical viewpoint on this solvable character and on the solution itself is that the Poland-Scheraga model is 
a Gibbs measure with only one body potentials and built on a one-dimensional process 
 \cite{cf:GB}. 
 
In \cite{cf:GO0,cf:GO,cf:NG} (see also \cite{cf:litan,cf:SNR} for a primitive version), a generalization of the Poland-Sheraga (gPS) model  has been introduced and the novelties are:
\begin{itemize}
\item The possibility of formation of non-symmetrical loops in the two strands (i.e., the contribution to a loop, in terms of number of bases, from the two strands is not necessarily the same).
\item  The two strands may be of different lengths.
\end{itemize}
These novelties are very substantial (we invite the reader to compare Figure~\ref{fig:dna} with \cite[Fig.~6]{cf:Fisher} or  \cite[Fig. 2.5]{cf:G}). Nevertheless, as already pointed out in \cite{cf:GO,cf:NG,cf:SNR}, the 
solvable character is preserved. However, that the novelties are really substantial is  witnessed by a  richer phenomenology (partially captured and understood in \cite{cf:EON,cf:NG}): in addition to the expected denaturation transition, the gPS model 
displays other transitions.  

Here we develop a mathematical analysis of the gPS model based on the observation 
that it  is a pinning model based on a  two-dimensional  renewal process. Much like for the original PS model, tools
from Renewal Theory allow going far toward a complete  understanding of the model. Nevertheless, as we will explain,
some important questions are still open and they correspond to open problems in the theory of two and higher dimensional renewal processes.

\subsection{The gPS model: biophysics version}
\label{sec:biop}
This subsection, as well as \S~\ref{sec:matching}, can be skipped if one is not focusing on the biophyisics set-up. 
The model we consider has been introduced in \cite{cf:GO}. The two DNA strands, of lengths $M$ and $ N\ge 1$ -- the length of course corresponds to the number of base bases  --
interact by forming some base pairs. We talk of $N$-strand, $M$-strand and of base $i$ of the $N-$ or $M-$strand
with the obvious meaning. An allowed configuration of our system
 is a collection
of base pairs 
\begin{equation}
\left( (i_1,j_1), (i_2,j_2), \ldots, (i_n,j_n)\right) \in \bbN^{2n}\, ,  \ \ \text{ with } \bbN=\{1,2, \ldots\}\, ,
\end{equation}
where $n\in \bbN:=\{1, 2, \ldots, N\}$ and 
\medskip

\begin{enumerate}
\item $(i_1, j_1)=(1,1)$;
\item  $i_k< i_{k'}$, as well as $j_k< j_{k'}$, for $1\le k < k' \le n$.
\end{enumerate}
\medskip

\begin{figure}[h]
\begin{center}
\leavevmode
\epsfxsize =12 cm
\psfragscanon
\psfrag{1}[c][l]{\tiny $1$}
\psfrag{2}[c][l]{\tiny $2$}
\psfrag{3}[c][l]{\tiny $3$}
\psfrag{4}[c][l]{\tiny $4$}
\psfrag{5}[c][l]{\tiny $5$}
\psfrag{6}[c][l]{\tiny $6$}
\psfrag{7}[c][l]{\tiny $7$}
\psfrag{8}[c][l]{\tiny $8$}
\psfrag{9}[c][l]{\tiny $9$}
\psfrag{d}[c][l]{\tiny $10$}
\psfrag{o}[c][l]{\tiny $12$}
\psfrag{q}[c][l]{\tiny $14$}
\psfrag{n}[c][l]{\tiny $19$}
\epsfbox{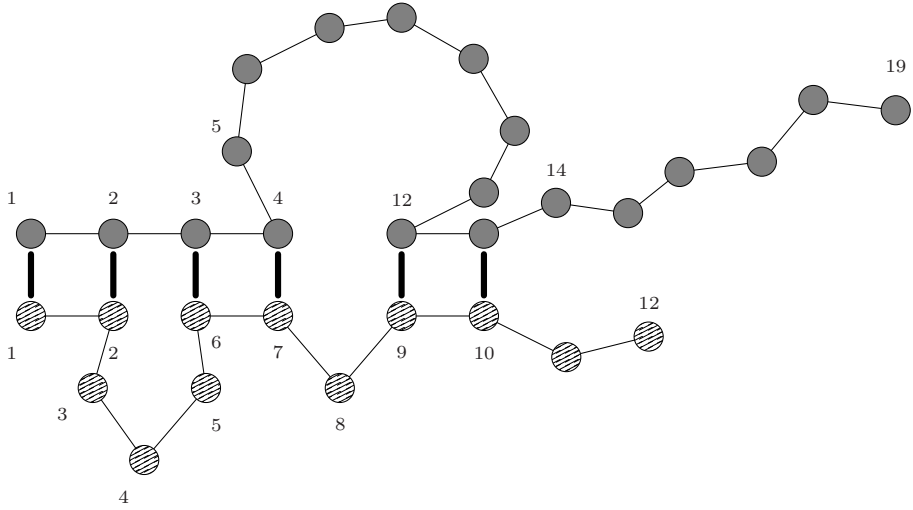}
\end{center}
\caption{\label{fig:dna} A representation of a trajectory of the  gPS model in the biophysics representation. The first strand  contains $12$ bases, the second strand $19$.
The six base pairs determining the configuration are $(1,1)$, $(2,2)$, $(6,3)$, $(7,4)$, $(9,12)$ and $(10,13)$.}
\end{figure}

The  first condition is simply saying that the first two bases form a pair and the second condition
is  imposed by the geometric constraint (see Figure~\ref{fig:dna}). 
The weight of every configuration is assigned by the following rules: 
\medskip

\begin{enumerate}
\item Each base pair is energetically favored and carries an energy $-E_b<0$;
\item A base which is not in pair it is either in a loop or in the free ends:
\begin{itemize}
\item It is in a loop if it is in $L_k:=\{i: \, i>i_k, \, i< i_{k+1}\}\cup \{j: \, j>j_k, \, j< j_{k+1}\}$
for some $k\in \{1, \ldots, n-1\}$: the loop $L_k$ has length $\ell_k:=(i_{k+1}-i_k)+ (j_{k+1}-j_k)-2$
and we associate to $L_k$ an entropy factor $B(\ell_k)$ with
\begin{equation}
B(\ell)\, :=\, s^\ell \ell^{-c}\, ,
\end{equation}
where $s\ge 1$ and $c>2$. There is also an energetic $E_l>0$ penalty associated to a loop.
\item The free ends have length $N-i_n$ and $M-j_n$ and to each free end
we associate the entropy term $A(\ell):=s^\ell (\ell+1)^{-\bar c}$ where $\bar c$ is another positive constant. 
\end{itemize}
\end{enumerate}
\medskip

As we will see  the value of $s$ is irrelevant.
The value of $\bar c$, chosen equal to $0.1$ in \cite{cf:NG}, is somewhat more relevant, but
what is very relevant is
the value of $c$: in \cite{cf:NG} is chosen equal to $2.15$.

These rules easily lead to a formula for the partition function, i.e. the sum of the weights over all allowed configurations,
of our system 
\begin{equation}
\label{eq:ZNM}
Z_{N}^M\, :=\, \sum_{i=0}^{N-1}\sum_{j=0}^{M-1} A(i) A(j) 
W_{N-i}^{M-j}\, ,
\end{equation}
where $W_l^r$ obeys the recursion relation ($\gb \ge 0$ is  proportional to the inverse of the temperature)
\begin{equation}
W_{m+1}^{r+1}\, =\, \exp(\gb E_b) W_{m}^{r}+  \exp(\gb (E_b-E_l))
\sumtwo{i,i':\, i+i'>0}{i<m, \, i'<r} B(i+i') W_{m-i}^{m-i'}\, ,
\end{equation}
with $W_1^1=1$ and $W_1^i=W_i^1$ for $i>1$. 

%We are interested in thermodynamical quantities and, notably, in the free energy density 
%\begin{equation}
%\label{eq:F0}
%F_\gamma (\gb)\, := \limtwo{N,M \to\infty:}{\frac M N \to \gamma} \frac 1N \log Z_N^M\, ,
%\end{equation}
%for $\gamma\ge 1$. Of course $F_\gamma (\gb)$ captures the exponential, or Laplace, asymptotic behavior of $Z_N^M$.

\subsection{The gPS model: renewal process viewpoint}
From a mathematical viewpoint we  take a more general viewpoint and we introduce a 
two-dimensional renewal pinning model. A discrete two-dimensional renewal issued from the origin is a random walk
$\tau=\{\tau_n\}_{n=0,1, \ldots}=(\tau^{(1)}, \tau^{(2)})=\{ (\tau^{(1)}_n, \tau^{(2)}_n)\}_{n=0,1, \ldots}$
where $\tau_0=(0,0)$ and, for $n\in \bbN:=\{1,2, \ldots\}$, $\tau_n$ is a sum of $n$ independent identically distributed 
random variables taking values in $\bbN^2$. So if we set $\mathrm{K}(n,m):= \bP(\tau_1=(n,m))$ then given
$\{(i_n,j_n)\}_{n=0,1,2,\ldots}$, with $(i_0,j_0)=(0,0)$, for every $k\in \bbN$
\begin{equation}
\bP\left(\tau_n= (i_n,j_n) \text{ for } n=1,2, \ldots, k\right)\, =\, \prod_{n=1}^k \mathrm{K}\left(i_n-i_{n-1},j_n-j_{n-1}\right)\, ,
\end{equation}
and, by construction, such a probability is zero unless the $i$'s and $j$'s form strictly increasing sequences.

We can then introduce for given $N$ and $M\in \bbN$ a 
pinning model  of length $(N, M)$ by forcing, i.e. conditioning, $\tau$ to visit $(N,M)$ and by penalizing ($h\le 0$)
or rewarding ($h\ge 0$) the number of renewals up to $(N,M)$. More formally,  we introduce the probability measure
$\bP_{N, M, h}$ by setting for every $k \in \bbN$ such that $k\le  \min(N,M)=: N\wedge M $ 
 and for every $\{(i_n,j_n)\}_{n=0,1, \ldots,k}$ with $(i_0,j_0)=(0,0)$,
$i_n-i_{n-1}>0$ as well as $j_n-j_{n-1}>0$ for $n=1, \ldots, k$ and $(i_k,j_k)=(N,M)$
\begin{equation}
\label{eq:Pc}
\frac{
\bP_{N, M, h}\left(\tau_n= (i_n,j_n) \text{ for } n=1,2, \ldots, k\right)}
{\bP\left(\tau_n= (i_n,j_n) \text{ for } n=1,2, \ldots, k\right)}\, :=\,
\frac 1{\cZ_{N, M, h}} \exp\left( h k\right)\, ,
\end{equation}
where $\cZ_{N, M, h}$ is the partition function (or normalization constant):
\begin{equation}
\cZ_{N, M, h}\, :=\, \bE \left[ 
\exp\left(
h \vert \tau \cap ([1, N] \times [1, M]) \vert 
\right); 
(N, M) \in \tau 
\right]\, ,
 \end{equation} 
 in which we are interpreting $\tau$ as a random subset of $\bbN^2$ and $\vert \cdot\vert $  denotes the cardinality.
Note that $\cZ_{0, 0, h}=1$, as well as $\cZ_{0, M, h}=\cZ_{N, 0, h}=0$.
Of course 
$\bP_{N, M, h}$ requires $\cZ_{N, M, h}\neq 0$ and 
whether this is the case or not depends on the inter-arrival distribution $\mathrm{K}(n,m)$ and, possibly, on $N$ and $M$.

\begin{SCfigure}[50]
\centering
\caption{\label{fig:dna-ren} A representation of a trajectory of the gPS model in the renewal process representation, corresponding to Figure~\ref{fig:dna}. 
The base pairs are now renewal points of a two-dimensional discrete renewal process: these points correspond to the (six) base pairs 
of Figure~\ref{fig:dna}, except that the first base pair is now $(0,0)$ and also all the other ones are translated down of $(1,1)$ with respect 
to Figure~\ref{fig:dna}. The renewal trajectory is drawn up to the renewal point $(9,12)$ and the trajectory up to this point
correspond to one of the possible trajectories of $Z^c_{9,12, h}$. The free ends, of lengths $2$ and $6$, are then represented as straight
lines that go till the boundary of the rectangle with opposite vertices $(0,0)$ and $(11,18)$. $Z^f_{11,18, h}$ is obtained by summing up 
with respect to the position of the last renewal point -- $(9,12)$ in this example -- the contribution of the constrained partition function
times the contribution due to the two free ends. 
}
\leavevmode
\epsfxsize =5 cm
\psfragscanon
\psfrag{0}[c][l]{\tiny $0$}
\psfrag{1}[c][l]{\tiny $1$}
\psfrag{9}[c][l]{\tiny $9$}
\psfrag{11}[c][l]{\tiny $11$}
\psfrag{12}[c][l]{\tiny $12$}
\psfrag{18}[c][l]{\tiny $18$}
\psfrag{base}[l][l]{\tiny $=$ base pair}
\psfrag{free}[l][l]{\tiny $=$ free end}
\epsfbox{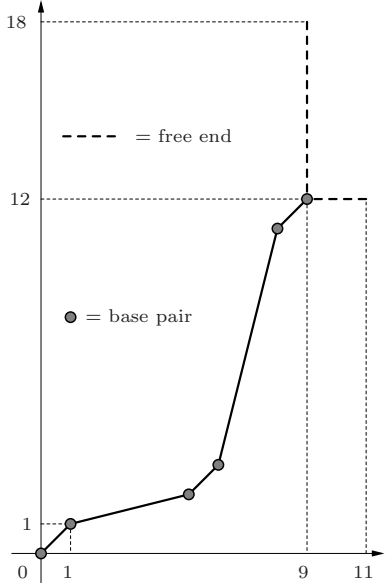}
\end{SCfigure}

\medskip

We have used the atypical notation $\cZ$ instead of $Z$ because the latter is going to be employed for the
model on which we really focus: we consider in fact a very special choice of the {\sl inter-arrival distribution} $\mathrm{K}(\cdot, \cdot)$, namely
$\mathrm{K}(n, m)=K(n+m)$ where 
 $K:\{2,3, \ldots\}\to (0, \infty)$ and
\begin{equation}
K(n) \, := \, \frac{L(n)}{n^{1 +\ga}} \,,
\end{equation}
where $L(\cdot)$ is a slowly varying function and $\ga \ge 1$ (see Appendix~\ref{sec:A} for the properties of slowly varying functions). With this definition, we have that $K(n)>0$ for every $n$: all statements generalize to 
the case in which $K(n)=0$ for finitely many $n$, but we make this choice to make a few proofs more concise. 
 We require $\sum_{n,m\in \bbN} K(n+m)=1$ and of course $\sum_{n,m\in \bbN} K(n+m)=
\sum_{m=1}^\infty m K(m+1)$.
We introduce then the {\sl constrained partition function}
$Z_{N, M,h}^c$ which coincides with $\cZ_{N, M, h}$ once the specific choice of the inter-arrival is made.
In an alternative explicit fashion 
\begin{equation}
\label{eq:Zc}
Z_{N, M,h}^c\, :=\, \sum_{n=1}^{N\wedge M} \sumtwo{\underline{l}\in \bbN ^n:}{\vert \underline{l}\vert =N}
\sumtwo{\underline{t}\in \bbN ^n:}{\vert \underline{t}\vert =M} \prod_{i=1}^n \exp(h) K\left( l_i +t_i\right)\,.
\end{equation} 
The {\sl free partition function} is then defined by
\begin{equation}
\label{eq:Zf}
Z_{N, M,h}^f\, :=\, \sum_{i=0}^N \sum_{j=0}^M K_f(i) K_f(j) Z_{N-i, M-j,h}^c\, ,
\end{equation}
where $K_f : \{0\}\cup \bbN \to (0, \infty)$ is defined as $K_f(n) := \overline{L}(n) / n^{\overline{\ga}}$ for every $n \ge 1$ and $K_f(0) = 1$ (an arbitrary choice: there is no loss of generality with respect to requiring just $K_f(0)>0$ and, once again, one
can even allow $K_f(n)=0$ for finitely many $n$, but we choose positivity for conciseness) with $\overline{\ga} \in \bbR$. The free partition function is the normalization associated to the probability
$\bP^f_{N, M, h}$ 
defined 
by setting for every $k \in \bbN$ such that $k\le  \min(N,M)$ 
 and for every $\{(i_n,j_n)\}_{n=0,1, \ldots,k}$ with $(i_0,j_0)=(0,0)$,
$i_n-i_{n-1}>0$ as well as $j_n-j_{n-1}>0$ for $n=1, \ldots, k$ and $i_k$ (respectively,  $j_k$)
that does not exceed $N$ (respectively, $M$)
\begin{multline}
\label{eq:Pf}
%\bP^f_{N, M, h}\left(\tau_n= (i_n,j_n) \text{ for } n=1,2, \ldots, k, \text{ and } \vert \tau \cap [1,N]\times [1,M]\vert =k \right)}
\bP^f_{N, M, h}\left(\tau  \cap [1,N]\times [1,M] \, =\, \left\{ (i_1, j_1),  (i_2, j_2), \ldots,  (i_k, j_k)\right\}
\right)
\, :=\\
\frac {\exp\left( h k\right)}{Z^f_{N, M, h}} 
K_f(N- i_k) K_f(M-j_k){\bP\left(\tau_n= (i_n,j_n) \text{ for } n=1,2, \ldots, k\right)}
\, .
\end{multline}

Here we introduce also the free energy density (the existence of the limit is proven in Section~\ref{sec:free})
\begin{equation}
\label{eq:F}
\tf_\gamma (h)\, := \limtwo{N, M \to\infty:}{\frac M N \to \gamma} \frac 1N \log Z^f_{N, M,h}\, ,
\end{equation}
for $\gamma>0$ or, as we will often do without loss of generality, for $\gamma\ge 1$.
The limit in  \eqref{eq:F} means: for every $\{(N_n, M_n)\}_{n=1,2, \ldots}$ with $\lim_n M_n/N_n=\gamma$
and $\lim_n N_n=\infty$.

\subsection{Matching biophysics and renewal process viewpoints} 
\label{sec:matching}
The first remark to make on the biophysics model is  that the dependence on $s$
of $Z_N^M$ (cf. \eqref{eq:ZNM}) is trivial:  $Z_N^M$ for a given value of $s$ coincides with 
$s^{N+M} \tilde Z_N^M$, where $ \tilde Z_N^M$ is $Z_N^M$ with $s=1$ and both $E_b$
and $E_l$ are decreased  of $2\log s$. We can therefore set $s=1$ without true loss of generality.
We then remark that we can match $\exp(-\gb E_b) Z_{N+1}^{M+1}$ and $Z^f_{N, M, h}$ by
an appropriate choice of $h$, $K(\cdot)$ and $K_f(\cdot)$: the fact that we consider the biophysics model
with lengths augmented by $1$ and renormalized by the factor $\exp(-\gb E_b)$  is just because in parallel with one-dimensional renewal pinning works 
we have chosen to start the renewal from {\sl time} (or {\sl renewal epoch}) zero, but without 
given an energetic reward to the base pair $0$. For matching is then made by observing that:
\medskip

\begin{enumerate}
\item The match for the free end terms $K_f(\cdot)$ and $A(\cdot)$ is easily made.
\item A {\sl minimal} inter-arrival step, that is $(1,1)$, in the renewal model corresponds to a base pair and contributes 
$K(2)\exp(h)$. It would then be matched to
$\exp(\gb E_b)$ in the biophysics model.
\item All other inter-arrival steps $(i,i')$, i.e. $i+i'=3,4,\ldots$, which give a contribution
$K(i+i')\exp(h)$ correspond to loops with $i+i'-2$ unpaired bases and the contribution in the biophysics
model is $\exp(\gb(E_b-E_l)) B(i+i'-2)$.
\end{enumerate}

\medskip

Then of course $c=1+\ga$ and $\bar c= \bar \ga$ and, at last,
the matching between the two models is done thanks to the ample freedom that we have in choosing 
$K(\cdot)$, except that we have required that $\sum_{n, m } K(n+m)=1$. This {\sl probability constraint}
in reality just corresponds to a shift in the parameter $h$. 
So renewal models include the biophysics ones: from a qualitative view point the matching
is immediate, form the quantitative one it requires some bookkeeping care. Explicit matchings are 
presented in Section~\ref{sec:ex}.

\subsection{The free  energy and the localization transition}  
We start from:

\medskip

\begin{proposition}
\label{th:propo1}
Recall \eqref{eq:F}. We have
\begin{equation}
\label{eq:Fc}
 \limtwo{N, M \to\infty:}{\frac M N \to \gamma} \frac 1N \log Z^c_{N, M,h}\, =\, \limtwo{N, M \to\infty:}{\frac M N \to \gamma}
 \frac 1N \log
  Z^f_{N, M,h}\, .
\end{equation}
\end{proposition}

\medskip

In practice, $ Z^c_{N, M,h}$ is a more fundamental quantity for our computations and we will first identify the free energy density
by looking at the exponential growth of this quantity and only after, in Section~\ref{sec:free}, we will match it with the 
exponential growth rate in the free case. 

\medskip

Note that we are just speaking of exponential {\sl growth} and not of {\sl decrease}. In fact $\tf_\gamma(h)\ge 0$ simply because
$Z^f_{N,M, h}\ge K_f(0)^2 Z^c_{N,M, h}=Z^c_{N,M, h}\ge \exp(h) K(N+M)$. As a matter of fact, this is a very important
issue because it is natural to set
\begin{equation}
\label{eq:defhc}
h_c\, :=\, \inf\{ h: \, \tf_\gamma (h)>0\}\, =\,  \max\{ h: \, \tf_\gamma (h)=0\}\, ,
\end{equation}
where the equality on the right comes from the fact that 
$\tf_\gamma(\cdot)$ is locally bounded, convex (hence continuous) and non-decreasing. 
These facts are evident
from the definitions, like the following two preliminary observations:
\medskip

\begin{enumerate}
\item
$ Z^c_{N, M,h}\le 1$ for $h\le 0$: hence $h_c\ge 0$;
\item we will see just below that $h_c=0$, but it is worth  pointing out that $h_c< \infty$ by elementary arguments. 
For example:
$h_c \le -\log K(2)$ because $Z^f_{N,M,h}\ge K_f(0) K_f(M-N) Z^c_{N, N, h}$ and 
$Z^c_{N, N, h}\ge (\exp(h)K(2))^N$. 
\end{enumerate}

\medskip

%Hence, by convexity,  $\tf_\gamma(h)=0$ if and only if  $h \le h_c$. Therefore
From \eqref{eq:defhc} we readily see that
$h_c$ is a non analyticity point of $\tf_\gamma(\cdot)$ and there is a phase transition of the system. 
It is easy to realize that this transition is  the denaturation (or localization/delocalization) transition:
$\partial_h \tf_\gamma(h)$ -- in case to be interpreted as, say, left derivative, but we will soon see that  
$\partial_h \tf_\gamma(h)$ exists except, in some cases, at $h=h_c$ -- is the density of base pairs (or contact fraction), which
is therefore positive, respectively zero, for $h>h_c$, respectively $h<h_c$.

The next result is much more quantitative about this transition: let us remark that, again by an elementary argument,
$\tf_\gamma(h) \ge  \tf_1(h)$ for every $\gamma\ge 1$.
All asymptotic statements in the next theorem are for $h \searrow 0$:

\medskip

\begin{theorem}
\label{th:F}
For every $\ga\ge 1$ and $\gamma\ge 1$ we have
$h_c=0$ and there exists $h_{c, \gamma}\in (0, \infty]$ such that $\tf_\gamma(\cdot)$ is real analytic in $(-\infty,0)\cup(0,h_{c, \gamma})$. If $h_{c, \gamma}< \infty$ then $h_{c, \gamma}$ is a non analyticity point.  
Moreover  if  $\sum_{n} n^2 K(n)< \infty$, a condition implied by $\ga>2$, we have 
\begin{equation}
\label{eq:fch}
\tf_\gamma(h) \sim \tf_1 (h) \, \sim \, c\,  h \, ,
\end{equation}
with $c^{-1}:=  {\frac 12\sum_{n=2}^\infty n(n-1) K(n)}$. If instead $\sum_{n} n^2 K(n)= \infty$, implied by  $\ga \in [1,2)$, 
there exists $c_{\ga,\gamma}\ge 1$ such that
\begin{equation}
\label{eq:F-GB}
\tf_\gamma(h) \sim c_{\ga,\gamma} \tf_1 (h)\ 
\ \text{ and } \ \tf_1 (h)\, \sim  \, 
L_\ga(h) h^{1/(\ga-1)} \end{equation}
where 
 $L_\ga(\cdot)$ is slowly varying  at $0$. In the case $\ga=1$, \eqref{eq:F-GB} should be interpreted 
 as $\tf_\gamma(h)=O(h^{1/\gep})$ for every $\gep>0$.
\end{theorem}
\medskip

In Section~\ref{sec:deloc}  $c_{\ga, \gamma}$ and $L_\ga(\cdot)$ are determined. The expression of $c_{\ga, \gamma}$
implicitly contains  nontrivial information on the system, see Proposition~\ref{th:Fprop}.

\medskip

We will see  that it may be that $h_{c,\gamma}=\infty$, for example  $h_{c,1}=\infty$ in full generality, but when 
$h_{c,\gamma}<\infty$, $h_{c,\gamma}$ may not be the only  critical point inside the localized regime (Theorem~\ref{th:2}). This means  that there is  more than one localized phase 
in the system: this is what we treat next, but we need to introduce more concepts and definitions.
By doing so we will  start outlining  the proof of Theorem~\ref{th:F}.

\subsection{Transitions in the localized regime}
\label{sec:fe2}
A crucial elementary observation is that
  \eqref{eq:Zc} can be  written as 
\begin{equation}
\label{eq:Zc-2}
\begin{split}
Z_{N, M,h}^c\, &=\, 
\exp\left((N+M) \tg\right)\sum_{n=1}^{N\wedge M} \sumtwo{\underline{l}\in \bbN ^n:}{\vert \underline{l}\vert =N}
\sumtwo{\underline{t}\in \bbN ^n:}{\vert \underline{t}\vert =M} \prod_{i=1}^n \tilde 
K_h\left( l_i ,t_i\right)\\
& =\, 
\exp\left((N+M) \tg\right) \bP\left( (N,M) \in \tilde \tau_ h \right)\, , 
\end{split}
\end{equation} 
where 
\smallskip

\begin{itemize}
\item
$\tilde K_h(n,m)= \exp(h -  (n+m) \tg)K(n+m)$ -- note that $\tilde K_h(n,m)$ is just function of $n+m$ --
and $\tg= \tg(h)$ is the only solution to 
\begin{equation}
\label{eq:hc}
\sum_{n,m}K(n+m)\exp(h-(n+m) \tg)\, =\, 1\, ,
\end{equation}
when such a solution exists (that is, when $h \ge 0$), and $\tg=0$ otherwise. We have therefore defined 
a function $\tg: \bbR \to [0, \infty)$.
\item 
   $\tilde \tau_h $ is the two-dimensional renewal  issued from $(0,0)$
with inter-arrival distribution  $\tilde K_h$: if $h<0$ then $\tilde K_h: \bbN^2\to [0, 1)$ is a sub-probability that we make a probability by defining 
$\tilde K_h(\infty):=1-\sum_{(n, m)\in \bbN^2} \tilde K_h(n,m)$ and, in this case, $\tilde \tau_h $ contains  a.s. a finite number of points (and $\{\infty\}$) and we refer to it as a {\sl terminating} renewal.
\end{itemize}

\medskip

$\tg(h)$ accounts for part of the free energy of the system and some basic features
are straightforward:

\medskip

\begin{lemma}
\label{th:g}
The function  $\tg(\cdot)$ 
is convex and real analytic except at $h=0$. 
\end{lemma} 
\medskip

\noindent
{\it Proof.}
First of all the function is well defined for $h \ge 0$ because the function
$\tg\mapsto \sum_{n,m}K(n+m) \exp(-(n+m)\tg)$, real analytic on the positive semi-axis, decreases from $1$ to $0$ as $\tg$ goes from $0$ to $\infty$.
In particular, $\tg(0)=0$. Analyticity
of $\tg(\cdot)$ on the positive semi-axis follows directly by the Inverse Function Theorem (for analytic functions: see for example \cite[Sec.~2.5]{cf:Primer}). Convexity can be proven directly
by differentiating the implicit expression \eqref{eq:hc}.  A less tedious proof of convexity 
can be achieved by recognizing that $\tg(h)$ is a free energy (see \eqref{eq:tflim0} below).
\qed 

\medskip

In view of \eqref{eq:Zc-2} we see that (limits are with $M\sim \gamma N$)
\begin{equation}
\label{eq:prelimZ}
\lim \frac 1N \log Z^c_{N,M, h} \, =\, (1+\gamma) \tg(h)
+ \lim \frac 1N \log 
\bP\left( (N, M)\in \tilde \tau_h \right)\, .
\end{equation}
Let us focus on $h>0$. 
The inter-arrivals of $\tilde \tau_h$ are exponentially integrable, in particular they have  finite mean: 
\begin{equation}
\label{eq:mu12}
\mu_h\, =\ \left(\mu_h^{(1)}, \mu_h^{(2)} \right) \, :=\,
 \bE [(\tilde\tau_h)_1]\, =\, \left(\sum_{n, m} n \tilde K_h(n,m), \sum_{n, m} m \tilde K_h(n,m)\right)\, . 
 \end{equation}
 Evidently $\mu_h^{(1)}=\mu_h^{(2)}$ and the Law of Large Numbers directly implies that $\tilde \tau_h$ stays close to the
 diagonal of the first quadrant. But then the event $\{ (N, M)\in \tilde \tau_h\}$ is a Large Deviation event for $\gamma \neq 1$
 and it contributes to the free energy: the singularities of the free energy for $h>0$ come from 
 this extra Large Deviation contribution and key word to understand these new transitions 
 is  \emph{Cram\'er regime}. The point in fact is whether or not  the Large Deviation event can be made typical
 by an exponential change of measure (a \emph{tilt}): the larger $\gamma$ is the more is possible 
 that a tilt does not suffice and the typical Large Deviation trajectories will not correspond to a tilt of the measure
 (in this case we say that we are outside of the Cram\'er regime). 
 On the other hand, the interaction strength directly impacts whether or not the process is in the Cram\'er regime.
 The formulas that follow, though probably rather obscure at this stage, precisely characterize the switching between
 Cram\'er  and non Cram\'er regimes.

\medskip

We introduce the convex function $q_h: \bbR^2 \to (0, \infty]^2$
\begin{equation}
\label{eq:q_h}
q_h( \gl)\,=\, q_h(\gl_1,\gl_2)\, :=\, \sum_{n,m} e^h K(n+m) 
\exp\left(-(\tg -\gl_1)n -(\tg -\gl_2)m \right)\, ,
\end{equation}
which is bounded in $(-\infty, \tg]^2$ and it is analytic in the interior of this domain. 
We set for $h>0$
\begin{equation}
\label{eq:gl1bar}
\overline {\gl}_1 (h)\, :=\, \sup\left\{ \gl_1<0:\, q_h( \gl_1, \tg (h))\le 1\right\}\, ,
\end{equation}
and,  since $q_h( \gl_1, \tg (h))$ increases continuously in $\gl_1$ from $q_h(-\infty, \tg (h))=0$
to $q_h(0, \tg (h))>1$,  $\gl_1 (h)$ is negative and it is characterized by  
$q_h(\overline{\gl}_1, \tg (h))=1$. Finally, we set, always for $h>0$
\begin{equation}
\label{eq:gammah}
\gamma_c (h)\, :=\,
\frac{\sum_{n,m} m K(n+m)\exp\left(- n \left(\tg(h) -\overline{\gl}_1(h) \right)\right)}{\sum_{n,m} n K(n+m)\exp\left(- n \left(\tg(h) - \overline{\gl}_1(h) \right)\right)}\, ,
\end{equation}
and both denominator and numerator are bounded because for $c>0$ if $\ga >1$
\begin{equation}
\label{eq:c>0}
\sum_{n,m} (n+m) K(n+m)e^{- cn}\, = \, \sum_{t=2}^\infty t K(t) \sum_{n=1}^{t-1} \exp(-cn)\, \le \frac{1}{e^c-1} 
\sum_{t=2}^\infty t K(t)\, < \infty\, .
\end{equation}
If $\ga=1$ it is easy to see that the denominator is bounded, but the numerator is $+\infty$ for every $h$,
so $\gamma_c(h)=\infty$ in this case. 
Here are some properties (see Section\ref{sec:deloc} for the proof):
\medskip

\begin{lemma}
\label{th:gammac}
Choose $\ga>1$. The function $\gamma_c:(0, \infty)\longrightarrow  (1, \infty)$ is real analytic  and 
\begin{equation}
\label{eq:gammac-th}
\gamma_c(0)\, 
:=\, \lim_{h \searrow 0} \gamma_c (h)\, =\, \frac 1{\ga-1} \vee 1\
 \ 
\textrm{ and } \ 
\gamma_c(\infty)\, =\, 
 \frac{\sum_m m K(1+m)}{\sum_m  K(1+m)}\,.
\end{equation}
\end{lemma}
\medskip

The examples worked out in  Section~\ref{sec:ex} show that $\gamma_c(\cdot)$ can have various behaviors: in particular, in general it is not monotonic. 

%Recall the notation we use for the renewal  $\tau =\{(0,0),(\tau^{(1)}_1, \tau^{(2)}_1), (\tau^{(1)}_2, \tau^{(2)}_2), \ldots \}$, so $\tau^{(1)}=\{0, \tau^{(1)}_1, \tau^{(1)}_2, \ldots\}$.

\medskip

\begin{theorem}
\label{th:2}
Fix $\gamma\ge 1$.  
$\tf_\gamma (\cdot)$ is analytic on  $\{h:\, h>0$ such that $\gamma_c(h)-\gamma \neq 0\}$
and $\tf_\gamma (\cdot)$ is not analytic for the values  $h>0$ at which $\gamma_c(h)-\gamma$
changes sign. However,
 $\tf'_\gamma(\cdot)$ is continuous  on the  positive semi-axis
\end{theorem}

\medskip

Theorem~\ref{th:2} is just a sample of the results we have and that can be gotten on these transitions 
that are transitions between localized regimes, because, by convexity of the free energy, the expected number of contacts does not decrease 
in $h$. In particular the \emph{tangential case} -- when $h \mapsto \gamma_c(h)-\gamma$ touches zero without changing sign
-- is treated in detail and while we can deal with most of the cases
we are unable to produce a concise statement that says 
 which zeros of $\gamma_c(h)-\gamma$ are critical points (some are not!)
and, in general, what is the precise order of the transition. This is due to the fact that these transitions, unlike the 
denaturation transition, 
do depend on the details of $K(\cdot)$ and, to a certain extent, one needs  to do a case by case study.
%and the results do not have a universal character. 
Examples and more considerations on all these issues 
are developed at the end of the introduction and  in Section~\ref{sec:ex}.

\subsection{Outline of the approach, sharp estimates and limit path properties}
\label{sec:sharp}

As we already mentioned, the cornerstone is  \eqref{eq:Zc-2}.
In fact 
 \eqref{eq:Zc-2} reduces  sharp, respectively Laplace, estimates on $Z_{N, M,h}^c$
to sharp, respectively Laplace,  estimates on the renewal function $\bP\left( (N,M) \in \tilde \tau_h \right)$.
A quick overview of the behavior of  $\bP\left( (N,M) \in \tilde \tau_h \right)$ is

\medskip

\begin{enumerate}
\item If $h<0$,  so $\tilde \tau_h $ is terminating, we will show that there exists $C_h > 0$ such that
\begin{equation}
\label{eq:terminating}
\bP\left( (N,M) \in \tilde \tau_h \right) \stackrel{N, M\to \infty} \sim C_h \tilde K_h(N,M)\,.
\end{equation} 
%An upper bound of $ \bP\left( (N,M) \in \tilde \tau_h \right)$ is also given for every $(N,M)$ in Proposition~\ref{th:SharpZ}.
 \item 
For $h>0$, recall \eqref{eq:mu12} and the discussion right after that formula, one can show that
 \begin{equation}
\label{eq:DBM0-1}
 \lim_{\gep \searrow 0} \lim_{t \to \infty}
 \frac1t \log\bP\left( \left( vt + \{x\in \bbR^2: \, \vert x \vert \le \gep N\}\right) \cap  \tilde \tau_h  \neq \emptyset \right) \, = \,  
 -  D_h(v)\, ,
\end{equation}
where
$D_h(\cdot)$ is a non-negative function defined in $\bbR^2$, but equal to $+\infty$ outside  of  the first quadrant.
We shall see that
$D_h(\cdot)$ is linear along rays, that is $D_h(sv)=sD_h(v)$, and 
$D_h(v)=0$ if and only if $v \propto \mu_h$ (for us $v \propto (1,1)$). 
Moreover $D_h(v_1,v_2)=D_h(v_2,v_1)$ so there is no loss
of generality in  sticking to $D(1, \gamma)$, $\gamma\ge 1$. 
We will then see that $\gamma \mapsto D(1,\gamma)$ is affine for $\gamma$ larger than a critical value $\gamma_c(h)>1$
or smaller than another critical value that coincides with $1/\gamma_c(h)$ (this symmetry follows directly by the symmetry
$Z^c_{N, M, h}=Z^c_{M, N, h}$, that, in turn, is a consequence of $\mathrm{K}(n,m)=\mathrm{K}(m,n)$. On the other hand 
$\gamma \mapsto D(1,\gamma)$ is strictly convex in the interval $(1/ \gamma_c(h), \gamma_c(h))$ -- the \emph{Cram\'er region} -- and, 
when $\gamma $ is in this interval, sharp asymptotic estimates can be obtained. Namely, 
  we have that
\begin{equation}
\label{eq:DBM0}
 \bP\left( (N, M) \in \tilde \tau_h \right)\,  \sim\,   \frac{c_v}{\sqrt{t}}
 \exp\left(- t D_h(v)\right)\, ,
\end{equation}
where $t= \sqrt{N^2+M^2}$, $v=(N,M)/t$,
 $c_v$ is a positive constant (which depends of course also on $h$ and $K(\cdot)$) 
 and  the asymptotic statement is for $t\to \infty$ and it is uniform provided that
  the 
  unit vector  $v$ is in a compact arc of  circle  subset of  the open arc of the unit circle that goes form $(\gamma_c(h),1) /\sqrt{1+ (\gamma_c(h))^2}$ to 
$(1,\gamma_c(h)) /\sqrt{1+ (\gamma_c(h))^2}$.
\end{enumerate}

\medskip

We remark that
by putting \eqref{eq:Zc-2} and \eqref{eq:DBM0} together we readily see  that we can make
a substantial step ahead with respect to \eqref{eq:prelimZ}:
\begin{equation}
\label{eq:limZ}
\limtwo{N,M\to \infty:}{M\sim \gamma N} \frac 1 N \log Z^c_{N, M, h} \, =\, (1+\gamma) \tg - D_h( 1, \gamma)\, .
\end{equation}
From this formula, since $\tg$ is (implicitly) determined by \eqref{eq:hc} and since 
$D_h(\cdot)$ has a variational formulation, we will be able to use it to establish Theorem~\ref{th:F}
and Theorem~\ref{th:2}.

\medskip

But with what we just outlined we can go beyond Laplace type estimates:
 \eqref{eq:terminating} and \eqref{eq:DBM0}  yield sharp $N \to \infty$ estimates
on $Z_{N, M,h}^c$ for $h<0$, with $M\sim \gamma N$, any $\gamma>0$. Same  for $h>0$, but only for $\gamma$ in the Cram\'er region. We state here 
the result for the free case, which is less immediate than the constrained one:

\medskip

\begin{theorem} 
\label{th:Zsharp}
We have the following sharp estimates for $M \sim \gamma N$ and $\ga >1$:
\begin{enumerate}
\item For $h>0$ and $\gamma\in (1/\gamma_c(h), \gamma_c(h))$ there exists $c_{\gamma, h}>0$ such that 
\begin{equation}
Z_{N, M,h}^f \, \stackrel{N \to \infty}{\sim} \,
\frac{c_{\gamma, h}}{\sqrt{N}}\exp\left(N\tf_{\frac M N}(h)\right)\, .
\end{equation} 
\item For $h<0$, if $\overline{\ga} < (1+\ga)/2$ we have
\begin{equation}
Z_{N, M,h}^f \, \stackrel{N \to \infty}{\sim} \, \frac{K_f(N)K_f(M)}{1-\exp(h)} \,.
\end{equation}
Moreover, if $\overline{\ga} > (1+\ga)/2$ 
\begin{equation}
Z_{N, M,h}^f \, \stackrel{N \to \infty}{\sim} \, \frac{ \exp(h){\left( \sum_{n \ge 0} K_f(n) \right)}^2 }{(1-\exp(h))^2} \mathrm{K}(N,M) \,.
\end{equation} 
\end{enumerate}
\end{theorem}

\medskip 

Sharp estimates on the partition function  lead to sharp control on path properties:
\medskip

\begin{theorem} 
\label{th:pathsharp}
Choose $\gamma>0$ and consider the case $M\sim \gamma N$ and $\ga >1$.
\begin{enumerate}
\item Let $(\cF_1, \cF_2):=\max\{\tau \cap [0, N]\times[0,M]\}$ be the last renewal epoch in $[0,N] \times [0,M]$. 
For $h<0$ and $\overline{\ga} < (1+\ga)/2$, the law of $(\cF_1, \cF_2)$ under $\bP^f_{N, M, h}$ -- a probability measure on $(\{0\} \cup \bbN)^2$ -- 
converges for $N \to \infty$ to the probability distribution that assigns to $(i,j)$ probability 
\begin{equation}
\label{eq:cF}
( 1 - \exp(h) ) \bP \left( (i,j) \in \tilde \tau_h \right) \,.
\end{equation}
Set $\cL_1:=N-\cF_1$ and $\cL_2:= M- \cF_2$. For $h<0$ and $\overline{\ga} > (1+\ga)/2$, the law of $(\cL_1, \cL_2)$ under $\bP^f_{N, M, h}$ -- a probability measure on $(\{0\} \cup \bbN)^2$ -- 
converges for $N \to \infty$ to the probability distribution that assigns to $(i,j)$ probability
\begin{equation}
\label{eq:cF-2}
\frac{1}{{\left( \sum_{n \ge 0} K_f(n) \right)}^2} K_f(i) K_f(j) \,.
\end{equation}
Moreover, for $h<0$ and $\overline{\ga} > (1+\ga)/2$, we have
\begin{equation}
\label{eq:L}
\lim_{L\to \infty}\lim_{N\to \infty} \bP^f_{N, M, h}( \tau \cap [L, N-L]\times [L, M-L]= \emptyset) = 1 \,.
\end{equation}
\item For $h>0$ and $\gamma\in (1/\gamma_c(h), \gamma_c(h))$
we have that both $\tf_\gamma(h) - \gamma \partial_{\gamma} \tf_{\gamma}(h)$ and $\partial_{\gamma} \tf_{\gamma}(h)$
are positive and 
 the law of
$(\cL_1, \cL_2)$ under $\bP^f_{N, M, h}$ -- a probability measure on $(\{0\} \cup \bbN)^2$ -- 
converges for $N \to \infty$ to the probability distribution that assigns to $(i,j)$ probability
\begin{equation}
\label{eq:cgammah8}
\frac {1} {C_{\gamma,h}}
K_f(i) \exp \left(- i  \left( \tf_\gamma(h) - \gamma \partial_{\gamma} \tf_{\gamma}(h) \right) \right) K_f(j) \exp \left( -j \partial_{\gamma} \tf_{\gamma}(h) \right) \, 
\end{equation}
with $C_{\gamma,h}>0$ the normalization constant. Moreover the law of $\tau$ under $\bP^f_{N,M, h}$ -- a probability on 
the  subsets of $(\{0\} \cup \bbN)^2$ -- converges in the same limit to the law of 
a positive recurrent two-dimensional renewal with inter-arrival law given by the function from $\bbN^2$ to $[0, 1)$
\begin{equation}
\label{eq:cgammah9}
(i, j) \mapsto K(i+j)\exp\left( - i  \left( \tf_\gamma(h) - \gamma \partial_{\gamma} \tf_{\gamma}(h) \right)  -j \partial_{\gamma} \tf_{\gamma}(h) \right)\,.
\end{equation}
\end{enumerate}
\end{theorem}
\medskip

Theorem~\ref{th:pathsharp} can be summed up as:
\begin{enumerate}
\item In the delocalized phase, $h<0$ (and  $\overline{\ga}\neq (1+\ga)/2$, see below), there is no contact in the bulk of the system and, according to whether
the $K_f(\cdot)$ exponent $\overline{\ga}$ is larger or smaller than $(1+\ga)/2$ the two strands are free except for 
$O(1)$ contacts all close to the origin, or the two strands get detached after finitely many contacts (all close to the origin) 
and they meet again at a $O(1)$ distance from $(N, M)$, terminating with two free ends  of length $O(1)$.
In the  case $\overline{\ga}=(1+\ga)/2$ the slowly varying corrections $L(\cdot)$ and $\overline L (\cdot)$  matter and we leave out
this rather cumbersome analaysis. 
\item In the localized phase ($h>0$) and for $\gamma$ in Cram\'er region  the process converges to a persistent renewal that we determine: this is similar to what happens in the one-dimensional case, but in this new set-up the limit process
has the expression \eqref{eq:cgammah9} which is much less straightforward than the corresponding one-dimensional case. A number of other results can be proven, in the spirit of the one-dimensional analogs
(see \cite{cf:CGZ} and \cite[Ch.~2]{cf:GB}), but we have chosen to  limit ourselves to Theorem~\ref{th:pathsharp}(2) and we signal that its proof, see \S~\ref{sec:pathproofs}, is more informative than just \eqref{eq:L}.
\end{enumerate}

\subsection{Open issues and perspectives}
\label{sec:open}
We do not treat a number of natural issues: we list and discuss them here.

\subsubsection{The non Cram\'er regime}
For $h>0$ and $\gamma\not \in (1/ \gamma_c(h), \gamma_c(h))$ we do not give
 sharp estimates. To our knowledge sharp estimates on the renewal function
 in this regime  are for the moment not available (the most advanced reference available appears to be \cite{cf:BM2}).
The issue is not a secondary one: it is at the heart of understanding the transitions and the different phases that
one observes in the localized regime. And what one expects is rather clear: 
for $\gamma$ in the Cram\'er region we have seen that the free ends are microscopic, i.e. $O(1)$, and the limit process
is  just a recurrent renewal; 
for $\gamma$
 in the interior of the complementary of the Cram\'er region instead a big loop should appear (showing that it is only one 
 is one of the issues) or the free ends should become macroscopic. 
 The analysis is certainly different for the constrained case, in which the expected big loop can be anywhere along the chain,
 and for the free case in which this \emph{excess of bases} may be in the bulk (again, the location would be uniformly distributed over the length of the chain) or in the free ends and, like for the delocalized case, the exponents $\ga$
 and $\overline \ga$ will certainly play a crucial role. This and very similar issues are widely discussed in the
 physical literature (\cite{cf:KBM,cf:NG}) and the analogy with Bose condensation is regularly invoked, but the analysis is far from being rigorous.
 
\subsubsection{Counting the transitions in the localized regime}
We present examples with zero, one or two transitions. Can there be more than two? Are they always finitely many? 

\subsubsection{Sharp estimates at criticality} If $\sum_n n^3 K(n) < \infty$
sharp estimates for $h=0$ are covered by \eqref{eq:DBM0}. We have not treated this case here
because it would be natural to consider the complete spectrum of loop exponents, but we meet again with the 
limitations of multivariate renewal theory. The gPS model demands control only on the special class of renewals with
inter-arrivals $\mathrm{K}(n,m)=K(n+m)$ and we hope that a ad hoc treatment will lead to progress. And of course the gPS
model is one more motivation for a more systematic study of multivariate renewals.

\subsubsection{Disordered interactions}
Here the issues are several: we stick to the one of disorder relevance at criticality (see the review of the literature in \cite[Ch.~4]{cf:G}), but there
are most probably intriguing questions also away from criticality (in analogy with \cite{cf:GT-alea}).
The effect of a disorder on the critical point is directly related to obtaining sharp estimates on the renewal function
of the underlying renewal process, at least if the disorder is introduced via
an IID family $\{h_{n,m}\}_{(n,m)\in \bbN^2}$ of random variables. In fact the tools 
developed for the basic disordered PS model (\cite{cf:G} and references therein) can be applied, but the problem is that 
sharp estimates are available only for the very particular case  of $\sum_n n^3 K(n) < \infty$.
And there is the issue that such an IID disorder is not the most suited for DNA modeling, but
if the aim is understanding the effect of noise on critical behaviors  this way of introducing 
the disorder is     certainly acceptable (and it is what has been done also in the bio-physical literature,
even sticking to the DNA/RNA set-up! See for example \cite{cf:BundHwa}). A more natural disorder is however 
obtained by assigning to each strand a sequence of, possibly IID, potentials $\{h^{(1)}_j\}_{j \in \bbN}$ and  
$\{h^{(2)}_j\}_{j \in \bbN}$ -- one can imagine the case in which the  two sequences are independent or 
the case in which they are (strongly) correlated -- and $h_{n,m}=h^{(1)}_n h^{(2)}_m$: note that correlations are introduced
with this product choice even if the two sequences are independent. This appears to be a very challenging model (see \cite{cf:BePo} and references therein for the issues that arise when correlations are introduced in the disorder sequence for the basic PS model).

\subsubsection{Related models}
 The gPS model is intimately related to the more complex RNA models for secondary structure: \cite{cf:EON}, where the vast literature is cited, is particularly interesting for us because  RNA models are linked with the gPS model. Models for circular DNA 
 \cite{cf:Yer} are also very much related to gPS, as pointed out for example in \cite{cf:GO0}. In \cite{cf:KBM} the authors focus
 on an issue (existence of one macroscopic loop) for circular DNA that is precisely the one that we face outside of the Cram\'er regime.  Finally, the gPS model can be seen as a toy model for interacting self-avoiding walks. In this direction we signal the 
 Zwanzig-Lauritzen model that has been tackled first by generating functions techniques (e.g. \cite{cf:BOPG}) and recently
 by probabilistic methods in \cite{cf:CPN,cf:NP}.
 
\subsection{Organization of the paper}
 In Section~\ref{sec:LDLL}, we present the results  on Large and Sharp Deviations  for bivariate renewal processes that we use. In Section~\ref{sec:constrained}, we introduce the constrained model and study the free energy in the localized and delocalized regime proving theorem~\ref{th:F}, Lemma~\ref{th:gammac} and Theorem~\ref{th:2}. At the end of this section, we work out explicitly some examples. Finally in Section~\ref{sec:free}, we prove Proposition~\ref{th:propo1} and we compute the sharp estimates of the free partition function and the path properties proving Theorem~\ref{th:Zsharp} and Theorem~\ref{th:pathsharp}. For the Dominated Converge Theorem we use the shortcut (DOM).

\section{Large Deviations and Local Limit Theorems for bivariate renewals}
\label{sec:LDLL}

We now give a number of results on the renewal $\tilde \tau_h$, $h>0$, defined at the beginning of Section~\ref{sec:fe2}. As it will be clear, they follow directly from various results that one can find in
 \cite{cf:BM,cf:BM1} where a more general case is treated (starting from the fact that we limit ourselves  to the  two-dimensional case). 
So let us start by introducing  the exponential moment generating function of $\xi:= (\tilde \tau_h)_1$ (we recall:
$\bP(\xi=(n,m))= \tilde K_h(n,m)$)
\begin{equation}
M (\lambda)\, := \, \bE \left[  \exp \left( \langle \lambda,\xi \rangle \right)\right] \, , 
\end{equation}
where $\gl\in \bbR^2$, $\langle \lambda, \xi\rangle =  \lambda_1 \xi_1 +   \lambda_2 \xi_2 $. 
From the definition of $ \tilde K_h(\cdot, \cdot)$ one readily sees that  $M(\gl)< \infty$ if and only if  both $\gl_1\le \tg$ and $\gl_2\le \tg$.
The Large Deviation  function $\gL(\cdot)$ corresponding to the random vector $\xi$ is  the Legendre transform of the function $ \log M (\cdot)$
\begin{equation}
\label{eq:Lambda}
\Lambda(\theta) \, =\,  \sup_{\lambda} \left\lbrace  \langle \lambda,\theta \rangle - \log  M (\lambda)   \right\rbrace\,,  
\end{equation}
where  $\theta \in\bbR^d$. Since we are after the renewal function of $\tilde \tau_h$ the Large Deviation function
of the inter-arrival random vector
 is just an intermediate step.
The asymptotic behavior of the renewal function is directly related to 
 the so-called {\sl second deviation function}, introduced and investigated  in \cite{cf:BM}: 
 \begin{equation}
\label{eq:inf}
D_h(\theta) \, =\, D(\theta) \, =\,  \inf_{s > 0} \frac{\Lambda(s \theta)}{s}\, ,
\end{equation}
where $\theta \in \bbR^d$ 
 and the notation with the subscript $h$ will be useful further on to remind the dependence on the parameter
 but at this stage it is rather superfluous. 
 
 \medskip
 
 \begin{rem}
 \label{rem:Dprop}
 From \eqref{eq:inf} one can see that $D(s\theta)=s D(\theta)$ for every $s\ge 0$.
In \cite[pp. 652-653]{cf:BM} a detailed analysis of $D(\cdot)$ is given, notably the fact that 
it is convex: for $p \ge 0, q \ge 0, p+q=1$ and $\theta , \eta \in \bbR^d$, 
\begin{equation}
D( p \theta + q \eta ) \le p D(\theta) + q D(\eta) \,.
\end{equation}  
We can immediately deduce from these properties that for every $\theta , \eta \in \bbR^d$, we have 
\begin{equation}
D(\theta + \eta) \le D(\theta) + D(\eta) \,.
\end{equation}
\end{rem}

\medskip

Here is an important step:

\medskip

\begin{proposition}
\label{th:D}
\cite[Theorem 1]{cf:BM}.
For every $\theta=(\theta_1, \theta_2) \in \bbR^2$
\begin{equation}
\label{eq:D}
 D(\theta) \, =\,  \sup_{\lambda \in A} \langle \lambda , \theta \rangle \,=\,   \sup_{\lambda \in \partial A} \langle \lambda , \theta \rangle \, ,
\end{equation} 
where 
$A$ is the closed convex set $\left\lbrace \lambda \in \bbR^2 : \, M(\lambda) \le 1 \right\rbrace$ and $ \partial A$ is the boundary of A. 
\end{proposition}

\medskip

%From \eqref{eq:D} one directly reads the convexity of $D(\cdot)$ (cf. Rem	ark~\ref{rem:Dprop}). 

 It is now practical to focus on the specific case we are considering, notably the fact that
$D(\theta)=\infty$ if $\theta$ is not in the first quadrant is an intuitive consequence of the fact that our process
has increments that have positive components and can be read out of the structure of $A$. 
 Let us make $A$
more explicit
\begin{equation}
\label{eq:Ah}
A\,=\, \left\{ \gl \in \bbR^2:\, q_h\left( \gl_1, \gl_2\right)\,\le\, 1 
\right\}\, ,
\end{equation}
where $q_h(\cdot)$ is defined in \eqref{eq:q_h}
and we recall that $h>0$ and 
$\tg=\tg(h)>0$ is chosen so that $(0,0)\in \partial A$. Note that $q_h\left( \gl_1, \gl_2\right)=q_h\left( \gl_2, \gl_1\right)$
and that $q_h(\cdot)$ is convex (this of course implies the convexity of  $A$) 
 and it is
symmetric with respect to the diagonal of the first and third quadrant.
\medskip

\begin{lemma}
\label{th:Ah}
We have
\begin{equation}
\label{eq:Ah.0}
A\, \subset \, \left\{ \gl \in \bbR^2:\, \gl_1 \le \tg \text{ and } \gl_2 \le \tg\right\}
\cap \left\{ \gl \in \bbR^2:\, \gl_2 \le -\gl_1\right \}\,.
\end{equation}
Moreover $\overline{\gl}_1\stackrel{\eqref{eq:q_h}}= \sup\{\gl_1<0: \, q_h(\gl_1, \tg)=1\}<-1$ and the equation $q_h(\gl_1, \gl_2)=1$
is uniquely solvable for $(\gl_1, \gl_2) \in [\overline{\gl}_1, \tg]^2$, defining the curve $\cW_h$, symmetric with respect to the diagonal of the first and third quadrant. $\cW_h$ is the graph
of a concave and  decreasing function $\tilde \gl_2: [\overline{\gl}_1, \tg] \to [\overline{\gl}_1, \tg]$  which satisfies $\tilde \gl_2(\overline{\gl}_1)=\tg$,
$\tilde \gl_2(\tg)=\overline{\gl}_1$ and $\tilde \gl_2(0)=0$. Moreover $\tilde \gl_2$ is analytic in the interior of its domain and 
\begin{equation}
\partial A \, =\, \left\{(\gl_1, \tg):\, \gl_1< \overline{\gl}_1\right\} \cup \left\{( \tg, \gl_2):\, \gl_2< \overline{\gl}_1\right\} \cup \cW_h \, .
\end{equation}
Finally, $q_h(\gl_1, \gl_2)<1$ for $(\gl_1, \gl_2)\in \partial A \setminus \cW_h$.
\end{lemma}

\medskip

\noindent
{\it Proof.} For this proof is practical to keep at hand part ({\sc left}) of Figure~\ref{fig:Bh}. %First of all the convexity of $q_h(\cdot)$ implies the convexity of $A$. 
The fact that $A\, \subset \, \left\{ \gl \in \bbR^2:\, \gl_1 \le \tg \text{ and } \gl_2 \le \tg\right\}$
is just the fact that $q_h(\gl_1,\gl_2)=\infty$ if $\gl_1 \wedge \gl_2 > \tg$. On the other hand this last observation, coupled with  $q_h(\gl_1, \gl_2)=
q_h(\gl_2, \gl_1)$, $q_h(0,0)=1$ and convexity of $A$, tells us that $A$ does not go above the line $\gl_2=-\gl_1$ and 
\eqref{eq:Ah.0} is established. Now, since $q_h(\cdot,\cdot)$ is separately non-creasing function (and even increasing where it is bounded),
$\partial A$ contains $\{\gl: \, q_h(\gl)=1\}$. So the issue is the solvability of $q_h(\gl)=1$ and it is straightforward to see that
 $q_h(\gl)=)$ has a (unique, by monotonicity) solution if $\gl_1 \in [\overline{\gl}_1, \tg]$ and this way we define a function
 $\tilde \gl_2(\cdot):  [\overline{\gl}_1, \tg] \mapsto [\overline{\gl}_1, \tg]$ (note that at this stage it is already clear that $\overline{\gl}_1\le -1$)
 which is analytic in the interior of its domain, by the analytic Implicit Function Theorem. It is actually immediate to check 
 that this function is not linear (for example, compute the second derivative at the origin), so  $\overline{\gl}_1<-1$. 
 Finally $q_h(\gl_1, \tg)<1$ for $\gl_1< \overline{\gl}_1$ and, by symmetry, $q_h( \tg, \gl_2)<1$ for $\gl_2< \overline{\gl}_1$.
 The proof is therefore complete.
\qed

\medskip

\begin{rem}
\label{rem:hdep}
Lemma~\ref{th:Ah} is given for fixed $h>0$, but  $\overline{\gl}_1$
depends also on $h$ and when we need to make this dependence explicit we write
$\overline{\gl}_1 (h)$. The same is true for the function $\tilde \gl_2(\cdot)$ and we write $\tilde \gl_{2, h}(\cdot)$. 
Note that the analyticity of  $\overline{\gl}_1(\cdot)$ on the positive semi-axis is a direct consequence 
of  the Analytic Implicit Function Theorem (\cite[Sec.~2.3]{cf:Primer}, and it is just a matter of analyticity in one variable).
If instead we consider the function $(\gl_1, h)\mapsto \tilde \gl_{2, h}(\gl_1)$, with $\tilde\gl_{2, h}(\gl_1)$ which is 
 obtained by 
solving for $\gl_2$ the equation $q_h(\gl_1, \gl_2)=1$ and we have seen that this requires $\gl_1 \in [  \overline{\gl}_1 (h), 
\tg(h)]$. Therefore, by the Analytic Implicit Function Theorem \cite[Sec.~2.3]{cf:Primer}, 
the function $(\gl_1, h)\mapsto \tilde \gl_{2, h}(\gl_1)$ is analytic (in two variables this time)
in the domain  $h>0$ and $\gl_1 \in (  \overline{\gl}_1 (h), 
\tg(h))$. 
\end{rem}

\medskip

From Proposition~\ref{th:D} and Lemma~\ref{th:Ah} we can derive a number of consequences, like for example
the fact that
 $D(\theta)< \infty$ for every
$\theta$ in the first quadrant (that here includes the two axes) and that, when $D(\theta)<\infty$, that is in the first quadrant, 
there are only two possibilities: either the supremum in the rightmost term in \eqref{eq:D}
is reached in the interior of $\cW_h$ or {\sl at the boundary}, that is  in $\{( \overline{\gl}_1, \tg), ( \tg,\overline{\gl}_1)\}$. We observe by direct inspection that 
if $\theta_2>\theta_1$ then if the supremum is not achieved in the interior, then it is achieved at  $( \overline{\gl}_1, \tg)$.
Moreover if $\theta_1= \theta_2>0$ the supremum is always achieved in the interior and, more precisely, at $(0,0)$ and therefore 
$D(\theta_1, \theta_1)=0$. This induces a partition of the first quadrant:
$\theta \in E_h$ if the supremum is achieved in the interior
and $\theta \in E_h^\complement$ if it is achieved at the boundary. 
 It is also useful to remark that $E_h$ is an open sector (\cite[p.~653]{cf:BM}, or it can be seen directly in our specific set-up): in our case it is also
symmetric with respect to the diagonal of the first quadrant, that is there exists $\gp \in (0, \pi/4)$ such that 
$E_h=\{ \theta \in [0, \infty)^2: \, \theta_2/ \theta_1\in (\tan (-\gp+\pi/4) , \tan (\gp+\pi/4) )\}$.

\medskip

\begin{rem}
\label{rem:Cramer}
\rm
Moreover,
\cite[Theorem 2]{cf:BM} tells us
 that if $\theta \in E_h$ then 
 the
infimum in  \eqref{eq:inf} is attained at a unique point $s(\theta)$, so $D(\theta)=\gL(\theta s(\theta))/ s(\theta)$,
and $\gL(\cdot)$ is analytic and strictly convex at $s(\theta)$.  
 This is what we may call the {\sl Cram\'er region} of parameters: such a region being the 
set of $\theta$'s for which the large deviations trajectories can be made typical by a suitable change of measure
({\sl tilting}).
\end{rem} 
\medskip

\medskip

We are now ready to state the  result that links $D(\cdot)$ and the renewal function.
Once again it is obtained by restricting to our context a result -- Theorem 5 -- in \cite{cf:BM}, see however Remark~\ref{rem:D} below.
We employ the notation $\gL''(\theta)$ for the Hessian matrix of $\gL$. 

\medskip

\begin{proposition}
\label{th:DBM}
The following representation  holds:
\begin{equation}
\label{eq:DBM}
 \bP\left( (N, M) \in \tilde \tau_h \right) \, = \,  
 \frac 1{t^{1/2}} \left( 
 \sqrt{
 \frac{\mathrm{det}\left(\gL''(s(v)v)
 \right)}
 {2\pi s(v)
 \left\langle v, \gL''(s(v)v) v \right\rangle}
 }
 + \gep(N,M)
 \right)
 \exp\left(- t D(v)\right)\, ,
\end{equation}
where $t=\sqrt{N^2+M^2}$, $v=(v_1,v_2)=(N,M)/t \, \in E_h\cap\{\theta:\, \vert \theta\vert =1\}$  and for every compact $ C \subseteq  E_h\cap\{\theta:\, \vert \theta\vert =1\}$
\begin{equation}
\lim_{t \to \infty}\suptwo{N, M: \sqrt{N^2+M^2}=t}{ \ \ v \in C}\gep(N,M)\, =\, 0\, .
\end{equation}
\end{proposition}

\medskip

Recalling \eqref{eq:Zc-2}, we see that Proposition~\ref{th:DBM} implies the sharp estimate
\begin{equation}
\label{eq:Zc-sharp}
Z_{N, M, h}^c \, =\, 
\frac{A(M/N)+ \tilde \gep (N, M)}{\sqrt{N}}
\exp \left((N+M) \tg(h)- N D_h(1,M/N)
\right) \, , 
\end{equation}
where $A(M/N)$ is equal to $(1+(M/N)^2)^{1/2}$ times the square root term in \eqref{eq:DBM}
(which depends only on $v$, which in turn is just a function of $M/N$) and 
$\tilde \gep (N, M)= (1+(M/N)^2)^{1/2}\gep (N, M)$.

\medskip

\begin{rem}
\label{rem:D}
\rm
In \cite{cf:BM} the factor $s(v)$ that appears just after $2\pi$ in \eqref{eq:DBM} has  the exponent $d+3=5$. This
formula appears also in \cite[p.~11]{cf:BM1}, with an additional oversight. The formula we give is in agreement with
\cite{cf:Doney} who covers only the case $v \propto \mu_h$ (cf. \eqref{eq:mu12}), but in greater generality than 
\cite{cf:BM,cf:BM1}.  Formula \eqref{eq:DBM} a priori requires some exponential decay of the inter-arrival law to make sure that
the Hessian is computed at an analyticity point of $\gL(\cdot)$. Actually, the analysis in \cite{cf:Doney} shows that this is not necessary
and \eqref{eq:DBM} still holds true for $h=0$ if $\sum_{n,m \ge 1} (n+m)^2 K(n+m)= \sum_{t\ge 2} t^2 (t-1) K(t)< \infty$. 
With the help of the notation in \cite{cf:Doney} we  remark that \eqref{eq:DBM} can be made slightly more readable 
if we observe that, as it is well known, the inverse $B(\theta)$ of $\Lambda''(\theta)$ is the covariance matrix of the {\sl tilted} random 
vector $X=(X^1, X^2)$ with $\bP(X=(n,m))\propto \tilde K_h(n,m)\exp(n \gl_1(\theta) +m \gl_2 (\theta)) / C_\theta$ with $C_\theta =  \sum_{n,m} K_h(n,m)\exp(n \gl_1(\theta) +m \gl_2 (\theta))$ and $\gl(\theta)$ is the optimal point of \eqref{eq:Lambda}. Therefore with $v=\theta/\vert \theta\vert$
\begin{multline}
\frac{\mathrm{det}\left(\gL''(\theta)
 \right)}
 { \left\langle v, \gL''(\theta) v \right\rangle}\, =\, \frac1{\mathrm{det}\left(B(\theta)\right)\left\langle v, (B(\theta))^{-1} v \right\rangle}
 \\ =\,  \left\langle v, \left(\begin{array}{cc} \bE [(X^2)^2]- \bE [X^2]^2 & -\bE[X^1X^2]+  \bE[X^1]\bE[X^2]
 \\ 
  -\bE[X^1X^2]+  \bE[X^1]\bE[X^2]
  & \bE [(X^1)^2]- \bE [X^1]^2 \\ \end{array}\right)
 v 
 \right\rangle^{-1}\,.
\end{multline}

\end{rem}
\medskip

A Local Limit Theorem, analog to Proposition~\ref{th:DBM}, for $ v \in E_h^\complement$ is 
available at the moment  (see \cite[Theorem 2.1]{cf:BM2}) only if the entries of $\gL^{''}(v)$
are all finite, and this is not always the case in our set-up, notably it is not the case if the exponent $\ga$ entering the definition
of $K(\cdot)$ is smaller than three. Nevertheless, the following weaker result will suffice for our purposes: 

\medskip

\begin{proposition}
\label{th:BM1}
For every $\theta$ 
\begin{equation}
\label{eq:BM1.0}
\limsup_{t \to \infty}
\frac 1t \log \bP\left( \lceil t \theta \rceil  \in \tilde \tau_h \right) \, \le\,  -  D\left(\theta \right)\, ,
\end{equation}
and if $\theta \in E_h^\complement$, with $\theta_1>0$ and $\theta_2>0$, then 
\begin{equation}
\label{eq:BM1.1}
\lim_{t \to \infty}
\frac 1t \log \bP\left( \lceil t \theta \rceil  \in \tilde \tau_h \right) \, =\,  -  D\left(\theta \right)\, .
\end{equation}
\end{proposition}
\medskip

Of course \eqref{eq:BM1.1} holds also for $\theta \in E_h$ as a immediate consequence of Proposition~\ref{th:DBM}.
\medskip

\noindent
{\it Proof.}
The upper bound is a direct consequence of the Large Deviations Principle  \cite[Theorem 4]{cf:BM}): 
\begin{equation}
\lim_{\gep \searrow 0} \lim_{t \to \infty} \frac{1}{t} \log 
\bP\left( t \{v\in \bbR^2:\, \vert v-\theta\vert \le \gep\} \cap \tilde \tau_h \neq \emptyset \right) \, = \,
-  D\left(\theta\right)\, .
\end{equation}

For the lower bound we assume without loss of generality that $\theta_2>\theta_1$ and 
we observe that if $\theta \in E_h^\complement$ -- a closed set (recall Remark~\ref{rem:Cramer}
and the explanation that precedes it) -- then either $\theta$ is in the boundary or in the interior of 
$E_h^\complement$. If it is in the interior   then there exists $\theta_2^\star< \theta_2 $ with $(\theta_1, \theta_2^\star)$ in the boundary of $  E_h^\complement$, so that
for every $\gep>0$ small we have $(\theta_1,\theta_2^\star-\gep) \in E_h$. If $\theta$ is in the boundary  of 
$E_h^\complement$ we directly set $\theta_2^\star< \theta_2 $. 
By Proposition~\ref{th:DBM}
\begin{equation}
\lim_{t\to \infty} \frac 1t \log \bP
\left( \left( \left\lceil t\theta_1\right\rceil -1,  
 \left\lceil t ( \theta_2^\star-\gep) \right\rceil \right) \in \tilde \tau _h \right)\, =\, - D_h \left( \left( \theta_1, \theta_2^\star-\gep
 \right) \right)\, .
\end{equation}
But 
\begin{equation}
\bP
\left(  \left\lceil t \theta \right \rceil \in \tilde \tau _h\right) \, \ge \, 
\bP
\left( \left( \left\lceil t\theta_1\right\rceil -1,  
 \left\lceil t ( \theta_2^\star-\gep) \right\rceil \right) \in \tilde \tau _h \right)
 \, \tilde K_h \left( 1,  \left\lceil t\theta_2 \right \rceil - 
 \left\lceil t(\theta_2^\star -\gep ) \right \rceil \right) \, ,
\end{equation}
and therefore we have 
\begin{equation}
\liminf_{t\to \infty} \frac 1t \log \bP
\left(  \left\lceil t \theta \right \rceil \in \tilde \tau _h\right) 
\, \ge \,  - D_h \left( \left( \theta_1, \theta_2^\star-\gep
 \right) \right) 
 - \tg(h) \left( \theta_2-\theta_2^\star + \gep\right)\,,
\end{equation}
and, by continuity of $D_h(\cdot)$, we get to
\begin{equation}
\liminf_{t\to \infty} \frac 1t \log \bP
\left(  \left\lceil t \theta \right \rceil \in \tilde \tau _h\right) 
\, \ge \,  - D_h \left( \left( \theta_1, \theta_2^\star
 \right) \right) 
 - \tg(h) \left( \theta_2-\theta_2^\star \right)\,.
\end{equation}
Since $\left( \theta_1, \theta_2^\star
 \right)\in E_h^\complement$ we have that 
 \begin{equation}
 D_h \left( \left( \theta_1, \theta_2^\star
 \right) \right) 
 + \tg(h) \left( \theta_2-\theta_2^\star \right) \, =\, 
 \overline \gl_1  \theta_1 + \tg(h)  \theta^\star + 
 \tg(h) \left( \theta_2-\theta_2^\star \right)\, =\, 
  \overline \gl_1  \theta_1 + 
 \tg(h) \theta_2 \, ,
 \end{equation}
and since $\theta\in E_h^\complement$ the rightmost term is $D_h(\theta)$ and we are done.
\qed

\section{The constrained model}

\label{sec:constrained}

\subsection{The free energy of the constrained model}
\label{sec:constr1}
Let us start by observing that
\begin{equation}
\label{eq:tflim0}
\tg(h)\, =\, \lim_{N \to \infty} \frac 1{2N}  \log Z^c_{N, N, h}\, =\, \frac 12 \tilde \tf_1(h)\,  .
\end{equation}
This can be seen, for $h\le 0$,  by the elementary argument presented  just after \eqref{eq:defhc}
and, for $h>0$,  by \eqref{eq:Zc-2} and Proposition~\ref{th:DBM}: in fact
$D(1,1)=0$.
\medskip

But in general $M\neq  N$ and, without loss of generality
we can assume $M\ge N$, hence $\gamma\ge 1$. The first result in this direction, which includes \eqref{eq:tflim0}, is:
\medskip

\begin{proposition}
\label{th:Fgc}
For every $\gamma\ge 1$ 
\begin{equation}
\label{eq:Fgc}
\tilde \tf_\gamma (h)\, :=\, \limtwo{N, M \to \infty:}{\frac MN \to \gamma}
\frac 1N \log Z^c_{N, M, h}\, =\,
\begin{cases}
0 & \text{ if } h\le 0\, ,\\
 (1+\gamma) \tg(h) - D_h(1,\gamma) & \text{ if } h>0\, ,
\end{cases} 
\end{equation}
and in fact we have also that for every $L>1$
\begin{equation}
\label{eq:Fgc-added}
\lim_{\gep\searrow 0} \lim_{N\to \infty }\sup_{\gamma\in [1, L]}\sup_{ M:\, \vert (M/N)-\gamma\vert \le \gep} 
\left \vert  \frac 1N \log Z^c_{N, M, h} - \tilde \tf_\gamma (h) \right \vert \, =\, 0\, .
\end{equation}
Moreover we have
\begin{equation}
\label{eq:FgcDh}
D_h(1,\gamma)\,=\, \max_{\gl \in B_h} \left(\gl_1+ \gamma \gl_2\right) \, ,
\end{equation}
with 
\begin{equation}
\label{eq:Bh}
B_h\, =\, \left \{ \gl :\, \overline{\gl}_1 \le \gl_1\le 0, \, 0 \le \gl_2 \le \tg \, , \sum_{n,m} \tilde K_h(n,m) \exp(\gl_1n +\gl_2 m) \, =1 
\right\}\, .
\end{equation}
Finally
$\tilde \tf_\gamma (h)>0$ for $h>0$ and in fact
\begin{equation}
\label{eq:febounds}
2\tg(h)\,  \le\,  \tilde \tf_\gamma (h)\,  \le\,  (1+\gamma) \tg(h)\, .
\end{equation}

\end{proposition}

\medskip

\noindent
{\it Proof.} It is slightly more practical to establish first a result which is a priori weaker than  \eqref{eq:Fgc}.
Consider first the limit for $N\to \infty$ and $M= \lfloor \gamma N \rfloor$ and let us establish the rightmost equality in \eqref{eq:Fgc} with this notion of limit:  call $\breve \tf_\gamma (h)$ this expression.
Again,  the case $h\le 0$ of \eqref{eq:Fgc} is treated by elementary methods just after the statement of Proposition~\ref{th:propo1}. For the case
$h>0$ we observe that
by \eqref{eq:Zc-2} it suffices to show that
%we have that \eqref{eq:Fgc} reduces to showing that 
\begin{equation}
\label{eq:DBMcor}
\lim_{N \to \infty}
\frac 1N \log \bP \left( (N, \lfloor \gamma N \rfloor) \in \tilde \tau_h \right)\, =\, -D_h(1, \gamma)\,.
\end{equation}
But \eqref{eq:DBMcor} is a direct consequence of Proposition~\ref{th:DBM} and Proposition~\ref{th:BM1}, so
the weaker version of   \eqref{eq:Fgc}
is established. Now we step to \eqref{eq:Fgc-added}: if we  establish \eqref{eq:Fgc-added} with $\tilde \tf_\gamma (h)$
replaced by $\breve \tf_\gamma (h)$, then 
  the original version of
 \eqref{eq:Fgc} holds, which implies that $\breve \tf_\gamma (h)=
 \tilde \tf_\gamma (h)$ and, in turn, that also  the original version of \eqref{eq:Fgc-added} holds. It is therefore a matter of comparing (uniformly) the limits along all sequences
 of $(M,N)$, with $\vert (M/N)-\gamma\vert \le \gep$, with the special case $M= \lfloor \gamma N \rfloor$.
 But this is obtained for example by exploiting 
 that if $M'> M$  we have that for any $C_1>1+\ga$ there exists  $C_2>0$ such that 
 \begin{equation}
 \label{eq:surg1}
 Z_{N, M, h}^c \, \le \, C_2 {M'}^{C_1} Z_{N, M', h}^c\, .
 \end{equation}  
 \eqref{eq:surg1} can be established by using $K(n)/K(n+m)\le 1/K(n+m)\le  C_2 m^{C_1}$  for $m\ge n$.
  More precisely we apply this inequality to the term
 $K(l_n+t_n)$ in \eqref{eq:Zc} to stretch the last renewal so that it matches the boundary constraint $(N,M')$, and then we 
 allow $n$ to go up to $N\wedge M'$, which in this case is $N$ anyways, and the new constant on $\underline{t}$ is
  $\vert \underline{t} \vert=M'$.
 Inequality \eqref{eq:surg1} can then be used to sandwich the partition functions 
 $Z_{N, M, h}^c$, with $\vert (M/N)-\gamma\vert \le \gep$ and $N$ sufficiently large, between
 $ Z_{N,\lfloor (\gamma -2\gep)N\rfloor, h}^c/N^{2C_1}$ and $Z_{N,\lfloor (\gamma +2\gep)N\rfloor, h}^c N^{2C_1}$.
 The continuity of $\gamma \mapsto \breve \tf_\gamma (h)$, and hence the uniform continuity and boundedness on compact sets -- in our case $[1, L]$ -- completes the argument and  the proof of \eqref{eq:Fgc} and \eqref{eq:Fgc-added}.

Let us check \eqref{eq:FgcDh}: it is of course a matter of replacing $\partial A$ in \eqref{eq:D}, Proposition~\ref{th:D},
by $B_h$. Much of the work has been done in Lemma~\ref{th:Ah}: we are just left with showing 
that we can restrict the supremum to $B_h$. First of all
the symmetry of $\partial A$ tells us that  $\gamma \ge 1$ implies that $\sup_{\partial A} (\gl_1+\gamma \gl_2)$ 
does not change if we restrict $\partial A$ to $\gl_1\le 0$. 
Morevover, by  Lemma~\ref{th:Ah}, if $\gl_1 < \overline{\gl}_1$, 
then $(\gl_1, \tg)\in \partial A$, but  since $\gl_1 + \gamma \tg <  \overline{\gl}_1+ \gamma \tg$ for $\gl_1 < \overline{\gl}_1$
we can actually neglect these points in taking the supremum.

We are left with the positivity of $\tilde \tf_\gamma (h)$ for $h>0$. This  follows directly by observing that
\begin{equation}
Z_{N, M,h }^c \, \ge\,  Z_{N-1, N-1,h }^c   \exp(h) K(M-N+2)\, ,
\end{equation}
which implies more than the positivity, that is  $\tilde \tf_\gamma(h) \ge 2 \tg(h)$. Finally, by exploiting also that $D_h(\cdot)\ge 0$ (for $h>0$) and 
that in any case $\tilde \tf_\gamma (h) =0$ for $h \le 0$, we see that \eqref{eq:febounds} holds and
the proof of Proposition~\ref{th:Fgc} is therefore complete.
\qed

\medskip

We can go beyond Proposition~\ref{th:Fgc} by exploiting the variational problem \eqref{eq:FgcDh}-\eqref{eq:Bh}.
Note that since $\tilde \tf_1(h)=2 \tg(h)$ we know that $\tilde \tf_1(\cdot)$ is analytic except at the origin, but 
for $\gamma>0$ the situation is more involved.

\medskip

\begin{proposition}
\label{th:reg} 
The function $\gamma_c: (0, \infty) \to (1, \infty)$, defined in
 \eqref{eq:gammah}, is real analytic.  Moreover
$\tilde \tf_\gamma(\cdot)$ is analytic on the positive semi-axis out of the set $\{h:\, 
\gamma_c(h)-\gamma=0\}$, but  $\tilde \tf_\gamma(\cdot)$ is not analytic at 
the values $h$ at which $\gamma_c(h)-\gamma$ changes sign. 
\end{proposition}

\medskip

\begin{rem}
\label{rem:wait}
Of course the regularity issue is not completely resolved by Proposition~\ref{th:reg}, both because 
it does not make clear wether or not every all the points in the discrete set $\{h:\, 
\gamma_c(h)-\gamma=0\}$ are non-analyticity points and because it does not specify the 
type of singularities. 
\end{rem}

\begin{figure}[h]
\begin{center}
\leavevmode
\epsfxsize =13.5 cm
\psfragscanon
\psfrag{0}[c][l]{\small $0$}
\psfrag{gl}[c][l]{\small $\overline{\gl}_1$}
\psfrag{G}[c][l]{\small $\tg$}
\psfrag{l1}[c][l]{\small $\gl_1$}
\psfrag{l2}[c][l]{\small $\gl_2$}
\psfrag{l2l}[c][l]{\small $\tilde\gl_2(\gl_1)$}
\psfrag{A}[c][l]{({\sc left})}
\psfrag{B}[c][l]{({\sc right})}
\epsfbox{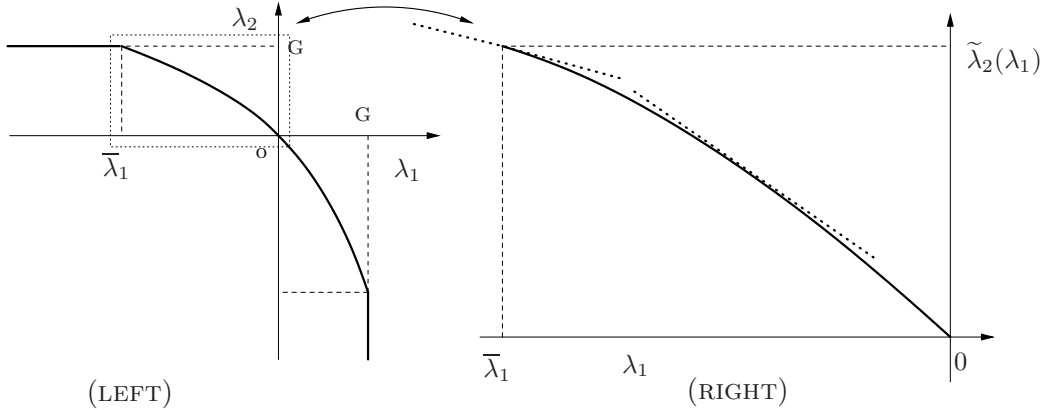}
\end{center}
\caption{\label{fig:Bh} In part ({\sc left}) $\partial A$, with $A=\{\gl:\, q_h(\gl_1, \gl_2)=1\}$, is represented by the
thick line. Notice the symmetries and the convexity of $A$. In part ({\sc right}) we zoom into the relevant part for our analysis.}
\end{figure}

\noindent
{\it Proof.} Recall the function $q_h(\cdot)$ and $\overline{\gl}_1$ from \eqref{eq:q_h} and \eqref{eq:gl1bar}. Recall also 
the definition  of  
$\gl_1 \mapsto \tilde \gl_2(\gl_1)$ in Lemma~\ref{th:Ah}:  $\tilde \gl_2(\gl_1)$ the only solution $\gl_2$ to $q_h(\gl_1, \gl_2)=1$, for $\gl_1\in [\overline{\gl}_1, \tg]$ (see Figure~\ref{fig:Bh}). Recall also that
 $\tilde \gl_2(\cdot)$ is concave -- i.e. convex down -- and analytic in $(\overline{\gl}_1, 0)$. We have that
$\tilde \gl_2'(\gl_1)$ equals $-\partial_{\gl_1}q_h(\gl_1, \gl_2)/\partial_{\gl_2}q_h(\gl_1, \gl_2)$ evaluated at 
 $\gl_2=\tilde \gl_2(\gl_1)$. Therefore
 \begin{equation}
 \label{eq:fgammah}
 \lim_{\gl_1 \searrow \overline{\gl}_1 } \tilde \gl_2'\left(\gl_1 \right)\, =\, 
 \tilde \gl_2'\left(\overline{\gl}_1^+\right)\, =\, - \frac{\sum_{n,m} n K(n+m)\exp\left(h- n \left(\tg(h) + \left\vert \overline{\gl}_1\right\vert \right)\right)}{\sum_{n,m} m K(n+m)\exp\left(h- n \left(\tg(h) + \left\vert \overline{\gl}_1\right\vert \right)\right)}\, \in (-1,0)\, .
\end{equation}
The denominator is bounded because of \eqref{eq:c>0} and 
the fact that the ratio is bounded below by $-1$ is just a consequence of
concavity and the fact that $\tilde \gl_2'(0)=-1$.   Moreover 
$ \tilde \gl_2'\left(\overline{\gl}_1^+\right)$ is a function of $h$: with the notations
in 
Remark~\ref{rem:hdep} we write rather $ \tilde \gl_{2,h}'\left(\overline{\gl}_1(h)^+\right)$ and a look at the right-hand side of 
\eqref{eq:fgammah}, recalling that $\overline{\gl}_1(\cdot)$ is analytic (cf. Remark~\ref{rem:hdep}), suffices to see that 
$h \mapsto  \tilde \gl_2'\left(\overline{\gl}_1^+\right)$ is analytic.  Since $\gamma_c(h)= -1/ \tilde \gl_2'\left(\overline{\gl}_1^+\right)$,
the properties of $\gamma_c(\cdot)$ claimed in the statement are proven.

We are now at the heart of the argument:
since the function to be maximized, $\gl_1+ \gamma \gl_2$ is constant on lines with slope $-1/\gamma$, the maximum is achieved
\begin{enumerate}
\item
at  a value $\gl_1\in (\overline{\gl}_1, 0)$ if $1/ \gamma>  -\tilde \gl_2'\left(\overline{\gl}_1^+\right)$: the value of $\gl_1$
is of course found by solving $\tilde\gl_2 '(\gl_1)=-1/\gamma$: this is what we call the Cram\'er regime (we have a tilted measure that makes typical the Large Deviation event that the renewal follows the slope $\gamma$);
\item
at $\gl_1= \overline{\gl}_1$
 if $1/ \gamma\le   -\tilde \gl_2'\left(\overline{\gl}_1^+\right)$: we are outside of the Cram\'er regime.
 \end{enumerate}
\smallskip

The situation is therefore that varying $h$ 
one may switch in and out of the Cram\'er regime.
Out of the Cram\'er regime  the free energy is actually equal to $\tg(h) - \overline \gl_1(h)$, 
and the function $h \mapsto \tg(h) - \overline \gl_1(h)$ is analytic for every $h>0$ (but it is not necessarily equal to the free energy
$\tilde \tf_\gamma (h)$!). 
In the  Cram\'er regime the free energy is strictly smaller than $\tg(h) - \overline \gl_1(h)$, just because the maximum 
of $\gl_1+ \gamma \gl_2$ is achieved on the boundary and $D_h(1, \gamma)$ contributes to the free energy with a negative
sign, cf. \eqref{eq:Fgc}. 
This explains why the changes of sign of 
$h\mapsto 1/ \gamma +\tilde \gl_2'\left(\overline{\gl}_1^+\right)$ are  non-analyticity points.
\qed

\subsection{Analysis of the free energy singularities in the localized regime}
\label{sec:singloc}

We now go deeper  into the analysis of the singularities for the non-analyticity points we have found, that is the values of $h>0$
for which $\gamma_c(h)-\gamma$ (cf. \eqref{eq:gammah}) changes sign. At the same time 
we tackle also the case in which $\gamma_c(h)-\gamma$ hits zero without crossing it.

For this and referring to the list of two items of the previous proof we introduce some notations.
\smallskip

\begin{enumerate}
\item In the Cram\'er regime we introduce 
the notation $ \left( \hat{\gl}_1(h,\gamma),\hat{\gl}_2(h,\gamma) \right) \in B_h$ 
 for the optimal point
 (for the sake of brevity, we omit the dependence on $\gamma$). It is of course found by 
solving $\tilde\gl_2 '(\gl_1)=-1/\gamma$ that yields  $\hat{\gl}_1(h)$ and 
$\hat{\gl}_2(h)= \tilde \gl_2\left(
\hat{\gl}_1(h)\right)$. Note that the values of $h$ corresponding to the Cram\'er regime is an union of 
open disjoint intervals: we call these intervals $\cI_{\textrm{1}}, \cI_{\textrm{2}},\ldots$. 
The analyticity of $\hat{\gl}_1(\cdot)$ and $\hat{\gl}_2(\cdot)$ in the intervals $\cI_j$ follows by Remark~\ref{rem:hdep} and by 
repeating the arguments in the same remark.
We also introduce  
\begin{equation}
\label{eq:hatc}
 \hat{c}_{\gamma}(h) \,:= \,  \left( \tg(h) - \hat{\gl}_1(h) \right) + \gamma \left( \tg(h) - \hat{\gl}_2(h) \right)\, ,
 \end{equation}
 for $h\in \cup_j \cI_j$ and of course $\tilde \tf_\gamma(h)=\hat{c}_{\gamma}(h) $ on this set.
 Note that
$\cup_j \cI_j$ can be alternatively characterized as $\{h: \, \gamma_c (h)-\gamma > 0\}$.
\item Out of the Cram\'er regime, that is for $h\in  \cI _{\textrm{b}}:= (0, \infty)\setminus  \cup_j \cI_j$, we introduce instead
\begin{equation}
\label{eq:ch}
c(h)\, := \, \tg(h) - \overline \gl_1(h)\, ,
\end{equation}
and $\tilde \tf_\gamma(h)=c(h)$ if and only if $h\in \cI _{\textrm{b}}$. Note that
$\cI _{\textrm{b}}$ can be alternatively characterized as $\{h: \, \gamma_c (h)-\gamma \le 0\}$. 
Note also that (cf.: the end of the proof of Prop.~\ref{th:reg})
$c(h)> \hat{c}_{\gamma}(h)$ in the interior on the intervals on which $ \hat{c}_{\gamma}(\cdot)$ is defined.
\end{enumerate}
\medskip

Observe moreover that    $ \partial \cup_j \cI_j= \partial \cI _{\textrm{b}}$, here  $\partial A$ denotes here the boundary of $A$
seen as subset of  $(0, \infty)$, and if $h \in  \partial \cI _{\textrm{b}}$  we have $\hat{c}_{\gamma}(h) = c(h)$ and $\gamma_c(h) = \gamma$. By differentiating \eqref{eq:gammah}  we obtain that
if $\gamma_c(h) = \gamma$
\begin{equation}
\label{eq:gammac'}
\gamma_{c}' (h)\, = \, -\frac{c'(h)}{\sum_{n,m} n K(n+m) \exp (-c(h) n )} \, \sum_{n,m} n (m- \gamma n) K(n+m) \exp (- n c(h)) \,,
\end{equation}
hence (with the the convention $\sign(0)=0$)
\begin{equation}
\label{eq:fgammac'}
\sign\left( \gamma_{c}' (h)\right)\, =\, -
\sign\left( \sum_{n,m} n (m- \gamma n) K(n+m) \exp (- n c(h)) \right)\, ,
\end{equation}
which is saying in particular that the sum in the right-hand side is zero if and only if $ \gamma_{c}' (h)$ is zero.

\medskip

In preparation of the next result, that investigates the regularity of the critical points in the positive semi-axis, 
it is useful to go through what may happen in the {\sl tangential cases}, namely   
when $\gamma_c(h_0)=\gamma$ and $\gamma_c'(h_0)=0$, for $h_0>0$. There are three different scenarios
\begin{enumerate}
\item  $\gamma_c(h_0+ \epsilon)-\gamma<0$ for every $\epsilon>0$ small, that is $h_0$ is a maximum, and in this case at $h_0$ there is no 
phase transition simply because $\tilde \tf_\gamma(h)=c(h)$ in a neighborhood of $h_0$, and $c(\cdot)$ is analytic in the positive semi-axis;
\item  $\gamma_c(h_0+ \epsilon)-\gamma>0$ for every $\epsilon>0$ small, that is $h_0$ is a minimum, and, as we will see in the next statement, 
$\tilde \tf_\gamma (\cdot)$ is at least $C^2$ at $h_0$, but we are not sure in general that $h_0$ is a critical point (see Remark~\ref{rem:postrefreg}), even if it looks very plausible that $h_0\in \partial \cI _{\textrm{b}}$ is a critical point of the free energy when $h_0$ is a minimum of $\gamma_c(\cdot)$;
\item $\gamma_c(\cdot)$ changes sign at $h_0$, that is $h_0$ is a saddle: in this case, as we have 
seen in Proposition~\ref{th:reg}, there is a transition at $h_0$ and, as we will point out in Remark~\ref{rem:postrefreg}, this transition is smoother than in the case
in which the derivative of $\gamma_c'(\cdot)$ is not zero.  
\end{enumerate}

\medskip 

\begin{proposition}
\label{th:refreg} 
Consider $h_0\in \partial \cI _{\textrm{b}}$, that is    $h_0>0$ such that $\gamma_c(h_0)=\gamma$. 
 The function $\tilde{\tf}'_{\gamma}(\cdot)$ is continuous at $h_0$,
 so a transition at $h>0$ is not of first order. If $\sum_{m} m^2 K(m)< \infty$,   $\tilde{\tf}^{''}_{\gamma}(\cdot)$ has a jump discontinuity (second order transition)
 at $h_0$ if and only if   $\gamma_c '(h_0)\neq 0$.
  If $\sum_{m} m^2 K(m)= \infty$, $\tilde{\tf}''_{\gamma}(\cdot)$ is continuous at $h_0$ (so the transition is of third order or more).
\end{proposition}
\medskip

Note that this statement says in particular that $\tilde{\tf}''_{\gamma}(\cdot)$ is continuous in full generality if
 $\gamma_c'(h_0)=0$.

\medskip

{\it Proof.} 
 Let us first prove that $\tilde{\tf}'_{\gamma}(h)$ is continuous at $h_0$. For $c(h)$, recall from Lemma ~\ref{th:Ah} that $ q_{h}\left( \overline{\gl}_1(h) , \tg(h)\right) = 1$ (see \eqref{eq:q_h} for the definition of $q_h(\cdot)$), therefore  
 \begin{equation}
 \frac{\partial}{\partial h} q_{h}\left( \overline{\gl}_1(h) , \tg(h)\right) \, = \, 0\, ,
 \end{equation}
  which directly implies that
\begin{equation}
\label{eq:c'}
  c'(h) \, \sum_{n,m} n K(n+m) \exp \left( - n c(h) \right)\, = \, \exp(-h) \,.
\end{equation}
For $ \hat{c}_{\gamma}(h)$, first replace $ \tg(h) - \hat{\gl}_1(h)$ by $ \hat{c}_\gamma (h) - \gamma \left( \tg(h) - \hat{\gl}_2(h)\right)$ (from \eqref{eq:hatc}) in $q_h(\cdot)$ to obtain $q_{h}\left( \hat{\gl}_1(h) , \hat{\gl}_2(h) \right)= \sum_{n,m} K(n+m) \exp\left(h - \hat{c}_{\gamma}(h)n - \left( \tg(h) - \hat{\gl}_2(h)\right) (m- \gamma n) \right)$. Recall that $ \tilde{\gl}_2'\left( \hat{\gl}_1(h)\right) = -1/{\gamma}$, which 
can be rewritten as
\begin{equation}
\label{eq:ID1}
\sum_{n,m} (m-\gamma n) K(n+m) 
\exp\left(h - \hat{c}_{\gamma}(h)n - \left( \tg(h) - \hat{\gl}_2(h)\right) (m- \gamma n) \right)\, = \, 0\, ,
\end{equation}
for every $h>0$. Keeping in mind this equality in evaluating $\frac{\partial}{\partial h} q_{h}( \hat{\gl}_1(h), \hat{\gl}_2(h)) = 0$, we obtain
\begin{equation}
\label{eq:wc'}
 \hat{c}'_{\gamma}(h) \, \sum_{n,m} n K(n+m) \exp\left(-\hat{c}_{\gamma}(h)n - \left( \tg(h) - \hat{\gl}_2(h)\right) (m- \gamma n) \right) \, = \, \exp(-h) \,.
\end{equation}
Since we have that $\left( \hat{\gl}_1(h_0),\hat{\gl}_2(h_0) \right) \, = \,\left( \overline \gl_1(h_0), \tg(h_0) \right)$ and $c(h_0) =  \hat{c}_{\gamma}(h_0)$ (from the definition of $ h_{c,\gamma}$), by evaluating \eqref{eq:c'} and \eqref{eq:wc'} at $h_0$ we get 
\begin{equation}
\label{eq:f=c}
\tilde{\tf}'_{\gamma}(h_0) \, = \,  c'(h_0) \, = \, \hat{c}'_{\gamma}(h_0) \, = \, {\left( \sum_{n,m} n K(n+m)\exp\left(h_0- n c(h_0) \right) \right)}^{-1}\, ,
\end{equation}
and  the continuity of $ \tilde{\tf}'_{\gamma}(h)$ at $h_0$ is proven.
\medskip

Now by differentiating  once again \eqref{eq:c'} we have

\begin{equation}
\label{eq:c''}
c''(h)= \frac{ - \exp(-h)+ {\left( c'(h)\right)}^2 \sum_{n,m} n^2 K(n+m) \exp \left( - n c(h) \right)}{\sum_{n,m} n K(n+m) \exp \left( - n c(h) \right)}  \, ,
\end{equation}
and by differentiating \eqref{eq:wc'} 
\begin{multline}
\label{eq:wc''}
\hat{c}''_{\gamma}(h)= \frac{ - e^{-h}+ 
\left( \hat{c}'_{\gamma}(h)\right)^2 
\sum_{n,m} n^2 K(n+m) 
e^{
-\hat{c}_{\gamma}(h)n - ( \tg(h) - \hat{\gl}_2(h)) (m- \gamma n) }
}
{\sum_{n,m} n K(n+m) e^{-\hat{c}_{\gamma}(h)n - \left( \tg(h) - \hat{\gl}_2(h)\right) (m- \gamma n) }}
 \\ + 
 \frac{ \hat{c}'_{\gamma}(h) \left( \tg'(h) - \hat{\gl}'_2(h) \right)
               \sum_{n,m} n(m-\gamma n)K(n+m) 
       e^{
               -\hat{c}_{\gamma}(h)n - \left( \tg(h) - \hat{\gl}_2(h) \right) (m- \gamma n)
          }
     }
{\sum_{n,m} n K(n+m) e^{-\hat{c}_{\gamma}(h)n - \left( \tg(h) - \hat{\gl}_2(h)\right) (m- \gamma n) }}
\, .
\end{multline}
We now observe that $c''(h_0)$ coincides with the first term in the right-hand side of \eqref{eq:wc''} evaluated at $h=h_0$. 
Therefore  $\tilde \tf''_\gamma(\cdot)$ is continuous at $h_0$ if and only if  the second term  in the right-hand side of \eqref{eq:wc''} vanishes at $h=h_0$. To clarify this issue we rewrite
 $ \tg'(\cdot) - \hat{\gl}'_2(\cdot)$ by exploiting the fact that by differentiating  \eqref{eq:ID1} with respect to $h$ we get
\begin{multline}
\label{eq:ID22}
\tg'(h) - \hat{\gl}'_2(h) \,=\\ - \hat{c}'_{\gamma}(h) \frac{\sum_{n,m} n(m-\gamma n)K(n+m) \exp\left( -\hat{c}_{\gamma}(h)n - \left( \tg(h) - \hat{\gl}_2(h)\right) (m- \gamma n) \right)}{\sum_{n,m}(m-\gamma n)^2 K(n+m) \exp\left( -\hat{c}_{\gamma}(h)n - \left( \tg(h) - \hat{\gl}_2(h)\right) (m- \gamma n) \right)} \,.
\end{multline}
Therefore going back to \eqref{eq:f=c}, \eqref{eq:c''} and \eqref{eq:wc''}, we see that if $\sum_{m} m^2 K(m) < \infty$ then
\begin{multline}
\label{eq:c''1}
\hat{c}''_{\gamma}(h_0) -   c''(h_0)\, = \, -
{\left( \tilde{\tf}'_{\gamma}(h_0)\right)}^3
\frac{ {\left( \sum_{n,m} n(m - \gamma n) K(n+m)\exp\left( h_0- n c(h_0) \right)\right)}^2}{\sum_{n,m} (m - \gamma n)^2 K(n+m)\exp\left( h_0- n c(h_0) \right)} \,,
\end{multline}
and if $\sum_{m} m^2 K(m) = \infty$ then $\hat{c}''_{\gamma}(h_0) -   c''(h_0)=0$.
So, if $\sum_{m} m^2 K(m) = \infty$ then  $\tilde \tf''_\gamma(\cdot)$ is continuous at $h_0$.
If $\sum_{m} m^2 K(m) < \infty$ then \eqref{eq:c''1} {\sl generically} tells us that $\tilde \tf''_\gamma(\cdot)$ has a jump discontinuity  at $h_0$, but the jump is zero if 
$\sum_{n,m} n(m - \gamma n) K(n+m)e^{ h_0- n c(h_0)}=0$ and, by \eqref{eq:fgammac'}, this is equivalent to 
$\gamma_c'(h_0)=0$.
The proof of Proposition~\ref{th:refreg} is therefore complete.
\qed

\medskip
 
\begin{rem}
\label{rem:postrefreg}
A sharper analysis of the singularity at $h_0\in \cI _{\textrm{b}}$  is possible
and one sees that the closer $\ga$ is to one the more the transition is regular. The general analysis however  is cumbersome
due also  to the fact that  the transition can gain regularity from cancellations that may appear and that depend 
on fine details: we have already found an instance of this
when in the proof of  Proposition~\ref{th:refreg} we have seen that if $\sum_{m} m^2 K(m) < \infty$ the second derivative of the free
energy has a jump at $h_0$ unless $\gamma_c'(h_0)=0$. 
%As we have anticipated in the list right before Proposition~\ref{th:efreg},since $\gamma_c'(h_0)=0$ says that $h_0$ can be a maximum, a minimum or a saddle point of $\gamma_c(\cdot)$,  the scenario in this case is really threefold. Notably 
These  cancellations are at the origin of the difficulties in resolving the issue in item (2) of
the same list.
\end{rem}

\medskip

\subsection{Free energy analysis of the delocalization transition}
\label{sec:deloc}

We complete now the proof of Theorem \ref{th:F} by studying the asymptotic behavior near $h_c$ of $\tilde{\tf}_{\gamma}(h)$:
we will show  in Section~\ref{sec:free} that $\tilde{\tf}_{\gamma}(h)= \tf_{\gamma}(h)$.

\medskip

We start by observing that, by 
 Proposition~\ref{th:Fgc},  \eqref{eq:Fgc} reduces to study the critical behavior of $\tg(h)$ and $D_h(1,\gamma)$. In the case
  $\gamma =1$, since $D_h(1,1)=0$, $\tf(h) = 2 \tg(h)$ and, since we have seen that
  the only singularity of $\tg(\cdot)$ is at the origin,  Theorem~\ref{th:F} reduces in this case to:
   
 \medskip
 
 \begin{lemma}
 \label{th:G}
 For $\ga>1$ we have
 \begin{equation}
 \tg(h) \stackrel{h \searrow 0}\sim
 \begin{cases}
 \frac 12 c h & \text{ if } \sum_{n} n^2 K(n) < \infty\, ,
 \\
  \frac 12 L_\ga (h)
 h^{1/(\ga-1)} & \text{ if }   \sum_{n} n^2 K(n) = \infty\, ,
 \end{cases}
 \end{equation}
 where $c$ is the same as for \eqref{eq:fch} and $ L_\ga (\cdot)$ is a slowly varying function.
 For $\ga=1$, $\tg(h)$ vanishes faster than any power of $h$. 
 \end{lemma}
 
 \medskip
 
 Implicit expressions for $L_\ga (\cdot)$ as well as a $\tg(\cdot)$ in terms of the inverse of a suitable slowly varying function
 in the case $\ga=1$
 can be found in the proof. 
 
 \medskip
 
 \noindent
 {\it Proof.}
  Actually since \eqref{eq:hc} can be written as
\begin{equation}
\label{eq:hc1}
\sum_{n \ge 2} (n-1) K(n) \exp( h- n \tg(h)) \, = \, 1 \,,
\end{equation}
we remark that $ \tg(h)$ is the free energy of the pinning model based on a one-dimensional renewal process and the proof is therefore just 
a revisitation of   \cite[Theorem~2.1]{cf:GB}. We give in any case a substantial part of the arguments here.

If $ \sum_n n^2 K(n) < \infty$, by (DOM) we have that $\sum_{n \ge 2} (n-1) K(n)\left( 1 - \exp(- n \tg(h))\right) \, \sim \, \tg(h) \sum_{n \ge 2} n(n-1) K(n) $ as $h \searrow 0$. Therefore by \eqref{eq:hc1}
\begin{equation}
\label{eq:g1}
\tg (h) \,  \sum_{n \ge 2} n (n-1) K(n)  \, \sim \, 1- \exp(-h) \, \sim \,h \,,
\end{equation}
and this of course proves the statement \eqref{eq:fch} for $\tilde{\tf}(h)$.

If $ \sum_n n^2 K(n) = \infty$ and $\ga \in (1,2)$, by \eqref{eq:hc1} and by Riemann sum approximation, one has 
\begin{multline}
\label{eq:RS3}
1 - \exp(-h)=
\sum_{n \ge 2} (n-1) K(n) \left( 1-\exp(- n \tg(h))\right)\\
\overset{h \searrow 0}{\sim} 
{\left( \tg(h) \right)}^{\ga -1} L\left( \frac{1}{\tg(h)}\right) 
\int_{0}^{\infty} \frac{1- \exp(-x)}{x^{\ga}} \, \mathrm{d}x \,,
\end{multline}
and
therefore
\begin{equation}
\label{eq:RS1}
\frac{1}{(\ga - 1)} \Gamma(2-\ga) {\left( \tg(h) \right)}^{\ga -1} L\left( 1/ \tg(h)\right)
 \, \overset{h \searrow 0}{\sim} \,
 h\, ,
\end{equation}
and $ \tg(h) \, \sim \,  \, \frac{1}{2} \, L_\ga (h) \, h^{1/(\ga-1)}$
where $ L_\ga (h) = 2 {\left( (\ga -1)/ { \Gamma(2 - \ga)}\right) }^{1/(\ga-1)} h^{-1/(\ga -1)} \hat{L}_{\ga}(h)$ and $\hat{L}_{\ga}(\cdot)$ is asymptotically equivalent to the inverse of $x \longmapsto x^{(\ga-1)} L(1/x)$.

For the case $\ga=2$ we restart with the right-hand side of  the first equality in \eqref{eq:RS3} which equals, up to an additive term $O(\tg)$  for $\tg \searrow 0$), to
\begin{equation}
\sum_{n\ge 2}\frac{L(n)}{n^2}\left( 1-\exp(-n\tg)\right) \, =\, 
\sum_{n= 2}^{\lfloor \gep/\tg \rfloor}\frac{L(n)}{n^2}\left( 1-\exp(-n\tg)\right) +
\sum_{n>\lfloor \gep/\tg \rfloor}\frac{L(n)}{n^2}\left( 1-\exp(-n\tg)\right)\, ,
\end{equation}
with $\gep>0$. By performing a Riemann sum approximation we see that the second term in the right-hand side is $O(\tg L(1/\tg))$.
For the first one instead we use that, by Taylor formula,  for every $\gd >0$ there exists $\gep>0$ such that
\begin{equation}
(1-\gd)\tg  \sum_{n= 2}^{\lfloor \gep/\tg \rfloor}\frac{L(n)}{n}\, \le \, 
\sum_{n= 2}^{\lfloor \gep/\tg \rfloor}\frac{L(n)}{n^2}\left( 1-\exp(-n\tg)\right)\, \le\, (1+\gd)\tg  \sum_{n= 2}^{\lfloor \gep/\tg \rfloor}\frac{L(n)}{n}\,,
\end{equation}
but the sum on the leftmost and rightmost term is asymptotic to $\int_1^{1/\tg} (L(t)/t) \dd t=: \breve{L}(1/\tg)$, which is  slowly varying and $L(x)=o(  \breve{L}(x))$ for
$x \to \infty$ \cite[Th. 1.5.9a]{cf:BGT}.
At this point we go back to the first equality in \eqref{eq:RS3} and we have
\begin{equation}
\label{eq:g2}
\tg(h) \tilde{L} \left( 1/( \tg(h)) \right) \stackrel{h\searrow 0} \sim  h \,.
\end{equation}
Since the right-hand side of the first equality in \eqref{eq:RS3} is an increasing function of $\tg$,
\eqref{eq:g2} can be asymptotically inverted 	and the case $\ga=2$ is complete.

For the case $\ga=1$ the computation is similar. Again we replace the term $(n-1)$ with $n$ in the right-hand side of the first equality in \eqref{eq:RS3}:
the error is $O(\tg)$. Then we are left with
$
\sum_{n\ge 2} (L(n)/n) (1-\exp(-n\tg))
$ and we split the sum into $n$ smaller and larger than $1/(\gep \tg)$. The first sum can be treated by Riemann approximation
yielding a term $O(L(1/\tg))$. The other term instead is asymptotic to $\check L(x):=\int_x^\infty (L(t)/t)$, which is slowly varying and $L(x)=o(   \check L(x))$
\cite[Th. 1.5.9b]{cf:BGT}. 
The right-hand side of the first equality in \eqref{eq:RS3} is an increasing function of $\tg$ that vanishes as $\tg$ tends to zero.
So $\breve L(x)$ can be chosen decreasing (to zero) and the slowly varying property implies that $\breve L (x) \ge x^{-\gep}$
for every $\gep>0$ and every $x$ sufficiently large. Hence $\breve L(1/\tg(h)) \ge \tg(h)^\gep$ for $h$ sufficiently small,
$\breve L(1/\tg(h))\sim h$ readily implies that $\tg(h)=O(h^{1/\gep})$ which is what we wanted to prove.
 %\begin{equation}
%\label{eq:g11}
%\tg(h) \, \overset{h \searrow 0}{\sim} \, 1/ \check{L}^{-1}(h) \,
%\end{equation}
 The proof of Lemma~\ref{th:G} is  therefore complete.
\qed 

\medskip

\noindent
{\it Proof of Lemma~\ref{th:gammac}.}
The analyticity of $\gamma_c(\cdot)$ has been proven in Proposition~\ref{th:reg} and the second statement of \eqref{eq:gammac-th} is trivial.
Let us prove the first statement of \eqref{eq:gammac-th} (keeping in mind that $\ga>1$).

If $\sum_{n \ge 1}n^2 K(n) < \infty$, recall the definition of $\gamma_c(\cdot)$ from \eqref{eq:gammah}, it is easy to see that 
\begin{equation}
\gamma_c(0) = \frac{\sum_{n,m \ge 1} m K(n+m)}{\sum_{n,m \ge 1} n K(n+m)} = 1 \,.
\end{equation}

Now if $\sum_{n \ge 1}n^2 K(n) = \infty$ and $\alpha \in (0,1)$, set $\tg(h) - \overline{\gl}_1(h) = x=o(1)$ (as $h \searrow 0$) and  remark that
\begin{equation}
\label{eq:gac0}
\sum_{n,m} n K(n+m) \exp \left( - x n \right)  \,=\,
\sum_{t \ge 2} K(t) \, \frac{e^{ x} (1- e^{- x t } )  + t e^{ x(1- t)} (1 - e^{x})}{{\left(e^{ x} - 1\right)}^2} \,.
\end{equation}
By Riemann Sum approximation, the right-hand side of \eqref{eq:gac0} is equivalent to (as $x \searrow 0$)
\begin{equation}
x^{\ga -2} L \left( 1/x\right) 
\int_{0}^{\infty} \frac{1 -  \exp(- y) (1+y) }{y^{\ga + 1}} \, \mathrm{d}y \,,
\end{equation}
therefore
\begin{equation}
\label{eq:gac01}
\sum_{n,m \ge 1} n K(n+m) \exp \left( - x n \right) \stackrel{x \searrow 0}\sim 
x^{\ga -2} \frac{\Gamma(2-\ga)}{\ga} L \left( 1/x\right) \,.
\end{equation}
Repeating the same argument leading to \eqref{eq:gac01}, we see that
\begin{equation}
\label{eq:gac02}
\sum_{n,m \ge 1} m K(n+m) \exp \left( - x n \right) \stackrel{x \searrow 0}\sim \\
x^{\ga -2} \Gamma(-\ga) L \left( 1/x\right) \,.
\end{equation}
By \eqref{eq:gammah},\eqref{eq:gac01} and \eqref{eq:gac02}, we get
\begin{equation}
\label{eq:gac03}
\gamma_c(h) \stackrel{h \searrow 0} \sim \frac{1}{\ga -1} \,.
\end{equation}

For the case $\ga = 2$, it suffices to show that
\begin{equation}
\label{eq:mncomp1}
\sum_{n,m} n K(n+m) \exp(-xn) \stackrel{x \searrow 0} \sim
\sum_{n,m} mK(n+m) \exp(-xn)\, .
\end{equation}
In fact, both terms are asymptotic to $\breve{L}(x)= \int_1^{1/x} L(t)/t$, slowly varying
by \cite[1.5.9b]{cf:BGT} and diverging at $\infty$ because $\sum_{n,m}nK(n+m)=\infty$. This can be seen by restarting from \eqref{eq:gac0}:
 the left-hand side of \eqref{eq:mncomp1} is equal to
\begin{multline}
\label{eq:gac0.1}
%\sum_{n,m} n K(n+m) \exp \left( - x n \right)  \,=\,
\frac{e^x}{{\left(e^x - 1\right)}^2}\sum_{t \ge 2} \frac{L(t)}{t^3}  \left((1- e^{- x t } )  + t e^{ -x t} (1 - e^{x})\right) \,=
\\
 \frac {1+O(x)}{x^2} \sum_{t \ge 2} \frac{L(t)}{t^3}  \left((1- e^{- x t } )  - tx e^{ -x t} \right)  +O(1)\, .
\end{multline}
A Riemann sum approximation shows that if the sum is limited to $n>\gep /x$, for an arbitrary $\gep>0$, yields an
$O(1)$ contribution.  For the term that is left we use  that, by Taylor formula, for every $\gd>0$ one finds a $\gep>0$ such that
\begin{equation}
\frac 1{x^2}
\sum_{t =2}^{\lfloor\gep/x\rfloor} \frac{L(t)}{t^3}  \left((1- e^{- x t } )  - tx e^{ -x t} \right)\, \le
\, \frac{1+ \gd}{x^2} \sum_{t =2}^{\lfloor\gep/x\rfloor} \frac{L(t)}{t^3} (tx)^2\stackrel{x \searrow 0}\sim\, (1+\gd)
 \sum_{t =2}^{\lfloor\gep/x\rfloor} \frac{L(t)}{t} \, ,
\end{equation}
and analogous lower bound with $1-\gd$. Since $\gd>0$ is arbitrary, the claimed asymptotic behavior
for the left-hand side of \eqref{eq:mncomp1} is established. The computation  for the
right-hand side is very similar and left to the reader. Therefore the proof is complete.
%by using \eqref{eq:kok} again, we see that 
%\begin{multline}
%\sum_{n,m \ge 1} n K(n+m) \exp \left( - \overline{a}_1 \tg(h) n \right) \stackrel{h \searrow 0}\sim \\\sum_{n \ge 0} \overline{\mathrm{K}}(n,0) \exp \left( - n \overline{a}_1 \tg(h)\right) - (\overline{a}_1 \tg(h) ) \sum_{n \ge 0} (n+1) \overline{\mathrm{K}}(n,0) \exp \left( - n \overline{a}_1 \tg(h)\right)\,.
%\end{multline}
\qed  

\medskip 

Recall now Proposition~\ref{th:Fgc} and in particular the fact that
$D_h(1, \gamma)$ is the result of an optimization problem, cf. \eqref{eq:FgcDh}, over
the set $B_h$ (cf. \eqref{eq:Bh}). As widely used and discussed in and right after Proposition~\ref{th:reg}, 
 the maximum can be achieved in the interior of $B_h$ (Cram\'er regime) or at the boundary (out of  Cram\'er regime).

\medskip

\begin{proposition}
\label{th:Fprop}
Choose $\gamma>1$.
If  $\sum_{n} n^2 K(n)< \infty$ we have that for $h$ small the system is outside of the Cram\'er regime and 
\begin{equation}
\label{eq:fch-2}
\tilde \tf_\gamma(h) \stackrel {h \searrow 0}\sim \tilde \tf_1 (h) \, =\, 2 \tg(h)\,  .
\end{equation}
%with $c^{-1}:=  {\frac 12\sum_{n=2}^\infty n(n-1) K(n)}$. 
If instead $\sum_{n} n^2 K(n)= \infty$ for $h$ small
the system is in the Cram\'er regime if $\gamma < 1/(\ga -1)$ and there exists $c_{\ga,\gamma} \in (1,\frac 12({\ga}^{1/(\ga-1)} \wedge (1+\gamma)))$ such that
\begin{equation}
\label{eq:fch-3}
\tilde \tf_\gamma(h) \sim {c_{\ga,\gamma}} \tilde \tf_1 (h)\, . 
%\ \text{ and } \ \tilde \tf_1 (h)\, \sim  \, \hat L(h) h^{1/(\ga-1)}  \,.
\end{equation}
If $\gamma >1/(\ga -1)$, the system is outside of this regime and 
\begin{equation}
\label{eq:fch-4}
\tilde \tf_\gamma(h) \sim   \frac{{\ga}^{1/(\ga-1)}}{2} \tilde \tf_1 (h) \,.
\end{equation}
If $\ga=1$ the system is in the Cram\'er regime for every $\gamma$ and
$\tilde \tf_\gamma(h)=O(h^{1/\gep})$ for every $\gep>0$.
The asymptotic behavior of $\tg(h)$
is given in Lemma~\ref{th:G}.
\end{proposition}
\medskip

\noindent{\it Proof.}
%Recall that $\tilde \tf_1(h) =  2 \tg(h)$, therefore the asymptotic behavior of $\tilde \tf_1(h)$ is an immediate consequence of Lemma~\ref{th:G}. So let us focus on 
Since $\gamma >1$ and 
 $D_h(1,\gamma) > 0$. 
  We recall Proposition~\ref{th:Fgc} and  for $(\lambda_1,\lambda_2) \in B_h$, we make the change of variables $ \tg(h) - \lambda_1 =  a_1 \, \tg(h)$, so $a_1\ge 1$, and $ \tg(h) - \lambda_2 =  a_2 \, \tg(h)$, so $a_2 \in [0,1]$.  
 With these new variables  we have
\begin{equation}
\label{eq:Da}
%\begin{split}
 D_h(1,\gamma)\,%&
 =\, \tg(h) \, \max_{ a \in \mathcal{B}_h} \left(1 - a_1+ \gamma (1- a_2)\right) 
 %\\&
 =\, \tg(h) \left( 
 1+ \gamma - \min_{ a \in \mathcal{B}_h} \left( a_1+ \gamma a_2\right)
 \right) 
 \,,
% \end{split}
\end{equation}
hence
\begin{equation}
\label{eq:DaF}
\tilde \tf_\gamma (h)\, =\, \min_{ a \in \mathcal{B}_h} \left( a_1+ \gamma a_2\right) \tg(h) \, ,
\end{equation}
with 
\begin{equation}
\label{eq:BBh}
\mathcal{B}_h\, :=\, \left \{ a :\, 1 \le a_1\le \frac{\overline{\gl}_1(h)}{\tg(h)},  \, 0 \le a_2 \le 1 \, , \sum_{n,m}  K(n+m) \exp\left( h - ( a_1 n + a_2 m)\tg(h) \right) \, =1 
\right\}\, .
\end{equation}
We set $ \Psi_h(a,\tg) \, := \, \sum_{n,m} K(n+m) \exp\left( - ( a_1 n + a_2 m)\tg(h) \right)$ and for every $ a \, \in \mathcal{B}_h$, we have that $\Psi_h(a,\tg(h)) \,=\, \exp(-h)$, which we will use asymptotically as
\begin{equation}
\label{eq:psi}
1 - \Psi_h(a,\tg(h)) 
\, \overset{h \searrow 0}{\sim} \,
h\, .
\end{equation}

If $ \sum_n n^2 K(n) < \infty$, since $\gamma_c(0) = 1$ (from Lemma~\ref{th:gammac}) and $\gamma > 1$, the system is outside of the Cram\'er regime for $h$ small (recall that the system is in the Cram\'er regime if and only if $\gamma_c(h)> \gamma$) and in the new variable means that the minimizers satisfies  $a_2=0$. Set $\overline{a} = (\overline{a}_1,0)$ with $\overline{a}_1 \tg(h) = \tg(h) - \overline{\gl}_1(h)$.
By (DOM) we have 
\begin{equation}
1 - \Psi_h(\overline{a},\tg(h))= 1 - \sum_{n,m \ge 1} K(n+m) \exp \left( - \overline{a}_1 n \tg(h) \right) 
\overset{h \searrow 0}{\sim} \,
\frac{1}{2} \overline{a}_1 \tg(h) \sum_{n \ge 2} n (n-1) K(n) \,,
\end{equation}
which, together with \eqref{eq:g1} and \eqref{eq:psi}, implies $ \overline{a}_1 \overset{h \searrow 0}{\sim} 2 $. Therefore , by  \eqref{eq:DaF},   $\tilde \tf_\gamma(h)\sim 2 \tg(h)$.

\medskip

When $ \sum_n n^2 K(n) = \infty$, we first treat the case $\ga \in (1,2)$: we have seen in Lemma~\ref{th:gammac} that in this case $\gamma_c(0) = 1/(\ga-1)$ which implies that we are in the Cram\'er regime for $h$ small if $\gamma < 1/(\ga-1)$ and outside if $\gamma >   1/(\ga-1)$  (and $a_2 = 0$). But let us consider the two regimes together for now and 
observe  that for $a_1\in [0,1]$, $a_2\ge 1$ and $a_1\neq a_2$
\begin{multline}
\label{eq:RS2}
1 - \Psi_h(a,\tg) \, = \, \\ \sum_{t \ge 2} K(t) 
\left( (t-1) -  \left( \frac{\exp\left( - t a_1 \tg(h) \right)- \exp\left( - t a_2 \tg(h)  - (a_1 -a_2 ) \tg(h)\right)}{\exp\left( - (a_1 -a_2 ) \tg(h)\right) -1}\right) \right)\, ,
\end{multline}
and by  Riemann sum approximation, the right-hand side is equivalent to (as $h \searrow 0)$
\begin{equation}
{\tg(h)}^{\ga -1}  \,L\left( 1/\tg(h)\right) \,
\int_{0}^{\infty} \frac{x -  \left( \exp(- a_1 x) - \exp(- a_2 x))/ (a_2 - a_1) \right) }{x^{\ga + 1}} \, \mathrm{d}x \,.
\end{equation}
The integral can be made explicit:
\begin{equation}
\label{eq:expli}
1 - \Psi_h(a,\tg)
\, \sim \,
b_\ga(a) \Gamma(- \ga)\, {\tg(h)}^{\ga -1}  \,L\left( 1/\tg(h)\right) \,,
\end{equation}
with 
\begin{equation}
\label{eq:balpha}
 b_\ga(a): = \frac{\left( {a_1}^\ga - {a_2}^\ga \right)}{  (a_1 - a_2)}\,=\, \ga \int_0^1 \left((1-t)a_2+ ta_1\right)^{\ga-1}
 \dd t\, .
 \end{equation} 
 By \eqref{eq:RS1} and \eqref{eq:psi} we get that for $a\in \cB_h$ and $a_1\neq a_2$
\begin{equation}
\label{eq:ca}
b_\ga(a) \, \sim \, \ga \, ,
\end{equation}
 but 
 the rightmost term in \eqref{eq:balpha}
  shows that the singularity in $a_1=a_2$ is removable and one directly verifies 
 that \eqref{eq:expli}, and therefore \eqref{eq:ca}, hold  also for $a_1=a_2(=1)$.
 
 A number of considerations are in order:
 \begin{enumerate}
 \item By recalling the convexity arguments used in 
 \S~\ref{sec:constr1}, we directly have that 
 the constraint in $\cB_h$ can be rewritten as $a_{2}= \tilde a_{2,h}(a_1)$, with $\tilde a_2(\cdot)$
 a decreasing convex function (this is the curve appearing in Fig.~\ref{fig:Bh}, up to the affine change of variable we performed).
 It will be more practical at this stage to write rather $a_1= \tilde a_{1,h}(a_2)$, and $\tilde a_{1,h}:[0,1]\to [1, \infty)$
 is also convex and decreasing with $\tilde a_{1,h}(0)=\overline{\gl}_1(h)/\tg(h) \sim \ga^{1/(\ga-1)}$ and
   $\tilde a_{1,h}(1)=1$. The choice in favor of $\tilde a_{1,h}(\cdot)$ over $\tilde a_{2,h}(\cdot)$ is because we prefer having a $h$ dependence in the image rather than in the domain.
 \item Since 
 $1-\Psi_h(\cdot, \tg)$ is concave, also $b_\ga(\cdot)$ is (this can also be verified directly) and the constraint in $\cB_h$,
 namely \eqref{eq:ca}, becomes in the limit $b_\ga(a) = \ga$, with $a_1\in [1,  \ga^{1/(\ga-1)}]$ and $a_2 \in [0,1]$. Such a constraint 
 can be repressed (like in \S~\ref{sec:constr1}) as $a_1= \tilde a_{1}(a_2)$, where $\tilde a_1:[0,1] \to [1,  \ga^{1/(\ga-1)}]$ is a convex decreasing function. This function is smooth (in fact, analytic, by the Implicit Function Theorem \cite{cf:Primer})
in the interior of the domain and one directly computes $\tilde a'_{1}(0^+)= -1/(\ga-1)$ and $\tilde a'_{1}(1^-)=-1$, which actually coincide with the limits as $h \searrow 0$ respectively of $\tilde a'_{1,h}(0^+)$ (this is precisely \eqref{eq:gac03})
 and $\tilde a'_{1,h}(1^-)=-1$ (in fact: $\tilde a'_{1,h}(1^-)=-1$ for every $h>0$ by symmetry).
 \item But \eqref{eq:ca}, that is the convergence of the constraint function for $h\searrow 0$, implies 
 $\lim_h a_{1,h}(a_2)=a_1(a_2)$ for every $a_2 \in [0,1]$, as well as the convergence of the derivative for
 $a_2 \in (0,1)$, because we are dealing with a sequence of convex functions and because the limit is differentiable. Note that 
 we have already pointed out that there is convergence of the derivatives also at the endpoints and that, by analyticity, 
 $a_1(\cdot)$ is strictly convex. This allows to conclude that 
the unique minimizer $\hat a(h)$ for   \eqref{eq:DaF} converges to the unique minimizer of 
$\min_{a \in \cB} a_1+ \gamma a_2$, with
\begin{equation}
\label{eq:BBhlim}
\mathcal{B}\, :=\, \left \{ a :\, 1 \le a_1\le \ga^{1/(\ga-1)},  \, 0 \le a_2 \le 1 \, , b_\ga(a)=\ga 
\right\}\, .
\end{equation}
The limit problem, like the approaching ones, can of course be rewritten in a one-dimensional form
using $\tilde a_{1}(\cdot) $ and $\tilde a_{1,h}(\cdot) $.
 \end{enumerate}
 
 \medskip 
 
At this point we can treat separately the non Cram\'er case, in which  
  $a_2 = 0$ (this happens when  $\gamma > 1/(\ga-1)$, but also when  $\gamma = 1/(\ga-1)$ by a limit procedure), so we have $a_1 = {\ga}^{1/(\ga-1)}$ and $\tilde \tf_\gamma (h) \sim  {\ga}^{1/(\ga-1)} \tg(h)$.
  
%If $\gamma < 1/(\ga-1)$ instead, $\tilde \tf_\gamma (h) \sim \min_{a \in \cB} (a_1+ \gamma a_2) \tg(h)$. Since $a_1(\cdot)$ is strictly convex and decreases from $a_1(0)= {\ga}^{1/(\ga-1)}$ to $a_1(1)= 1$, one can easily see (use $\min_{a \in \cB} (a_1+ \gamma a_2)) $ that $ \tilde \tf_\gamma (h) \sim  c_{\ga,\gamma} \tg(h)$ for some $ c_{\ga,\gamma} \in ( {\ga}^{1/(\ga-1)} \wedge (1+\gamma), {\ga}^{1/(\ga-1)} \vee (1+\gamma))$.
 
  If $\gamma < 1/(\ga-1)$ instead, $\tilde \tf_\gamma (h) \sim \min_{a \in \cB} (a_1+ \gamma a_2) \tg(h)$,
  so $c_{\ga, \gamma}=\frac 12 \min_{a \in \cB} (a_1+ \gamma a_2)$
  and \eqref{eq:fch-3} is proven.
   Since $a_1(\cdot)$ is strictly convex and decreases from $a_1(0)= {\ga}^{1/(\ga-1)}$ to $a_1(1)= 1$, one can easily see 
   (just evaluate $ a_1+ \gamma a_2$ at $a=(1,1)$ and $a=( {\ga}^{1/(\ga-1)},0)$)
 that  $c_{\ga,\gamma} \in ( 1,\frac 12({\ga}^{1/(\ga-1)} \wedge (1+\gamma)))$. 
 
 \medskip

\begin{rem}
\label{th:lag}
For $\ga \in (1,2)$ we  can use  the Lagrange multiplier method to solve the limit optimization problem $\min_{a \in \cB} a_1+ \gamma a_2$.
With $s$ as multiplier we have
\begin{equation}
\label{eq:lag}
 \begin{cases}
(a_2 - a_1) \, = \, s \, \ga \, ({a_1}^{\ga - 1} -1 ) \, ,
 \\
\gamma (a_2 - a_1) \, = \, s \, \ga \,(1 - {a_2}^{\ga - 1}  ) \, ,
 \end{cases}
 \end{equation}
which implies that (divide the two equations in \eqref{eq:lag} and use
$b_\ga(a)=\ga$)
\begin{equation}
\label{eq:a2}
 {a_2}^{\ga -1} \, =\,  1 - \gamma ( {a_1}^{\ga -1} -1) \,.
\end{equation}
In particular, and consistently with what precedes, 
 if $\gamma > 1/( {a_1}^{\ga -1} -1) \ge 1/(\ga-1)$ no solution is found and we are out of the Cram\'er regime  $(a_1,a_2) = ({\ga}^{1/(\ga-1)},0)$. If $\gamma <   1/(\ga-1)$ instead a solution is found and in fact we are in the Cram\'er regime.
% Here is the explicit solution for $\ga=3/2$:
%\begin{equation}
%a_1 = \left(\frac{(\gamma+1) (2 \gamma -1)}{2 (\gamma^2 - \gamma +1)}\right)^2 \, \hspace{1cm} \, {a_2} \, =\, \left( \frac{(\gamma+1) ( 2 -\gamma)}{2 (\gamma^2 - \gamma +1)} \right)^2\,.
%\end{equation}
\end{rem}

\medskip

For the case $\ga=2$, from Lemma~\ref{th:gammac}, we have $\gamma_c(0)=1$, which implies that for every $\gamma >1$, we are outside of the Cram\'er regime for $h$ small and $a= \overline{a}=(\overline{a}_1,0)$. Observe that 
\begin{equation}
\label{eq:overa}
1 - \Psi_h(\overline{a},\tg) = \sum_{t \ge 2} \frac{L(t)}{t^3} \left( (t-1) - \frac{\exp( - \overline{a}_1 \tg(h) (t-1) )-1}{\exp( - \overline{a}_1 \tg(h) )-1} \right) \,,
\end{equation}
and follow the same procedure adopted for the case $\ga =2$ in Lemma~\ref{th:G}: split the sum in 
\eqref{eq:overa} into $t$ larger than $\gep / \overline{a}_1 \tg$, for an arbitrary $\gep>0$ (which yields $O(1)$ contribution) and to $t$ smaller than $\gep / \overline{a}_1 \tg$ to obtain
\begin{equation}
1 - \Psi_h(\overline{a},\tg) \stackrel{h \searrow 0} \sim \frac{1}{2} \overline{a}_1 \tg(h) \tilde{L} \left( 1/( \overline{a}_1 \tg(h)) \right) \,.
\end{equation}
We know that the left-hand side is equivalent to $h$ from \eqref{eq:psi}, and using the property of slowly varying functions \eqref{eq:psv} ( $\tilde{L} \left( 1/( \overline{a}_1 \tg(h))\right) \sim \tilde{L} \left( 1/(  \tg(h)) \right)$), then by \eqref{eq:g2} we get $ \overline{a}_1 \sim 2$. Therefore by \eqref{eq:DaF}, we obtain $\tilde \tf_\gamma(h)\sim 2 \tg(h)$.

\medskip

For the case $\ga = 1$, we do not strive for the sharp prefactor, since we did not go for the sharp behavior of $\tg(h)$, cf. Lemma~\ref{th:G}.
We just use \eqref{eq:DaF} and observe that $\min_{a \in \cB_h} a_1+\gamma a_2$ is bounded (because the supremum 
is over a set that is bounded uniformly in $h$, with $h$ in a right neighborhood of zero). Hence in this case $\tilde \tf_\gamma(h)=O(h^{1/\gep})$
for every $\gep>0$. The proof is therefore complete.
\qed

\subsection{The bio-physics model and other examples}
\label{sec:ex}

In this section we make $ \gamma_{c}(\cdot)$ explicit for some choices of $K(\cdot)$, starting with \cite{cf:NG}.

\medskip

\subsubsection{The bio-physics model}
We refer to  Section~\ref{sec:biop}. We have seen that the geometric constant $s$ is irrelevant then we take $s=1$ and $B(l)= 1/l^{c}$. Recall that $E_b>0$ is the binding energy and $E_l>0$ is the loop initiation cost. Set
\begin{equation}
K_h(n) := \frac{c_K}{n^{c}} \exp \left( h (E_b -E_l \mathbf{1}_{n>2})- n \tg(h) \right) \,,
\end{equation}
with $c_K = \sum_{n,m \ge 1} 1/(n+m)^c = \sum_{n \ge 2} (n-1)/n^c$ is the normalization constant and $\tg(h)$ is the only solution to $\sum_{n,m} K_h(n+m) = 1$. 

\begin{figure}[ht]
\label{fig:biop}
\centering
\scalebox{0.7}{\includegraphics*{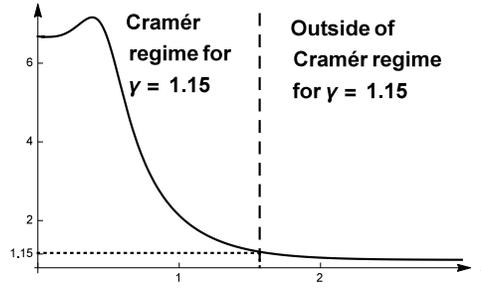}}
\caption{A representation of the function $\gamma_c(\cdot)$ for the distribution $K_h(\cdot)$ with $c=2.15$, $E_b=6$ and $E_l=3$. The horizontal dashed line corresponds to the value of $\gamma=1.15$ and the vertical one to the critical point $h_{c,\gamma}\simeq 1.676$.}
\end{figure}

In Figure~\ref{fig:biop}, a plot of $\gamma_c(\cdot)$ for values of $E_B,E_l,c$ and $\gamma$ chosen in \cite{cf:NG}, shows that the system exhibits an unique transition at $h_{c,\gamma} \simeq 1.676$ (the vertical dashed line). For $\gamma = 1.15$, observe that the system is in the Cram\'er regime for $h<h_{c,\gamma}$ and outside of the Cram\'er regime if $h \ge h_{c,\gamma}$.
In \cite{cf:NG}   the transition at $h_{c,\gamma}$ is described as between a phase with microscopic free strands at the end of 
the polymer and macroscopic free strands.

\subsubsection{A more explicitly solvable class: basic case}
We exploit the inter-arrival law used in \cite{cf:GB-ejp} in the one-dimensional case.
Recall that Euler's Gamma function $ \Gamma (z) := \int_{0}^{\infty} t^{z-1} \exp(- t) \,\mathrm{d}t\,$ defines an analytic function on $ \lbrace z \in \mathbb{C}: Re(z) > 0 \rbrace$ and that for $z \in \mathbb{C} \setminus \lbrace 0,-1,-2,... \rbrace$, it verifies $\Gamma(z+1) = z \Gamma(z)$ (therefore $\Gamma(n+1)= n!$ for $n \in \mathbb{N})$. Recall also that the Taylor coefficients of the function $ (1 - z)^{\ga}$ for $ \vert z \vert < 1$ and $ \ga \in \mathbb{R} \setminus \lbrace 0,1,2,...\rbrace$, is known exactly
\begin{equation}
\label{eq:z}
\sum_{n \ge 0} \frac{\Gamma(n-\ga)}{ n!} z^n = \Gamma(-\ga) (1-z)^{ \ga} \,,
\end{equation} 
and asymptotically from Stirling's formula we have $ \Gamma(n-\ga) / n! \overset{n \to \infty}{\sim} 1/n^{1+\ga}$. Note that the first terms of the series in \eqref{eq:z} have alternating signs.

First, let us suppose that $\ga \in (1,2)$ and set 
\begin{equation}
K_1(n) \, := \, \frac{\Gamma(n-\ga)}{\Gamma(-\ga)n!}\,, \hspace{0.2 cm} \text{for} \hspace{0.2 cm} n \ge 2 \,
\end{equation}
and from \eqref{eq:z} we have that $ \sum_{n \ge 2} (n-1) K_1(n) = 1$ and $ \gamma_{c}(h) = 1 /(\ga -1)$. This implies that there is no transition in the localized phase and for every $ h > 0$, the free energy is either -- recall \eqref{eq:ch} --
$c(h)$ (if $ \gamma \ge \gamma_c(h)$: out of the Cram\'er regime) or  $\hat{c}_{\gamma}(h)$ (if $\gamma <  \gamma_c(h) $, the Cram\'er regime, recall
\eqref{eq:hatc}).

\subsubsection{A more explicitly solvable class: general case}
\begin{figure}[ht]
\centering
   \begin{minipage}[b]{.46\linewidth}
     \scalebox{0.7}{\includegraphics*{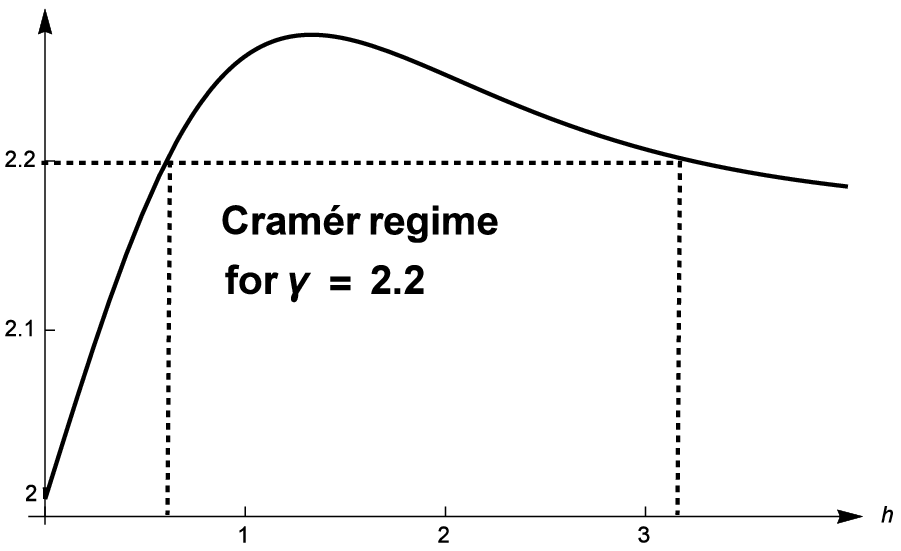}}
    \end{minipage}
     \quad
   \begin{minipage}[b]{.46\linewidth}
    \scalebox{0.7}{\includegraphics*{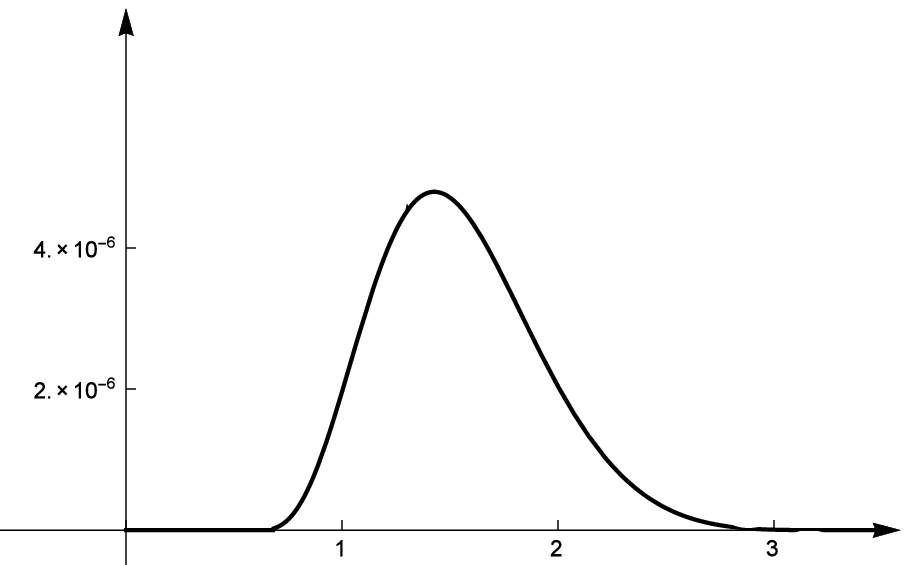}}
   \end{minipage}\hfill
\caption{\label{fig:alpha1} On the left  the function $\gamma_c(\cdot)$ for the inter-arrival distribution $K_2(\cdot)$, with $\ga= 1.5$ and $ \kappa = 0.01$: for $\gamma = 2.2$. There are two values of $h$ (vertical dashed lines) such that $\gamma_c(h) = \gamma$ and therefore there are two transitions. On the right we plot the difference $ c(\cdot) - \tilde \tf_{\gamma}(\cdot)$, that is $ c(\cdot) -\hat c_\gamma(\cdot)$ (recall \eqref{eq:hatc} and \eqref{eq:ch}), which makes clear the presence of the two transitions. The resolution of the graph does not allow to appreciate the positivity of such a difference for example at $h=3$ where it is about $6\times 10^{-9}$.}
\end{figure}

More general explicit cases can be built by modifying $K_1(\cdot)$ on a   finite number of sites.
Let us choose now $K_2(2) = K_1(2)$, $K_2(3)= \kappa$ and normalize the rest 
\begin{equation}
K_2(n) = c_{\kappa} \, K_1(n) \,, \hspace{0.2 cm} \text{for} \hspace{0.2 cm} n \ge 4 \,
\end{equation} 
with $ c_{\kappa} = \left( 1 - K_1(2) - 2 \kappa \right) / \left( 1 - K_1(2)-  2 K_1(3) \right)$ and $\kappa \in [0,1)$. With this choice of the inter-arrival $K(\cdot)$ there are  two transitions in the localized phase 
if an only if $\gamma>2$, see Figure~\ref{fig:alpha1}: the system here is in the Cram\'er regime only for intermediate values of $h$.
\noindent

\begin{figure}[ht]
\centering
   \begin{minipage}[b]{.46\linewidth}
     \scalebox{0.7}{\includegraphics*{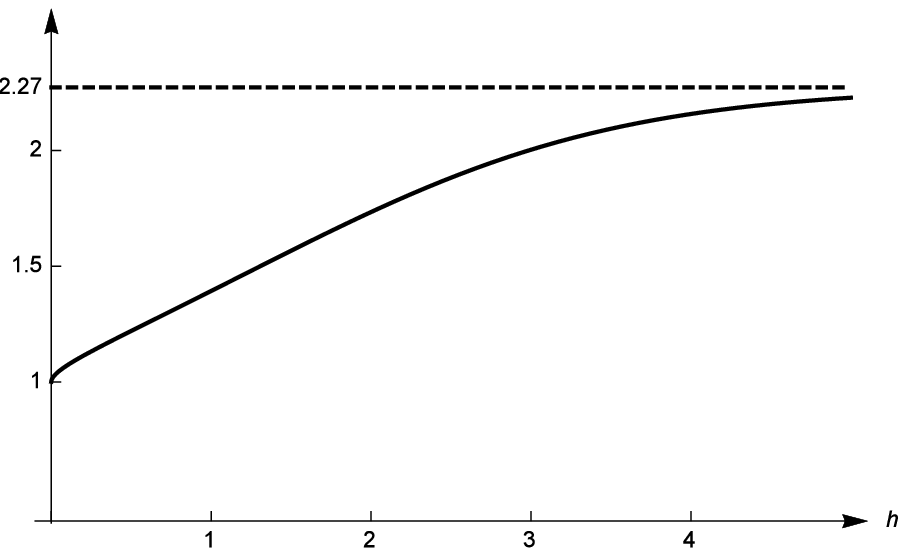}}
    \end{minipage}
     \quad
   \begin{minipage}[b]{.46\linewidth}
    \scalebox{0.7}{\includegraphics*{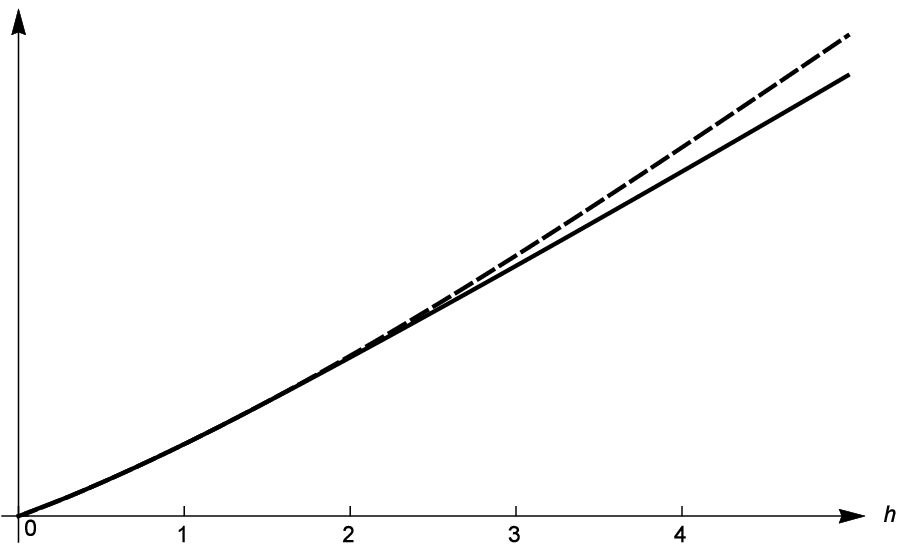}}
   \end{minipage}\hfill
\caption{\label{fig:alpha22} For the inter-arrival distribution $K_3(\cdot)$, the function $\gamma_c(\cdot)$ (on the left) is stricly increasing (for $\ga = 2.5, \varrho = 0.02$), this implies that if $\gamma_0 := \lim_{h \to \infty} \gamma_c(h) \le \gamma$ (dashed line represent the value of $\gamma_0$), there is no transition in the system. On the right, we present the free energy $\tilde \tf_\gamma(h)$ and $c(h)$ (dashed line, recall \eqref{eq:ch}) for $\gamma = 1.5$, the critical point $h_{c,\gamma}$ is represented by the dashed vertical line\,. }
\end{figure}

We then present an example with $\ga \in (2,3)$: Proposition~\ref{th:Fprop} shows that in this case for $h$ small the system
is outside the Cram\'er regime.
 We take  $K_3(2) = \varrho$ and $K_3(n) = c_{\varrho} K_1(n)$ for $n \ge 3$ with $ c_{\varrho} = (1 -\varrho )/\left( \sum_{n \ge 3} (n-1) K_1(n)\right)$. A look at the Figure~\ref{fig:alpha22} shows that, if $\gamma < 2.27\ldots$, the system is outside of the Cram\'er regime for $h$ is below a critical value $h_{c,\gamma}$. For larger values of $\gamma$ the system is outside of the Cram\'er regime for every $h>0$.

\section{The free model}
\label{sec:free} 

\begin{SCfigure}[10]
\centering
\caption{\label{fig:DandF} We plot the Large Deviations functional $D_h(1, \gamma)$ and the free energy
$\tilde \tf _\gamma (h)$, that coincides with $\tf_\gamma (h)$, as functions of $\gamma$ for a given $h>0$. 
Relevant features are the  convexity of the first and concavity of the second, which become strict in the 
Cramer region or interval $(1/\gamma_c(h), \gamma_c(h))$. Both functions are real analytic except at the boundary of the Cramer  interval
and they are affine functions outside of this interval. This and a number of other properties can be extracted from
the variational expression \eqref{eq:D} for $D_h(1, \gamma)$, the analysis in \S\ref{sec:constrained}  and formula
 \eqref{eq:Fgc} for $\tf_\gamma(h)$.
  In particular we have that $\partial_\gamma D_h(1, \gamma)$ is equal to $\tg(h)$ for $\gamma> \gamma_c(h)$ -- in \eqref{eq:gammah} for
  an explicit expression for $\gamma_c(h)$ -- and to $-\tg(h)/\gamma_c(h)+ D_h(1, \gamma_c(h))$ for  $\gamma<1/ \gamma_c(h)$. Note also that the fact that the lower endpoint of Cramer interval 
  is $1/\gamma_c(h)$   follows from the symmetry $D_h(v_1,v_2)=D_h(v_2,v_1)$ and 
  the scaling behavior $D_h(cv)=cD_h(v)$ for $c>0$ and $v=(v_1,v_2)$, for every $v_1$ and $v_2$ positive.
  The free energy becomes constant above $\gamma_c(h)$ and it is equal to
  $\tg(h)+ \vert \overline{\gl}_1(h)\vert$ (recall \eqref{eq:gl1bar}).
}
\leavevmode
\epsfxsize =7 cm
\psfragscanon
\psfrag{0}[c][l]{\small $0$}
\psfrag{1}[c][l]{\small $1$}
\psfrag{F}[c][l]{\small $\tilde\tf_\gamma(h)$}
\psfrag{D}[c][l]{\small $D_h(1,\gamma)$}
\psfrag{ga}[c][l]{\small $\gamma$}
\psfrag{gac}[c][l]{\small $\gamma_c(h)$}
\psfrag{gab}[c][l]{\small $\frac1{\gamma_c(h)}$}
\epsfbox{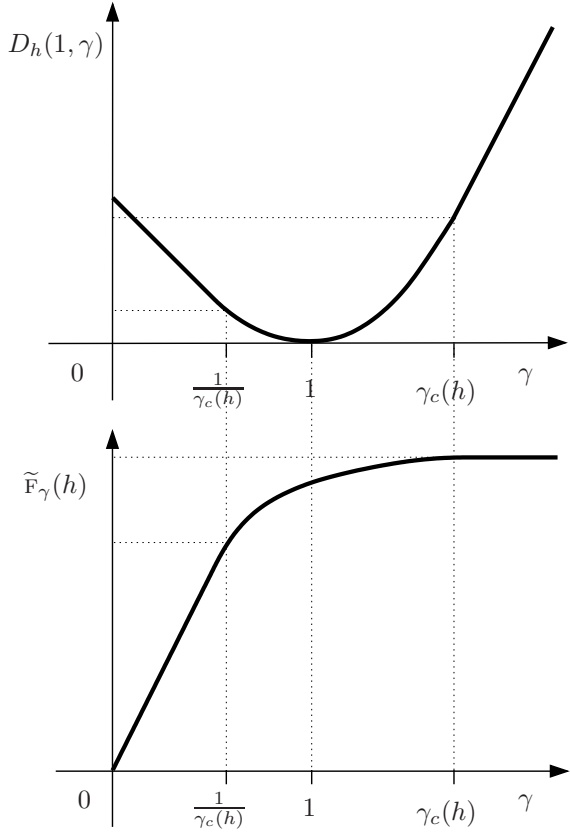}
\end{SCfigure}

In this section we use the notation $\gamma_N:= M/N$ and $\gamma'_N:= M'/N'$.
Recall that we assume $M\ge N$, but $M'\in \{0, \ldots, M\}$ may be smaller than $N'
\in \{0, \ldots, N\}$. We use the short-cut of saying that $\gamma$ is in the Cramer region
if $(1, \gamma)\in E_h$, i.e. if $(1, \gamma)$ is in the Cramer region.
Figure \ref{fig:DandF}  and its caption  sum up properties of $D_h(1, \cdot)$ and of $\tilde \tf_\cdot(h)$ that will come handy in the remainder.

\subsection{A preliminary result and free energy equivalence}
We start with a preliminary important lemma on the constrained free energy: 
a minimal part of its strength will be used to prove, just below (Proposition~\ref{th:f=c}),  that this free energy coincides with the one
of the free model. The full strength of this lemma is used in \S~\ref{sec:Zfsharp}.

\medskip

\begin{lemma}
\label{th:prelimZf}
There exists $a(h)> 0$ such that for 
$\gamma'_N \ge \gamma_N$ we have
\begin{equation}
\label{eq:prelimZf1}
\tilde\tf_{\gamma_N}(h) - \frac{N'} N \tilde\tf_{\gamma'_N}(h) \, \ge \, a(h)
\left(1 -\frac{\gamma_N}{\gamma'_N}\right)+ \frac{\gamma_N}{\gamma'_N} {\tilde\tf_{\gamma'_N}(h)} \left(1-\frac{M'}M\right)\, .
\end{equation}

For $\gamma'_N \le \gamma_N$ we have
\begin{equation}
\label{eq:prelimZf2}
\tilde \tf_{\gamma_N}(h) - \frac{N'} N \tilde \tf_{\gamma'_N}(h) \, \ge \, \partial_\gamma\tilde \tf_{\gamma}(h)\big \vert_{\gamma=\gamma_N}
(\gamma_N -\gamma'_N) + \tilde \tf_{\gamma'_N}(h) \left(1- \frac{N'}N\right)  \, ,
\end{equation}
where $\partial_\gamma\tilde \tf_{\gamma}(h)\ge 0 $ and   for every $\gep\in (0,  \gamma_c(h))$ we have  $\inf_{\gamma\le  \gamma_c(h)-\gep}\partial_\gamma\tilde \tf_{\gamma}(h)>0$.

In particular for every $\gep$ as above and every $L>0$  there exists $a_{\gep,L}(h)>0 $ such that 
\begin{equation}
\label{eq:prelimZf3}
\tilde \tf_{\gamma_N}(h) - \frac{N'} N \tilde \tf_{\gamma'_N}(h) \, \ge \, a_{\gep,L}(h)
\left(
\left\vert\gamma_N -\gamma'_N\right\vert  +
\ind_{\gamma'_N\ge \gamma_N} \gamma_N \left(1- \frac{M'}M\right)  
+
\ind_{\gamma'_N< \gamma_N}  \left(1- \frac{N'}N\right)  
\right)\, ,
\end{equation} 
for $\gamma_N \le \gamma_c(h)-\gep$ and $\gamma'_N \in [1/L, L]$.
\end{lemma}

\medskip

\noindent
{\it Proof.}
$\gamma'_N \ge \gamma_N$ we write
\begin{equation}
\label{eq:prelimZf1-1}
\tilde \tf_{\gamma_N}(h) - \frac{N'} N \tilde \tf_{\gamma'_N}(h) \, = \, \left( \tilde \tf_{\gamma_N}(h)- \frac{\gamma_N}{\gamma'_N}
\tilde \tf_{\gamma'_N}(h)\right)
+ \frac{\gamma_N}{\gamma'_N} {\tilde \tf_{\gamma'_N}(h)} \left(1-\frac{M'}M\right)\, ,
\end{equation}
and it suffices to focus on the first term between parentheses in the right-hand side and by using  $\gamma= \gamma_N$ and $\gamma'= \gamma'_N$ to keep expressions short, since $\gamma\mapsto \tilde\tf_\gamma (h)$ is concave 
%and strictly so for $\gamma< \gamma_c(h)$ 
we have $\tilde \tf_{\gamma'}(h)\le \tilde \tf_{\gamma}(h)+\partial_\gamma \tilde \tf_{\gamma}(h)(\gamma'-\gamma)$ so
\begin{equation}
\label{eq:int67.1}
\tilde \tf_{\gamma}(h)- \frac{\gamma}{\gamma'}
\tilde \tf_{\gamma'}(h)\, \ge \,
\left(1- \frac \gamma{\gamma'} \right) \left( \tilde \tf_{\gamma}(h)
-  \gamma \partial _\gamma \tilde \tf_{\gamma}(h)\right)\, ,
\end{equation}
and the right-hand side is non-negative because   $\gamma\mapsto \tilde \tf_{\gamma}(h)$ is concave
and $\tilde \tf_0(h)=0$, so 
 $\tilde \tf_{\gamma}(h)
-  \gamma \partial _\gamma \tilde \tf_{\gamma}(h)$ is non decreasing in $\gamma$ and it is even increasing if   $\gamma> 1/\gamma_c(h)$. This implies in particular that
\begin{equation}
a(h)\, :=\, \inf_{\gamma\ge 1}\left(
 \tilde \tf_{\gamma}(h)
-  \gamma \partial _\gamma \tilde \tf_{\gamma}(h)\right)\,=\, \tilde \tf_{1}(h)
-   \partial _\gamma \tilde \tf_{\gamma}(h)\big\vert_{\gamma=1}\, >\, 0\, ,
\end{equation}
and \eqref{eq:prelimZf1} is proven.

\medskip

Let us turn then to $\gamma'_N\le \gamma_N$ and  \eqref{eq:prelimZf2}.
This time we write
\begin{equation}
\label{eq:prelimZf2-1}
\tilde \tf_{\gamma_N}(h) - \frac{N'} N \tilde \tf_{\gamma'_N}(h) \, = \, \left( \tilde \tf_{\gamma_N}(h)- 
\tilde \tf_{\gamma'_N}(h)\right)
+ \left( 1- \frac{N'} N \right) {\tilde \tf_{\gamma'_N}(h)} \, ,
\end{equation}
so it suffices to bound from below the first term in parentheses in the right-hand side as claimed in \eqref{eq:prelimZf2},
but this is a direct consequence of concavity in $\gamma$ of the free energy.
%: 
%neglecting like before the subscript $N$ and using again the concavity we see that
%for $\gamma'\le \gamma$
%\begin{equation}
 %\tilde \tf_{\gamma_N}(h)- 
%\tilde \tf_{\gamma'_N}(h)\, \ge \, (\gamma-\gamma')\partial _\gamma  \tilde \tf_{\gamma_N}(h) \, ,
%\end{equation}

The validity of \eqref{eq:prelimZf3} follows from \eqref{eq:prelimZf1} and \eqref{eq:prelimZf2} by elementary considerations.
\qed

\medskip
We are now ready to show that free and constrained models have the same free energy: 

\medskip

\begin{proposition}
\label{th:f=c}
For every $h \in \bbR$ we have that the limit that defines $\tf_\gamma (h)$, i.e. \eqref{eq:F}, exists and 
\begin{equation}
\tf_\gamma (h) \, = \, \tilde{\tf}_\gamma(h) \,.
\end{equation}
\end{proposition}

\medskip

\noindent
{\it Proof.} Recall that we treat a result for $N\to \infty$ and $M \sim \gamma N$.
First of all observe that, since  $Z_{N, M, h}^f \ge Z_{N, M, h}^c$ we just need to worry about the upper bound. 
For the upper bound we use Proposition~\ref{th:Fgc} and, more precisely the uniform estimate \eqref{eq:Fgc-added}.
As a matter of fact the supremum with respect to $\gamma$ in  \eqref{eq:Fgc-added}
can be extended to $\gamma \in [1/L, L]$: the choice to restrict to $\gamma\ge 1$ in that proof was just to conform to the convention chosen from the beginning (but too restrictive for this proof as we will see). But the
extension to $\gamma  \in [1/L, L]$ can also be obtained from the result for $\gamma \in [1, L]$
by exploiting the symmetry of the expressions for the exchange of $M$ and $N$. Therefore we have that 
for every $\gd >0$
and every $L>1$ there exists $N_0$ such that 
\begin{equation}
\label{eq:forf=c.1}
Z_{N', M', h}^c \, \le \, \exp \left(N' \tilde \tf_{\gamma'_N}(h)+ \gd
\right)\, ,
\end{equation}
for every $N'$ and $M'$ such that $N'\ge N_0$ and $\gamma'_N=M'/N' \in [1/L, L]$.
Or, equivalently (possibly changing the value of $N_0$), \eqref{eq:forf=c.1}
holds for every $M'\wedge N'\ge N_0$ and $\gamma'_N \in [1/L, L]$. 
But then thanks to Lemma~\ref{th:prelimZf}
\begin{equation}
\tilde \tf_{\gamma_N}(h) - \frac{N'} N \tilde \tf_{\gamma'_N}(h) \, \ge \, 0\, ,
\end{equation}
so 
\begin{equation}
\label{eq:forf=c.2}
Z_{N', M', h}^c \, \le \, \exp \left(N \tilde \tf_{\gamma_N}(h)+ \gd
\right)\, ,
\end{equation}
for the same values of $M'$ and $N'$.
Therefore by using  $K_f(n) \le n^C$ (for some $C>0$) we have 
\begin{multline}
\label{eq:pwf7}
Z^f_{N,M,h} \, \le \, N^{2C}\sum_{N', M'} Z_{N', M', h}^c \, \le (1+\gamma)N^{2C+2} \exp \left(N \tilde \tf_{\gamma_N}(h)+ \gd
\right)
+\\ N^{2C} N_0^2 \max_{N',M'\le N_0}  Z_{N', M', h}^c 
+ N^{2C} \sumtwo{N',M': \, N' \wedge M' \ge N_0  }{N'\le N, M'\le M, \gamma'_N\not\in [1/L,L]} Z_{N', M', h}^c\, .
\end{multline}
Since the maximum in the second line is just a constant we are left with controlling the last sum.
But, again, for every $L'>0$ -- we want $L'<L$ but close to $L$ so we set $L'=2L/3$ -- we have  $Z^c_{N,\lfloor NL'\rfloor,h}\le \exp(N \tilde \tf_{L'}(h)+ \gd)$ for $N$
sufficiently large and we can use this result to bound the growth of $Z^c_{N',M',h}$, with 
$\gamma_N' \ge L$ by the usual comparison argument obtained by modifying the last inter-arrival.
The net result is that for every $L'<L$  (we choose $L-L'$ small) and every $\gd>0$ we can find $N_0$ such that for 
$N'\ge N_0$
\begin{equation}
\label{eq:pwf8}
Z^c_{N',M',h} \stackrel{\gamma'_N \ge L} \le   \exp\left(N' \left(\tilde \tf_{L'}(h)+ \gd\right)\right)\,.
\end{equation}
But $\gamma'_N \ge L$ implies $N'\le M/L \le (3/2)\gamma N/L$ (here $3/2$ can be replaced by any number larger than one), so 
$N' \tilde \tf_{L'}(h) \le (3\gamma/2) N  (\tilde \tf_{L'}(h)/L')$.
Since $\tilde \tf_\gamma (h)= \tilde \tf_{\gamma_c(h)} (h)$ for $\gamma\ge \gamma_c(h)$
we see that for $2L/3 \ge \gamma_c(h)$
\begin{equation}
N' \tilde \tf_{L'}(h) \, \le\,  \left( \frac 32 \gamma   \tilde \tf_{\gamma_c(h)} (h)\right) \frac {N}{ L'}
 \, \le\,  \left(  \gamma   \tilde \tf_{\gamma_c(h)} (h)\right) \frac {N}{ L}
\, ,
\end{equation}
and we can take $L' < L$ and $N' \le M/L \le c \gamma N/L$ with $c>1$.
so by choosing $L$ sufficiently large we can bound the expression in \eqref{eq:pwf8}
by $\exp(2\gd N)$. Finally, by a symmetry argument, we see that the contribution of the terms with
$\gamma'_N \le 1/L$  in the last term in \eqref{eq:pwf7} is smaller than the contribution of the terms
with $\gamma'_N \ge L$, and we are done. 
\qed

\subsection{Sharp estimates on $Z^f_{N, M, h}$}
\label{sec:Zfsharp}

\subsubsection{The case $h>0$: Proof of (1) in Theorem~\ref{th:Zsharp}}

Recall the definition \eqref{eq:Zf} of $Z_{N, M,h}^f$ and that we work with $M \sim \gamma N $, $\gamma$ in the Cramer region.
For every $a\in (0,1)$
\begin{multline}
\label{eq:mlin1}
\sum_{i=\lfloor a N\rfloor}^{N} \sum_{j=0}^M K_f(i) K_f(j) Z_{N-i, M-j,h}^c \, =\\
\sum_{i=\lfloor a N\rfloor}^{N} \sum_{j=0}^M K_f(i-\lfloor a N\rfloor)\frac{K_f(i)}{K_f(i-\lfloor a N\rfloor)} K_f(j) Z_{N-i, M-j,h}^c\,, 
\end{multline}
and since $K_f(\cdot)$ is positive and it has an asymptotic power law behavior there exists $C>0$ such that (for $N$ sufficiently large)
\begin{equation}
\label{eq:mlin2}
\begin{split}
\sum_{i=\lfloor a N\rfloor}^{N} \sum_{j=0}^M K_f(i) K_f(j) Z_{N-i, M-j,h}^c \,&\le \,  N^C
\sum_{i=0}^{N-\lfloor a N\rfloor} \sum_{j=0}^M  K_f(i) K_f(j) Z_{N- \lfloor a N\rfloor-i, M-j,h}^c
\\
&=\, N^{C} Z_{N-\lfloor a N\rfloor, M,h}^f\, .
\end{split}
\end{equation}
Note that the leftmost term in \eqref{eq:mlin1} and \eqref{eq:mlin2} is $Z_{N, M,h}^f$ if $a=0$. On the other hand 
\begin{equation}
\label{eq:obs4.0}
N^{C} Z_{N-\lfloor a N\rfloor, M,h}^f \stackrel{N \to \infty}
\asymp \exp\left(N(1-a)\tf_{\gamma_N/(1-a)}(h)\right)\, ,
\end{equation}
(we say that $f(x)\stackrel{x \to a} \asymp g(x)$ if $f(x) = O(g(x))$ and $g(x) = O(f(x))$ as $x \to a$) and we observe that 
\begin{equation}
\label{eq:obs4.1}
\tf_{\gamma_N} (h)\, >\, (1-a) \tf_{\gamma_N/(1-a)}\,,
\end{equation}
which follows from 
\eqref{eq:prelimZf1} by choosing $M'=M$ so  $N'/N= (\gamma_N/\gamma'_N)$ so that one obtains
$\tf_\gamma(h)-(\gamma/\gamma') \tf_{\gamma'}(h)>0$ for every $\gamma'>\gamma\ge 1$, and
this inequality becomes \eqref{eq:obs4.1} if we choose $\gamma'=\gamma_N/(1-a)$.

At this point we observe that,
since $K_f(0)=1$,  we have 
$Z_{N, M,h}^f\ge Z_{N, M,h}^c \asymp \exp(N \tf_{\gamma_N}(h))$ and therefore, by \eqref{eq:mlin2}--\eqref{eq:obs4.1} we see 
that for every $a\in (0,1)$
  there exists $q>0$ such that
\begin{equation}
\label{eq:Zfq0}
Z_{N, M,h}^f\,
= \, 
\left( 1+ O\left( \exp(- q N ) \right)\right)
\sum_{i=0}^{\lfloor a N\rfloor} \sum_{j=0}^M K_f(i) K_f(j) Z_{N-i, M-j,h}^c \, .
\end{equation}

A parallel, somewhat easier, argument can be put at work 
when we restrict the summation in the definition \eqref{eq:prelimZf1} of $Z_{N, M, h}^f$ to $j \ge \lfloor a N \rfloor$.
We have to 
 use again Lemma~\ref{th:prelimZf}: \eqref{eq:prelimZf1} for $N'=N$ simply becomes the fact that $\gamma \mapsto\tf_{\gamma}(h)$
 is (strictly) increasing for $\gamma < \gamma_c(h)$ and this allows to conclude
that  for every $a\in (0,1)$ there exists $q>0$ such that
 \begin{equation}
 \label{eq:Zfq}
Z_{N, M,h}^f\,
= \, 
\left( 1+ O\left( \exp(- q N ) \right)\right)
\sum_{i=0}^{\lfloor a N\rfloor} \sum_{j=0}^{\lfloor  aN\rfloor} K_f(i) K_f(j) Z_{N-i, M-j,h}^c \, .
\end{equation}

\medskip 

With \eqref{eq:Zfq0} we see that we can restrict  the sum in \eqref{eq:Zf} to a small (since we can choose $a>0$ small)
macroscopic square. We want now to show that we can  restrict almost to a microscopic square: a microscopic square would be
a square of size that does not diverge with $N$. The result we shall now prove is 
 \begin{equation}
 \label{eq:Zflog}
Z_{N, M,h}^f\,
= \, 
\left( 1+ O\left( \exp(- (\log N)^{3/2} ) \right)\right)
\sum_{i=0}^{\ell_N} \sum_{j=0}^{\ell_N} K_f(i) K_f(j) Z_{N-i, M-j,h}^c \, ,
\end{equation} 
where
\begin{equation}
\label{eq:ellN}
\ell_N \, :=\, 
\lfloor  (\log N)^2 \rfloor\, .
\end{equation}
For this choose $a$ small so that $(M-j) / (N-i)$ is in the Cramer region for all values
of $i$ and $j$ in the summation in \eqref{eq:Zfq}. We can then apply
\eqref{eq:Zc-sharp} and, more precisely, the following consequence of   \eqref{eq:Zc-sharp}: for every $N$ sufficiently large
\begin{equation}
\label{eq:fromZc-sharp1}
Z_{N', M', h}^c \, \le \, 
C
\exp \left(N' \tf_{\gamma'_N}(h)\, ,
\right) \, , 
\end{equation}
where $C>0$  and $N-\lfloor  aN\rfloor \le N'
\le N$, $M-\lfloor  aN\rfloor\le M' \le M= \gamma_N N $. 
For this we exploit \eqref{eq:prelimZf3} of Lemma~\ref{th:prelimZf}: since $i$ and $j$ are (macroscopically)
small we have that there exists $c=c(h,a)>0$ such that
\begin{equation}
\label{eq:prelimZf1s}
N\tf_{\gamma_N}(h) - N'  \tf_{\gamma'_N}(h) \, \ge \, c 
\left( N
\left\vert \gamma'_N -\gamma_N\right \vert +  
\ind_{\gamma'_N \ge \gamma_N}\left(M-M'\right)
+\ind_{\gamma'_N < \gamma_N}\left(N-N'\right)\right)\, ,
\end{equation}
where $N'=N-i$, $M'=M-j$,  $N$ is sufficiently large
and both $i$ and $j$ in their range of summation. But now we recall that we aim at \eqref{eq:Zflog}
and therefore if $i< \ell_N$
then $j \ge \ell_N$ and the same is true if we exchange $i$ and $j$. Therefore, omitting the constant $c$,
the right-hand side of  \eqref{eq:prelimZf1s} is equal to $\gamma'_N (N-N')\ge  (N-N') \ge \ell_N/(2\gamma)$ for $\gamma'_N \ge \gamma_N$. For $\gamma'_N \le \gamma_N$ instead 
the right-hand side of  \eqref{eq:prelimZf1s} is equal to $M-N\gamma'_N+N-N'$
which, on one hand, it bounded below by
$M-N\gamma_N+N-N' = N-N'$. On the other hand it is equal to $(M-M') -(N-N') (\gamma_N' -1)$
which is bounded below by $(M-M') -(N-N') (\gamma_N -1)$. For $\gamma'_N \le \gamma_N$ we have
$M-M'\ge  \ell_N$ so either $N-N' \le   \frac 1{2 (\gamma_N-1)}\ell_N$, so
 $(M-M') -(N-N') (\gamma_N -1)\ge \frac 12 \ell_N$,
or $N-N' >   \frac 1{2 (\gamma_N-1)}\ell_N$. Hence 
the  right-hand side of  \eqref{eq:prelimZf1s} is bounded below for $\gamma_N'\le \gamma_N$ by
\begin{equation}
\frac c2 \ell_N \min\left(
  \frac 1{ \gamma_N-1}, 1\right)\, ,
\end{equation}
so, recalling the lower bound found for $\gamma_N'\ge \gamma_N$, we see that if we set   $c_\gamma := \frac 12  \min(1/\gamma,1)$ we have 
\begin{equation}
\label{eq:prelimZf1r}
N\tf_{\gamma_N}(h) - N'  \tf_{\gamma'_N}(h) \, \ge \, c_{\gamma_N} \ell_N\, ,
\end{equation}
for $((N-N'), M-M'))\in ( [0, aN\rfloor)^2\setminus [0, \ell_N)^2)\cap \bbZ^2$.
But then \eqref{eq:Zflog} becomes evident from \eqref{eq:fromZc-sharp1}, \eqref{eq:prelimZf1r}
and the fact 
 the summation is on less than $a^2 N^2=O(N^2)$ sites: so the total contribution by summing over the sites in small macroscopic square minus the almost microscopic square is $O(N^2\exp(N\tf_{\gamma_N}(h)- c_{\gamma_N} \ell_N))$.
 On the other hand $Z_{N, M,h}^f\ge Z_{N, M,h}^c\ge AN^{-1/2}\exp(N\tf_{\gamma_N}(h))$ for 
 some $A>0$ and $N$ large, cf. \eqref{eq:Zc-sharp}, and \eqref{eq:Zflog} is proven.

\medskip

The question of the sharp estimates on $Z^f$ is then reduced to find the leading behavior of 
\begin{equation}
\sum_{i=0}^{\ell_N} \sum_{j=0}^{\ell_N} K_f(i) K_f(j) Z_{N-i, M-j,h}^c \, .
\end{equation}
Observe  that $(N-i,M-j)/ \vert (N-i,M-j)\vert$ is close to $(1, \gamma)$ 
and hence it is in a compact subset of $J$ of $E_h$ for every $ i,j \in [0,\ell _N]$. So 
%we have $(N-i,M-j)\in E_h$, then by \eqref{eq:Zc-sharp}
%\begin{equation}\label{eq:sharpZ_N-i}
%Z^{c}_{N-i, M-j,h} = \frac{A((M-j)/(N-i)) + \tilde \epsilon (N-i,M-j)}{\sqrt{N-i}} \exp \left( (N-i) \tf_{(M-j)/(N-i)} (h) \right) \,.
%\end{equation}
%Therefore
we have %($\ell_N$ defined in \eqref{eq:ellN})
\begin{multline}
\label{eq:Zflog1}
 \sum_{i=0}^{\ell_N} \sum_{j=0}^{\ell_N}  K_f(i) K_f(j) Z^{c}_{N-i,M-j,h} 
\stackrel{N \to \infty}\sim  \frac{A(\gamma) \exp \left( \tf_{\gamma_N}(h) N \right)}{\sqrt{N}} \\
 \sum_{i=0}^{\ell_N} \sum_{j=0}^{\ell_N}  K_f(i) K_f(j) \exp \left( - N \left( \tf_{\gamma_N}(h) - \tf_{(M-j)/(N-i)} (h) + \frac{i}{N} \tf_{(M-j)/(N-i)} (h) \right) \right) \,.
\end{multline}
Taylor expansion yields that there exists $\tilde\gamma\in J$
\begin{multline}
\label{eq:Taylor1}
\tf_{\gamma_N}(h) - \tf_{(M-j)/(N-i)} (h) \,=\\ \left( \gamma_N - \frac{M-j}{N-i} \right) \partial_\gamma \tf_{\gamma}(h)\big \vert_{\gamma=(M-j)/(N-i)} + \left( \gamma_N - \frac{M-j}{N-i}\right)^2 \partial^2_\gamma \tf_{\gamma}(h)\big\vert_{\gamma=\tilde\gamma}\,, 
\end{multline}
and since
\begin{equation}
\label{eq:ga-ga0}
\gamma_N - \frac{M-j}{N-i}  \, =\, 
% \gamma \left( 1 - \frac{1-\frac{j}{M}}{1-\frac{i}{N}} \right) = \gamma_N \left( 1 - (1 - \frac{j}{M}) ( 1 + \frac{i}{N}) \right) + O\left(\frac{\ell_N^2}{N^2}\right)=
 \frac{j}{N} - \gamma_N \frac{i}{N} + O\left(\frac{\ell_N^2}{N^2}\right)\, = \, O\left(\frac{\ell_N}{N}\right)  \,,
\end{equation}
 the second term in the right-hand side of \eqref{eq:Taylor1} is $O\left({\ell_N^2}/{N^2}\right)$ and finally,
using also that $\gamma_N-\gamma=O(1/N)$, we get to
\begin{multline}
\label{eq:Zfe}
\tf_{\gamma_N}(h) - \tf_{(M-j)/(N-i)} (h) + \frac{i}{N} \tf_{(M-j)/(N-i)} (h) 
= \\
\frac{1}{N} \left( j \partial_{\gamma} \tf_{\gamma}(h) + i  \left( \tf_\gamma(h) - \gamma \partial_{\gamma} \tf_{\gamma}(h) \right) \right)  + O\left(\frac{\ell_N^2}{N^2}\right)\,.
\end{multline}
%where we have used the fact that for fixed h, $\partial_\gamma \tf_{\gamma}(h)\big \vert_{\gamma=(M-j)/(N-i)} = \partial_\gamma \tf_{\gamma}(h)$. 
Recall now that $ \tf_\gamma(h) - \gamma \partial_{\gamma} \tf_{\gamma}(h)$ is positive for $\gamma>\gamma_c(h)$, so that 
 \eqref{eq:Zfe} implies  the double sum in the right-hand side of \eqref{eq:Zflog1} 
 converges to 
\begin{equation}
\label{eq:Cgah}
C_{\gamma, h} \, :=\,  \left( \sum_{i=0}^\infty K_f(i) \exp \left(- i  \left( \tf_\gamma(h) - \gamma \partial_{\gamma} \tf_{\gamma}(h) \right) \right) \right) \left( \sum_{j=0}^\infty K_f(j) \exp \left( -j \partial_{\gamma} \tf_{\gamma}(h) \right) \right)\, ,
\end{equation}
and
\begin{equation}
\label{eq:sharpZf}
Z_{N, M,h}^f\,
\sim \, 
A(\gamma) C_{\gamma,h}   \frac{\exp \left( \tf_{\gamma_N}(h) N \right)}{\sqrt{N}} \,.
\end{equation}
This completes the proof of $(1)$ in Theorem~\ref{th:Zsharp} with $ c_{\gamma,h} = A(\gamma) C_{\gamma,h}$.
\qed

\medskip

\subsubsection{The case $h < 0$: Proof of (2) in Theorem~\ref{th:Zsharp}}

In the delocalized phase, $ \tilde \tau_h$ is terminating and $ Z^c_{N,M,h} = \bP\left( (N,M) \in \tilde \tau_h \right)$.
Recall that $\tau$ is recurrent, i.e. $\sum_{n,m} \mathrm{K}(n,m)=1$, then
\begin{equation}
\tilde K_h(\infty)  \, = \, 1 - \sum_{(n, m)\in \bbN^2} \tilde {\mathrm{K}}_{h}(n,m)= 1 - \exp(h) > 0 \,.
\end{equation} 
For the renewal function we write
\begin{equation}
\label{eq:GF}
\bP \left( (N,M) \in \tilde \tau_h \right) %= \sum_{j=0}^{\infty} \bP^{j*} \left( (\tilde \tau_h)_1  (N,M) \right) 
\,=\,  \sum_{j=0}^{\infty} \tilde {K}_{h}^{j*}(N,M) = \sum_{j=0}^\infty \exp(jh) {\mathrm{K}}^{j*}(N,M) \, ,
\end{equation}
where $\tilde K_h^{j*} (\cdot, \cdot)$ is the j-fold convolution of $ \bP(.)$.

\medskip
Take $\ga >1$, we will use the following estimate:

\medskip

\begin{lemma}
\label{th:Q}
There exists $c >0$ such that
for every $j \in \bbN$ and for every $(N,M)$
\begin{equation}
\label{eq:ind}
\mathrm{K}^{j*} (N,M) \le j^c \, \mathrm{K}(N,M) \,.
\end{equation}
Moreover,
\begin{equation}
\label{eq:Q1}
\lim_{N,M \to \infty} \frac{\mathrm{K}^{j*}(N,M)}{\mathrm{K}(N,M)} \,=\, j \,.
\end{equation}
\end{lemma}

\medskip

\noindent
{\it Proof.}
It is clear that \eqref{eq:ind} holds for $j=1$. Assume that it is true for $j< 2s$ and we want to show it for $j=2s$. Observe that
\begin{multline}
\mathrm{K}^{2s*} (N,M) \le 2 \sum_{n=1}^{\lfloor N/2 \rfloor} \sum_{m=1}^{\lfloor M/2 \rfloor} \mathrm{K}^{s*} (n,m) \mathrm{K}^{s*}(N-n,M-m) + \\
 \left( \sum_{n=1}^{\lfloor N/2 \rfloor} \sum_{m=\lfloor M/2 \rfloor+1}^{M} + \sum_{n=\lfloor N/2 \rfloor+1}^{N} \sum_{m=1}^{\lfloor M/2 \rfloor} \right) \mathrm{K}^{s*} (n,m) \mathrm{K}^{s*}(N-n,M-m):= Q_1 + (Q_2+Q_3) \,.
\end{multline}
First we have
\begin{equation}
Q_1 \le 2 s^c  \sum_{n=1}^{\lfloor N/2 \rfloor} \sum_{m=1}^{\lfloor M/2 \rfloor} \mathrm{K}^{s*} (n,m) \mathrm{K}(N-n,M-m) \,,
\end{equation}
and from \eqref{eq:psv}, we see that there exists $c_1 >0$ such that $L(ux) \le c_1 L(x)$ for every $u \in [ 1/2,1]$ and $x \ge 1$, therefore $\mathrm{K}(N-n,M-m) \le c_1 2^{1+\ga} \mathrm{K}(N,M)$ and
\begin{equation}
Q_1 \le c_1 2^{2+\ga - c} (2s)^c \mathrm{K}(N,M)  \sum_{n,m} \mathrm{K}^{s*} (n,m) \,.
\end{equation}

Now observe that $Q_2(N, M)=Q_3(M,N)$, so in view of the bound we are after 
it suffices to consider $Q_2$. Remark then that 
there exists $c_2 >0$ such that
\begin{equation}
\label{eq:Q2}
Q_2 \,\le\,   s^c  \sum_{n=1}^{\lfloor N/2 \rfloor} \sum_{m=\lfloor M/2 \rfloor+1}^{M} \mathrm{K}(n,m) \mathrm{K}^{s*}(N-n,M-m) \le c_2 2^{- c} (2 s)^c K(M) \,,
\end{equation}
and 
\begin{equation}
\label{eq:Q22}
Q_3 \, \le\,   s^c  \sum_{n=1}^{\lfloor N/2 \rfloor} \sum_{m=\lfloor M/2 \rfloor+1}^{M} \mathrm{K}(N-n,M-m) \mathrm{K}^{s*}(n,m)
 \le c_2 2^{- c} (2 s)^c K(N) \,.
\end{equation}
It suffices then to prove that there exists $c_3 >0$ such that $ K(N) \wedge K(M) \le c_3 \mathrm{K}(N,M)$. By elementary arguments we see that this follows if we can show that for every $x,y \ge 0$ and for every slowly varying function $L(\cdot)$, there exists $ c_{L(\cdot)} >0$ such that
\begin{equation}
\label{eq:stat4}
x L(x) \vee y L(y) \, \ge\, c_{L(\cdot)} (x+y) L(x+y) \,.
\end{equation}
By symmetry it suffices to consider the case $y\ge x$ and in this case it suffices to show that 
$yL(y) \ge c_{L(\cdot)} (x+y) L(x+y)$. Since $(x+y)\in [y, 2y]$ we can apply \eqref{eq:psv} to see that there exists
$c'_{L(\cdot)}>0$ such that $L(y)\ge c'_{L(\cdot)}  L(x+y)$ for every $y\ge x$ and this implies the desired inequality, and therefore \eqref{eq:stat4},
with $c_{L(\cdot)} =c'_{L(\cdot)}/2$.

Therefore $ Q_2 + Q_3 \le c_4 2^{1- c} (2 s)^c \mathrm{K}(N,M)$ and 
\begin{equation}
\mathrm{K}^{2s*} (N,M) \le ( c_1 2^{2+\ga-c} + c_4 2^{1-c}) (2 s)^c \mathrm{K}(N,M) \,,
\end{equation}
and if $c = 1+ \log_2 (c_1 2^{1+\ga} +c_4)$, we obtain \eqref{eq:ind} for $j=2s$. The procedure 
can be repeated   for $j=2 s +1$ with minor  changes. Therefore \eqref{eq:ind} is proven.
\medskip

For what concerns  \eqref{eq:Q1}, serve that it  holds for $j=1$. Assume that it is still valid to $j=s$ and write 
\begin{equation}
\label{eq:dsum}
\frac{\mathrm{K}^{(s+1)*}(N,M)}{\mathrm{K}(N,M)} = \sum_{n=1}^{N-1}  \sum_{m=1}^{M-1} \frac{\mathrm{K}^{s*}(n,m) \mathrm{K}(N-n,M-m)}{\mathrm{K}(N,M)}\,.
\end{equation}
Split the double sum in \eqref{eq:dsum} to four terms:
\begin{equation}
S_1 = \sum_{n=1}^{\lfloor N/2 \rfloor}  \sum_{m=1}^{\lfloor M/2 \rfloor} \mathrm{K}^{s*}(n,m) \frac{\mathrm{K}(N-n,M-m)}{\mathrm{K}(N,M)}
\end{equation}
\begin{equation}
S_2 =  \sum_{n=1}^{\lfloor N/2 \rfloor -1}  \sum_{m=1}^{\lfloor M/2 \rfloor-1} \mathrm{K}(n,m) \frac{\mathrm{K}^{s*}(N-n,M-m)}{\mathrm{K}(N,M)}
\end{equation}
\begin{equation}
S_3 =  \sum_{n=1}^{\lfloor N/2 \rfloor }  \sum_{m=1}^{\lfloor M/2 \rfloor-1} \mathrm{K}^{s*}(n,M-m) \frac{ \mathrm{K}(N-n,m)}{\mathrm{K}(N,M)}
\end{equation}
\begin{equation}
S_4 =  \sum_{n=1}^{\lfloor N/2 \rfloor -1}  \sum_{m=1}^{\lfloor M/2 \rfloor} \mathrm{K}^{s*}(N-n,m) \frac{\mathrm{K}(n,M-m)}{\mathrm{K}(N,M)}
\end{equation}
For any fixed $n$ and $m$ and as $ N,M \to \infty$, the two ratios in $S_1$ and $S_2$ converge respectively to 1 and to s (recall that $\mathrm{K}(n,m)=K(n+m)$) and the two ratios are uniformly bounded (from the uniform convergence property of the slowly varying functions \eqref{eq:psv} for the ratio in $S_1$ and from \eqref{eq:ind} for the ratio in $S_2$). Then by  (DOM), we obtain $S_1+ S_2 \to 1+s$ as $N,M \to \infty$. 

Since $S_3$ and $S_4$ are essentially the same quantity when we exchange $M$ and $N$, we just focus on $S_3$.
If we first assume that $M\ge N$ (hence $M+N \in [M, 2M]$), 
from \eqref{eq:ind} and \eqref{eq:psv}  we obtain that
\begin{multline}
\label{eq:S3}
S_3 \,\le\,  
s^c \sum_{n=1}^{\lfloor N/2 \rfloor }  \sum_{m=1}^{\lfloor M/2 \rfloor-1} \mathrm{K}(n,M-m) \frac{ \mathrm{K}(N-n,m)}{\mathrm{K}(N,M)} \, \le\\ c_5 s^c 
\sum_{n=1}^{\lfloor N/2 \rfloor }  \sum_{m=1}^{\lfloor M/2 \rfloor-1}\frac{L(M+n)}{(M+n)^{1+\ga}}
\frac{L(N+m)}{(N+m)^{1+\ga}}
\frac{(N+M)^{1+\ga}}{L(N+M)}\, 
\le \\
c_6 s^c 
\sum_{n=1}^{\lfloor N/2 \rfloor }  \sum_{m=1}^{\lfloor M/2 \rfloor-1}
\frac{L(N+m)}{(N+m)^{1+\ga}}
\, \le \, 
 \frac{c_6 s^c}2 N \sum_{n>N}\frac{L(n)}{n^{1+\ga}}\, \le \, c_7 L(N)N^{1-\ga}\, .
\end{multline}
By repeating the argument for $N \ge M$ we obtain that $S_3 =O(L(M)M^{1-\ga})$ in this case. 
Therefore, since $\ga>1$, $\lim_{N, M \to \infty}S_3=0$  by the basic properties of 
slowly varying functions and an elementary argument. 
\qed

\medskip

To prove the part $(2)$ in theorem~\ref{th:Zsharp}, we need to know the sharp estimates in the constrained case for $h<0$:
\begin{proposition}
\label{th:SharpZ}
If $h<0$, for every $(N,M)$, there exists $c_h >0$ such that 
\begin{equation}
\label{eq:sharpZ_c1}
Z^c_{N,M,h} \le c_h \mathrm{K}(N,M) \,
\end{equation}
where $c_h = \sum_{j=0}^{\infty} j^c \exp(j h)$. Moreover
\begin{equation}
\label{eq:sharpZ_c}
Z^c_{N,M,h} \stackrel {N,M \to \infty}\sim \frac{\exp(h)}{{\left( 1 - \exp(h) \right)}^{2}} \, \mathrm{K}(N,M)\,.
\end{equation} 
\end{proposition}

{\it Proof.}
From \eqref{eq:GF}, we have
\begin{equation}
\label{eq:pf}
\frac{Z^c_{N,M,h}}{\mathrm{K}(N,M)} = \frac{\bP \left( (N,M) \in \tilde \tau_h \right)}{\mathrm{K}(N,M)} = \sum_{j=0}^{\infty} \exp(j h) \frac{\mathrm{K}^{j*}(N,M)}{\mathrm{K}(N,M)} \,,
\end{equation}
and using \eqref{eq:ind}, for fixed $j$, we see that the ratio is bounded above by $j^c$, therefore we obtain \eqref{eq:sharpZ_c1}.

For \eqref{eq:sharpZ_c}, the ratio in \eqref{eq:pf} converges to $j$ as $N,M \to \infty$ from \eqref{eq:Q1} and bounded from \eqref{eq:ind}, then by  (DOM) we get
\begin{equation}
\lim_{N,M \to \infty} \frac{Z^c_{N,M,h}}{\mathrm{K}(N,M)}  =  \sum_{j=0}^{\infty} j \exp(j h)= \frac{\exp(h)}{{\left( 1 - \exp(h) \right)}^{2}} \,,
\end{equation}
\qed

\medskip

We are now ready to prove the sharp estimate of $Z^f_{N,M,h}$:
\begin{proposition}
\label{th:sharpZf}
Suppose that $M \sim \gamma N$. 
For $h <0$, as $N \to \infty$
\begin{itemize}
\item If $ \overline{\alpha}< (1+\alpha)/2$, we have
\begin{equation}
\label{eq:sharpZ_f}
Z^f_{N,M,h} \, \sim \,  \frac{K_f(N) K_f(M)}{1- \exp(h)} \,.
\end{equation}
\\
\item If $ \overline{\alpha} > (1+\alpha)/2$, we have
\begin{equation}
\label{eq:sharpZ_f1}
Z^f_{N,M,h} \, \sim \,  \frac{ \exp(h){\left( \sum_{n \ge 0} K_f(n) \right)}^2 }{(1-\exp(h))^2} \mathrm{K}(N,M) \,.
\end{equation}
\end{itemize}
\end{proposition}

\medskip

\noindent
{\it Proof.}
Let us write
\begin{equation}
\label{eq:ff}
\frac{Z^f_{N,M,h}}{K_f(N)K_f(M)} = \sum_{n=0}^N \sum_{m=0}^M Z^c_{n,m,h} \frac{K_f(N-n)K_f(M-m)}{K_f(N)K_f(M)} \,.
\end{equation} 
We split the last sum into 
\begin{multline}
\label{eq:T4}
T_1+T_2+T_3+T_4 = \\
\left( \sum_{n=0}^{\lfloor N/2 \rfloor} \sum_{m=0}^{\lfloor M/2 \rfloor}
+ \sum_{n={\lfloor N/2 \rfloor} +1}^N \sum_{m={\lfloor M/2 \rfloor}+1}^M 
+ \sum_{n=0}^{\lfloor N/2 \rfloor} \sum_{m={\lfloor M/2 \rfloor}+1}^M  
+ \sum_{n={\lfloor N/2 \rfloor} +1}^N \sum_{m=0}^{\lfloor M/2 \rfloor} \right) \\
Z^c_{n,m,h} \frac{K_f(N-n)K_f(M-m)}{K_f(N)K_f(M)} \,.
\end{multline}
For $T_1$, for fixed $n$ and $m$, the ratio in \eqref{eq:T4} converges to $1$ and by  \eqref{eq:psv}, this ratio is bounded. Then by  (DOM), Fubini-Tonelli Theorem and the fact that $\mathrm{K}(\cdot,\cdot)$ is a discrete probability density we obtain (recall \eqref{eq:Zc-2} and in this case $\tg(h)=0$)
\begin{multline}
\label{eq:sumZ}
\lim_{N,M \to \infty} T_1 = \sum_{n=0}^{\infty} \sum_{m=0}^{\infty} Z^c_{n,m,h} = \sum_{n=0}^{\infty} \sum_{m=0}^{\infty} \sum_{j=0}^{\infty} \tilde {\mathrm{K}}_h^{j*}(n,m)=\\
\sum_{n=0}^{\infty} \sum_{m=0}^{\infty} \sum_{j=0}^{\infty} \exp(j h)  \mathrm{K}^{j*}(n,m)
= \frac{1}{1- \exp(h)} \,.
\end{multline}

For $T_2$, if $ \overline{\alpha} \in  [1, (1+\alpha)/2)$ and  $\sum_{n \ge 0} K_f(n) < \infty$,  using \eqref{eq:sharpZ_c}, we obtain (recall that $M \sim \gamma N$)
\begin{equation}
T_2  \le c_8 \sum_{n= \lfloor N/2 \rfloor+1}^N \sum_{m= \lfloor M/2 \rfloor+1}^M \mathrm{K}(n,m) \frac{K_f(N-n)K_f(M-m)}{K_f(N)K_f(M)}  \, =\, O\left( N^{2 \overline{\alpha} - \alpha - 1} \frac{L(N)}{(\overline L  (N))^2}\right)\, ,
\end{equation}
 so  $ T_2  \to 0$ as $N \to \infty$.

If $\overline{\alpha} \le 1$ and $ \sum_{n \ge 0} K_f(n) = \infty$, using once again \eqref{eq:sharpZ_c}  with \eqref{eq:recall}, we get
\begin{multline}
T_2  \,\le \,c_8 \sum_{n= \lfloor N/2 \rfloor+1}^N \sum_{m= \lfloor M/2 \rfloor+1}^M \mathrm{K}(n,m) \frac{K_f(N-n)K_f(M-m)}{K_f(N)K_f(M)}\, \le   \\
c_9 \mathrm{K}(N,M)   \frac{\sum_{n= \lfloor N/2 \rfloor+1}^N K_f(N-n)
\sum_{m= \lfloor M/2 \rfloor+1}^MK_f(M-m)}{K_f(N)K_f(M)}
 \,=\, O \left(N^{1-\ga}L(N) \right),
\end{multline}
and  $ T_2  \to 0$ as $N \to \infty$ in this case too.

Let us look at $T_3$ (the argument for $T_4$ is identical). We have 
\begin{multline}
\label{eq:T3.g}
T_3 \,=\,  \sum_{n=0}^{\lfloor N/2 \rfloor} \sum_{m=\lfloor M/2 \rfloor +1}^{M} Z^c_{n,m,h} \frac{K_f(N-n)}{K_f(N)} \frac{K_f(M-m)}{K_f(M)} \,\le \\
c_{10}  \sum_{n=0}^{\lfloor N/2 \rfloor}  K(M+n)  \sum_{m=\lfloor M/2 \rfloor +1}^{M} \frac{K_f(M-m)}{K_f(M)}\, 
\le\, c_{11} \frac{L(M)}{M^\ga} \frac{\sum_{m=0}^M K_f(m) }{K_f(M)}  \, 
.
\end{multline} 
If $ \overline{\alpha} \in  [1, (1+\alpha)/2)$ and $\sum_m K_f(m)<\infty$ then $T_3=O(N^{\overline{\alpha} - \ga}L(N)/ \overline{L}(N))$ and it tends to zero because $\ga>1$ implies $\ga >\overline{\ga}$.
If instead $\overline{\alpha} \le 1$ and $\sum_m K_f(m)=\infty$ the last ratio in the rightmost term in \eqref{eq:T3.g}
is $O(M)=O(N)$. Hence $T_3= O(N^{1-\ga})$ and $T_3$ tends to zero also in this case. 
The proof of \eqref{eq:sharpZ_f} is therefore complete.

\medskip 

Now for $ \overline{\alpha} > (1+\alpha)/2$ (which implies that $\overline{\ga} >1$ since $\ga >1$), we write 
\begin{equation}
\label{eq:ff1}
\frac{Z^f_{N,M,h}}{\mathrm{K}(N,M)} = \sum_{n=0}^N \sum_{m=0}^M Z^c_{n,m,h} \frac{K_f(N-n)K_f(M-m)}{\mathrm{K}(N,M)} \,.
\end{equation} 

We split the last sum to $U_1+U_2+U_3+U_4$ as in \eqref{eq:T4}. Since $\sum_{i \ge 0} \sum_{j \ge 0} Z^c_{i,j,h} = 1/(1 - \exp(h))< \infty$ (see \eqref{eq:sumZ})  we have 
\begin{equation}
U_1 \, =\, O \left(
\frac{K_f(N)K_f(M)}{\mathrm{K}(N,M)} 
\right)\, =\, O \left( N^{1 + \ga - 2 \overline{\ga} } \frac{\left( \overline{L}(N)\right)^2}{L(N)}\right) 
\, =\, o(1) \,.
\end{equation}
Also the terms $U_3$ and $U_4$ give a vanishing contribution. Let us see it for $U_3$ (the computation is identical for $U_4$):
\begin{multline}
U_3 \,=\,  \sum_{n=0}^{\lfloor N/2 \rfloor} \sum_{m=\lfloor M/2 \rfloor +1}^{M} \frac{ Z^c_{n,m,h}}{\mathrm{K}(N,M)} K_f(N-n) K_f(M-m) \, \le \\ c_{12} \sum_{n=0}^{\lfloor N/2 \rfloor} K_f(N-n) \sum_{m=\lfloor M/2 \rfloor +1}^{M} K_f(M-m)
\, \le \, c_{13} \sum_{n \ge N/2} K_f(n)\, =\, O( \overline{L}(N) N^{1 - \overline{\ga}})
\,.
\end{multline} 
The relevant contribution comes from $U_2$:  
\begin{equation}
\begin{split}
U_2 \,&=\,  \sum_{n=\lfloor N/2 \rfloor+1}^{N} \sum_{m=\lfloor M/2 \rfloor +1}^{M} \frac{ Z^c_{n,m,h}}{\mathrm{K}(N,M)} K_f(N-n) K_f(M-m) \\
&=\, \sum_{n=0}^{N-\lfloor N/2 \rfloor-1} \sum_{m=0}^{M-\lfloor M/2 \rfloor -1} \frac{ Z^c_{N-n,M-m,h}}{\mathrm{K}(N,M)} K_f(n) K_f(m)\, . 
\end{split}
\end{equation}
But the ratio in the last term is bounded, cf.  \eqref{eq:sharpZ_c}, and in fact  
 \eqref{eq:sharpZ_c} tells us also that for every $m$ and $n$
\begin{equation}
\lim_{N,M \to \infty} \frac{Z^c_{N-n,M-m,h}}{\mathrm{K}(N,M)}\, =\, \frac{ \exp(h)} {(1- \exp(h))^2}\, ,
\end{equation}
which, by applying (DOM),  implies
\begin{equation}
\lim_{N,M \to \infty } U_2 = 
 \frac{\exp(h)}{{\left( 1 - \exp(h) \right)}^{2}} {\left( \sum_{n \ge 0} K_f(n) \right)}^2 \,,
\end{equation}
and  completes the proof of \eqref{eq:sharpZ_f1} and, in turn, the proof of Proposition~\ref{th:SharpZ}.
\qed

\medskip

\subsection{Path properties: proof of Theorem~\ref{th:pathsharp}}

\label{sec:pathproofs}

In this section, we suppose that $M \sim \gamma N$ and $\ga >1$.

\medskip

{\it Proof of (1) in Theorem~\ref{th:pathsharp}.}
We first consider the case $h<0$.
Recall that $(\cF_1,\cF_2)$% := \max \lbrace \tau \cap [0,N] \times [0,M] \rbrace $ 
 is the last renewal epoch in $[0,N]\times [0,M]$. If $\overline{\ga} < (1+ \ga)/2$, for 
fixed $i$ and $j$ (so we can assume
$i<N$ and $j<M$) 
 we have
\begin{equation}
\label{eq:maxf}
\begin{split}
\bP^f_{N,M,h} \left( (\cF_1,\cF_2) =(i,j) \right) & = \frac{Z^c_{i,j,h} K_f(N-i) K_f(M-j)}{Z^f_{N,M,h}} \\
& \hskip -.3 cm
\stackrel{N \to \infty}\sim (1-\exp(h))Z^c_{i,j,h} \frac{ K_f(N-i) K_f(M-j)}{K_f(N) K_f(M)} \,,
\end{split}
\end{equation}
where the estimation follows from \eqref{eq:sharpZ_f}. Since  $i$ and $j$ are $O(1)$
 \eqref{eq:psv} 
 the ratio in the rightmost term in \eqref{eq:maxf} converges to one. Hence it suffices  to prove that $ (1-\exp(h)) \sum_{i,j} Z^c_{i,j,h} =1$, but this is done in \eqref{eq:sumZ}.
%\begin{equation}
%\label{eq:sumz}
%\sum_{i \ge 0} \sum_{j \ge 0} Z^c_{i,j,h} = \frac{1}{1- \exp(h)} \,,
%\end{equation}
We the recall  that $Z_{i,j,h} = \bP \left( (i,j) \in \tilde \tau_h \right)$ and \eqref{eq:cF} is proven. 

Now recall that  $(\cL_1, \cL_2)\, :=\, \left(N- \cF_1,M-\cF_2 \right)$. If $\overline{\ga} > (1+ \ga)/2$, for fixed $i$ and $j$  %(in this case from \eqref{eq:sharpZ_f1} we see that the free ends are microscopics and a jump of order $N$ appears in the system, i.e. $i$ and $j$ are $o(N)$), 
by using \eqref{eq:sharpZ_c} and \eqref{eq:sharpZ_f1} we see that
\begin{multline}
\label{eq:cL1}
\bP^f_{N,M,h} \left( (\cL_1,\cL_2) =(i,j) \right) \, = \, \frac{Z^c_{N-i,M-j,h} K_f(i) K_f(j)}{Z^f_{N,M,h}} \\
\stackrel{N \to \infty} \sim \frac{\mathrm{K}(N-i,M-j)}{\mathrm{K}(N,M)} \frac{K_f(i) K_f(j)}{{\left( \sum_{n \ge 0} K_f(n) \right)}^2} \,.
\end{multline}
The proof of \eqref{eq:cF-2} is therefore complete by observing that by \eqref{eq:psv} the first ratio in \eqref{eq:cL1} converges to one.

We are then left (for $h<0$) with \eqref{eq:L}. Here we  prove  more: consider
 $(\cE_1,\cE_2) := \max\{\tau \cap [0, \lfloor N/2 \rfloor ]\times[0,\lfloor M/2 \rfloor]\}$ under $\bP^f_{N, M, h}$ for $\overline{\ga} > (1+ \ga)/2$. If $i$ and $j$ are fixed, by \eqref{eq:sharpZ_f1} we have 
\begin{multline}
\label{eq:cE}
\bP^f_{N,M,h} \left( (\cE_1,\cE_2) =(i,j) \right) \,  =\,  \exp(h) \sum_{s \ge \lfloor N/2 \rfloor} \sum_{t \ge \lfloor M/2 \rfloor } \frac{Z^c_{i,j,h} \mathrm{K}(s-i,t-j) Z^f_{N-s,M-t,h}}{Z^f_{N,M,h}}
\\
 \stackrel{N \to \infty}\sim \frac{(1-\exp(h))^2 Z^c_{i,j,h}}{{ \left( \sum_{n \ge0} K_f(n)\right)}^{2}} \sum_{s \ge \lfloor N/2 \rfloor} \sum_{t \ge \lfloor M/2 \rfloor } Z^f_{N-s,M-t,h}  \frac{\mathrm{K}(s-i,t-j)}{\mathrm{K}(N,M)} \,.
\end{multline}
By making the change of variable $(s,t) \to (N-s, M-t)$ we see that   for 
 $s$ and $t$ of $O(1)$, the very last ratio in \eqref{eq:cE} converges to one and the same ratio is bounded by \eqref{eq:psv}
 in all the range of the sum.
 Hence, by the (DOM), the expression in \eqref{eq:cE} converges to 
\begin{equation}
\frac{(1-\exp(h))^2 Z^c_{i,j,h}}{{ \left( \sum_{n \ge0} K_f(n)\right)}^{2}} \sum_{s,t \ge 0} Z^f_{s,t,h} \,,
\end{equation}
and observe that
\begin{multline}
\label{eq:sumzf}
\sum_{s,t \ge 0} Z^f_{s,t,h} = \sum_{s,t \ge 0} \sum_{n=0}^s \sum_{m=0}^t Z^c_{n,m,h} K_f(s-n) K_f(t-m) \\
= \sum_{n,m \ge 0} Z^c_{n,m,h} \sum_{s \ge n} \sum_{m \ge t} K_f(s-n) K_f(t-m) = \frac{{(\sum_{n \ge0} K_f(n))}^{2}}{(1-\exp(h))}  \,,
\end{multline}
where the last equality follows from \eqref{eq:sumZ}. Therefore the law of $(\cE_1,\cE_2)$ under $\bP^f_{N,M,h}$ converges for $N \to \infty$ to the probability distribution that assigns  to $(i,j)$ probability $(1-\exp(h)) Z^c_{i,j,h}$ (which is correctly normalized, like  \eqref{eq:cF}, by \eqref{eq:sumZ}).

Let $(\cC_1,\cC_2) := \min \{\tau \cap [ \lfloor N/2 \rfloor,N ] \times [ \lfloor M/2 \rfloor,M ]\} $ be the first renewal epoch in $[\lfloor N/2 \rfloor,N ] \times [ \lfloor M/2 \rfloor,M ] $ and set
 $(\cH_1,\cH_2) := \left( N - \cC_1,M-\cC_2 \right)$. For fixed $i$ and $j$ by \eqref{eq:sharpZ_f1} we have 
\begin{equation}
\label{eq:cH}
\begin{split}
\bP^f_{N,M,h} \left( (\cH_1,\cH_2) =(i,j) \right) & = \exp(h) \sum_{s \ge \lfloor N/2 \rfloor} \sum_{t \ge \lfloor M/2 \rfloor } \frac{Z^f_{i,j,h} \mathrm{K}(s-i,t-j) Z^c_{N-s,M-t,h}}{Z^f_{N,M,h}}\\
& \stackrel{N \to \infty}\sim \frac{(1-\exp(h))^2 Z^f_{i,j,h}}{{ \left( \sum_{n \ge0} K_f(n)\right)}^{2}} \sum_{s \ge \lfloor N/2 \rfloor} \sum_{t \ge \lfloor M/2 \rfloor } Z^c_{N-s,M-t,h}  \frac{\mathrm{K}(s-i,t-j)}{\mathrm{K}(N,M)} \,.
\end{split}
\end{equation}
The second ratio converges to $1$ (same argument as above) and  therefore  \eqref{eq:cH} converges to
\begin{equation}
\label{eq:pH}
\frac{(1-\exp(h)) Z^f_{i,j,h}}{{ \left( \sum_{n \ge0} K_f(n)\right)}^{2}}  \,,
\end{equation}
and from \eqref{eq:sumzf}, we see that this expression adds up ($i, j\ge 0$) to one.
%\begin{equation}
%\sum_{i,j \ge 0} \frac{(1-\exp(h)) Z^f_{i,j,h}}{{ \left( \sum_{n \ge0} K_f(n)\right)}^{2}} =1\,.
%\end{equation}
The law of $(\cH_1,\cH_2)$ converges as $N \to \infty$ to the probability distribution that assigns  to $(i,j)$ probability \eqref{eq:pH}.
We conclude that the contacts are either close to $(0,0)$ or to the last renewal epoch: therefore we get \eqref{eq:L}.

\qed

\medskip

\noindent
{\it Proof of (2) in Theorem~\ref{th:pathsharp}.}
In this case $h >0$ and $\gamma \in (1/\gamma_c(h), \gamma_c(h))$:
the positivity of $\tf_\gamma(h) - \gamma \partial_{\gamma} \tf_{\gamma}(h)$ and $\partial_{\gamma} \tf_{\gamma}(h)$
is a direct consequence of the strict concavity of $\tf_\cdot(h)$ in the Cram\'er region (see caption of Figure~\ref{fig:DandF}).
Then choose  $(i,j)$ with non-negative integer entries. By \eqref{eq:Zc-sharp} and \eqref{eq:sharpZf} we have
\begin{multline}
\label{eq:cL}
 \bP^f_{N,M,h} \left( (\cL_1,\cL_2) = (i,j) \right) \, = \, \frac{Z^c_{N-i,M-j,h} K_f(i) K_f(j)}{Z^f_{N,M,h}} 
\stackrel{N \to \infty} \sim 
\\
\frac{A \left( (M-j)/(N-i) \right) \sqrt{N} }{C_{\gamma,h} A(\gamma) \sqrt{N-i}}   \exp \left( \tf_{(M-j)/(N-i)} (h) (N-i) - \tf_{M/N}(h) N \right)K_f(i) K_f(j)  \,.
\end{multline} 
 Observe now that the ratio in the rightmost term  of \eqref{eq:cL} converges as $N \to \infty$ to $1/C_{\gamma,h}$ (defined in \eqref{eq:Cgah}) and using \eqref{eq:Zfe} and the fact that $(M-j)/(N-i)$ is close to $\gamma$, \eqref{eq:cL} converges to
\begin{equation}
  \frac{1}{C_{\gamma,h}} \exp \left( - j \partial_{\gamma} \tf_{\gamma}(h) - i  \left( \tf_\gamma(h) - \gamma \partial_{\gamma} \tf_{\gamma}(h) \right) \right) K_f(i) K_f(j) \,,
\end{equation}
where $C_{\gamma,h}$ is defined in \eqref{eq:Cgah}. Hence \eqref{eq:cgammah8} is proven.

\medskip

To complete the proof of Theorem~\ref{th:pathsharp}, we need to show that for $h >0$ and if $ M \sim \gamma N$ such that $ \gamma \in (1/ \gamma_c(h), \gamma_c(h))$, we have
$
\bP^f_{N,M,h}$ converges, for $N\to \infty$, to the law of a bivariate renewal
with the inter-arrival probability given in \eqref{eq:cgammah9}. 
For this, fix a
$k \in \bbN$  and for every $\{(i_n,j_n)\}_{n=0,1, \ldots,k}$ with $(i_0,j_0)=(0,0)$, $i_1 < i_2< \dots < i_k < N$ and $j_1 < j_2 <...< j_k < M$, we have 
\begin{multline}
\label{eq:Pcf}
\bP^f_{N, M, h}\left(\tau_n= (i_n,j_n) \text{ for } n=1,2, \ldots, k \right)
 = \\
\left( \prod_{n=1}^k  K(i_n-i_{n-1} + j_n-j_{n-1})\right) \exp(h) \frac{Z^f_{N-i_k,M-j_k,h}}{Z^f_{N,M,h}} \,.
\end{multline}
Since $\gamma_k = (M-j_k)/(N-i_k)$ close to $\gamma$ (in the Cram\'er region), by \eqref{eq:sharpZf}, we see that the ratio in \eqref{eq:Pcf} is equal to
\begin{equation}
1+o(1))
\exp \left( - N ( \tf_{\gamma_N}(h) - \tf_{\gamma_k}(h) + \frac{i_k}{N} \tf_{\gamma_k}(h)  ) \right)\,. 
\end{equation} 
Following the same procedure used for \eqref{eq:Zfe}, we see that the exponent converges to 
\begin{equation}
\label{eq:exp}
  - i_k( \tf_\gamma(h) - \gamma \partial_\gamma \tf_\gamma(h)) -j_k \partial_\gamma \tf_\gamma(h)  \,,
\end{equation}
and therefore the left-hand side of \eqref{eq:Pcf} converges to 
\begin{equation}
\prod_{n=1}^k  K(i_n-i_{n-1} + j_n-j_{n-1})  \exp \left(h  - (i_n-i_{n-1})( \tf_\gamma(h) - \gamma \partial_\gamma \tf_\gamma(h)) - (j_n-j_{n-1}) \partial_\gamma \tf_\gamma(h) \right) \,.
\end{equation}
We are left with proving that \eqref{eq:cgammah9} is a probability distribution. Recall from \eqref{eq:FgcDh} and \eqref{eq:Bh} that
\begin{equation}
D_h(1,\gamma) = \max_{(\lambda_1,\lambda_2) \in B_h} (\lambda_1 + \gamma \lambda_2) = \hat{\lambda}_1(\gamma) + \gamma \hat{\lambda}_2(\gamma) \,,
\end{equation}
with of course $(\hat{\lambda}_1(\gamma),\hat{\lambda}_2(\gamma)) \in B_h$ (keep in mind that they depend also on $h$) i.e. 
\begin{equation}
\label{eq:hat12}
\sum_{n,m} K(n+m)\exp \left( h   - n (\tg(h) -\hat{\lambda}_1(\gamma)) - m (\tg(h) -\hat{\lambda}_2(\gamma)) \right) = 1 \,,
\end{equation}
therefore $\tf_\gamma(h) = (\tg(h) - \hat{\lambda}_1(\gamma)) + \gamma ( \tg(h) - \hat{\lambda}_2(\gamma))$. Replace $\tg(h) -\hat{\lambda}_1(\gamma)$ in \eqref{eq:hat12} by $\tf_{\gamma}(h) - \gamma ( \tg(h) - \hat{\lambda}_2(\gamma))$, we obtain 
\begin{equation}
\label{eq:hat12_1}
\sum_{n,m} K(n+m)\exp \left( h - n \tf_{\gamma}(h) -  (m-\gamma n) (\tg(h) -\hat{\lambda}_2(\gamma))  \right) = 1 \,,
\end{equation} 
and recall from \eqref{eq:ID1} that 
\begin{equation}
\label{eq:ID11}
\sum_{n,m} (m- \gamma n ) K(n+m) \exp \left( h - n \tf_{\gamma}(h) -  (m-\gamma n) (\tg(h) -\hat{\lambda}_2(\gamma))  \right) = 0 \,.
\end{equation}
Differentiating \eqref{eq:hat12_1} with respect to $\gamma$ and using \eqref{eq:ID11}, we obtain that
\begin{equation}
\left( \partial_\gamma \tf_\gamma(h) - (\tg(h) - \hat{\lambda}_2(\gamma)) \right) \sum_{n,m} n K(n+m) \exp \left( - n \tf_{\gamma}(h) - (m-\gamma n) (\tg(h) -\hat{\lambda}_2(\gamma))  \right) \,,
\end{equation}
is zero, and this
 implies directly that 
\begin{equation}
\partial_\gamma \tf_\gamma(h) = \tg(h) - \hat{\lambda}_2(\gamma) \,,
\end{equation}
\begin{equation}
\tf_\gamma(h) - \gamma \partial_\gamma \tf_\gamma(h) = \tg(h) - \hat{\lambda}_1(\gamma) \,.
\end{equation}
Therefore from \eqref{eq:hat12}, we get that \eqref{eq:cgammah9} is a probability distribution.
%\begin{equation}\sum_{i,j} K(i+j) \exp \left(h -i (\tf_\gamma(h) - \gamma \partial_\gamma \tf_\gamma(h)) - j \partial_\gamma \tf_\gamma(h) \right)=1 \,.\end{equation}
\qed 

\medskip

\appendix

\section{Slowly and regularly varying functions}
\label{sec:A}

 $L:[0, \infty)\to (0, \infty)$ is a {\sl slowly varying} function at $\infty$ if is measurable and if $\lim_{x\to \infty} \frac{L(ux)}{L(x)}=1$ for every $u>0$.
 The function $x \mapsto L(1/x)$ is slowly varying at zero if $L(\cdot)$ is slowly varying at $\infty$.
  It can be shown that this convergence holds uniformly in $u$ \cite[Th.~1.2.1]{cf:BGT}: for every $0< c_1 <c_2 < \infty$ 
\begin{equation}
\label{eq:psv}
\lim_{x\to \infty} \sup_{u \in [c_1,c_2]} \left| \frac{L(ux)}{L(x)} - 1 \right| =0 \,.
\end{equation}
A function of the form $x \mapsto x^a L(x)$, $a\in \bbR$, is said to be {\sl regularly varying (at $\infty$) of exponent $a$}. Analogous definition for regularly varying at zero. 
 Examples of slowly varying function include logarithmic functions (of course the trivial example is the constant) like $ a (\log(x))^b$ as $ x \to \infty$ with $ a >0$ and $b \in \bbR$. 
We refer to \cite{cf:BGT} for the full theory of slowly and regularly varying functions: we just recall some basic important facts. 
First of all that both $L(x)$ and $1/L(x)$ are $o(x^\gep)$ for every $\gep>0$, which directly implies that if $f(\cdot)$ is regularly
varying with exponent $a$ and $g(\cdot)$ is regularly varying with exponent $b<a$, then $g(x)=o(f(x))$. 

We will often use  that for $\gb>0$ \cite[Sec.~1.5.6]{cf:BGT}
 \begin{equation}
 \label{eq:recall}
 \sum_{n \ge N} \frac{L(n)}{ n^{1+\beta}} 
 \stackrel{N \to \infty}\sim \frac{L(N)} { \beta  N^{\beta  }}\
 \ \
 \text{ and } \ \ \ 
 \sum_{n=1}^N \frac{L(n)}{ n^{1-\beta}} \stackrel{N \to \infty}\sim \frac{L(N)} { \gb N^{-\beta}} \, ,
 \end{equation} 
which can be proven by Riemann sum approximation. We will often use Riemann sum approximations involving regularly
varying functions also beyond \eqref{eq:recall} and the central tool to control these approximations are the so called 
{\sl Potter bounds} \cite[Th.~1.5.6]{cf:BGT}.

Another important issue is about asymptotic invertibility of regular functions: a regular function of exponent $a>0$ (respectively $a<0$) is asymptotically
equivalent to an increasing  (respectively decreasing) function. Moreover the inverse of a monotonic regularly varying function of exponent $a\neq 0$
is a regularly varying function of exponent $1/a$. In different terms, if $f(\cdot)$ is regularly varying of exponent $a\neq 0$, then there exists 
$g(\cdot)$ regularly varying of exponent $1/a$ such that $f(g(x))\sim g(f(x))\sim x$
\cite[Sec.~1.5.7]{cf:BGT}. Occasionally we use  other properties of slowly varying functions and we refer directly to 
\cite{cf:BGT}.


\begin{thebibliography}{99}

%\bibitem{cf:B}
%J. Bertoin, \emph{Renewal theory for embedded regenerative sets}, Ann.Probab. {\bf 27} (1999), 1523-1535

\bibitem{cf:BePo}
Q. Berger and  J. Poisat, \emph{On the critical curve of the pinning and copolymer models in correlated Gaussian environment}, Elect. J. Probab. {\bf 20} (2015), 35 pages. 

\bibitem{cf:BGT}
N. H. Bingham, C. M. Goldie and J. L. Teugels, \emph{Regular variation}, Cambridge University Press, Cambridge (1987).

\bibitem{cf:BM}
A. A. Borovkov and A. A. Mogul'skil', \emph{The second rate function and the asymptotic problems of renewal and hitting the boundary for multidimensional random walks}, Sibirsk. Mat. Zh. {\bf 37} (1996) No. 4, 745-782.


\bibitem{cf:BM1}
A. A. Borovkov and A. A. Mogul'skil', \emph{Integro-local limit theorems including large deviations for sums of random vectors \uppercase\expandafter{\romannumeral 2}}, Theory Probab.Appl. {\bf 45} (2001), 3-22.

\bibitem{cf:BM2}
A. A. Borovkov and A. A. Mogul'skil', \emph{On large deviations of sums of independent random vectors on the boundary and outside oh the cram\'er zone \uppercase\expandafter{\romannumeral 1}}, Theory Probab.Appl. {\bf 53} (2009), No. 2, 301-311.

\bibitem{cf:BOPG} R.~Brak, A.~L.~Owczarek, T.~Prellberg and A.~J~Guttmann, \emph{Finite-length scaling of collapsing directed walks}, Phys. Rev. E {\bf 48} (1993), 2386-2396.

\bibitem{cf:BundHwa} R. Bundschuh and T. Hwa, \emph{Statistical mechanics of secondary structures formed by random RNA sequences}, Phys. Rev. E {\bf 65} (2002), 031903 (22 pages).

\bibitem{cf:CGZ}
F. Caravenna, G. Giacomin and L. Zambotti,
\emph{Sharp asymptotic behavior for wetting models in (1+1)-dimension}, 
Elect. J. Probab. {\bf 11} (2006), 345-362.

\bibitem{cf:CPN}
P.~Carmona, N.~P\'etr\'elis and G.~B.~Nguyen,  \emph{Interacting partially self-avoiding walk, from phase transition to the geometry of the collapse phase}, Ann. Probab., to appear, http://arxiv.org/abs/1306.4887

\bibitem{cf:Doney} 
 R. A. Doney, \emph{An analogue of the renewal theorem in higher dimensions},
Proc. London Math. Soc. (3) {\bf 16} (1966), 669-684.

\bibitem{cf:EON}
T.~R.~Einert, H.~Orland and R.~R.~Netz,
\emph{Secondary structure formation of homopolymeric single-stranded nucleic acids including force and loop entropy: implications for DNA hybridization}, Eur. Phys. J. E {\bf 34} (2011), 55 (15 pages).

\bibitem{cf:Fisher}
M.~E.~Fisher, \emph{Walks, walls, wetting, and melting},
J. Statist. Phys. {\bf 34} (1984), 667-729.

\bibitem{cf:GO0}
T.~Garel and H.~Orland, \emph{On the role of mismatches in DNA denaturation}, arXiv:cond-mat/0304080

\bibitem{cf:GO} T.~Garel and H.~Orland, \emph{Generalized Poland-Scheraga model for DNA hybridization},  Biopolymers {\bf 75} (2004), 453-467.


\bibitem{cf:GB} G. Giacomin, \emph{Random polymer models}, 
Imperial College Press, World Scientific, 2007. 

\bibitem{cf:GB-ejp} G. Giacomin, 
\emph{Renewal convergence rates and correlation decay for homogeneous pinning models}, 
Elect. J. Probab. {\bf 13} (2008), 513-529. 
 
\bibitem{cf:G} G. Giacomin, \emph{Disorder and critical phenomena through basic probability models},Lecture notes from the 40th Probability Summer School held in Saint-Flour, 2010. Lecture Notes in Mathematics {\bf  2025}, Springer, 2011.

\bibitem{cf:GT-alea} G. Giacomin and F. L. Toninelli, \emph{The localized phase of disordered copolymers with adsorption}, 
ALEA-Latin American J. Probab. Math. Stat. {\bf 1} (2006), 149-180.

\bibitem{cf:dH}
F.~den Hollander, \emph{Random polymers}, Lectures from the 37th Probability Summer School held in Saint-Flour, 2007. Lecture Notes in Mathematics {\bf 1974},  Springer-Verlag, 2009. 

\bibitem{cf:KBM}
A.~Kabak\c{c}\i $\breve{\mathrm{o}}$glu, A. Bar and D. Mukamel,
\emph{Macroscopic loop formation in circular DNA denaturation},
Phys. Rev. E {\bf 85} (2012), 051919.


\bibitem{cf:KMP}
Y. Kafri, D. Mukamel and L. Peliti, \emph{Why is the DNA denaturation transition first order?}, 
Phys. Rev. Lett. {\bf 85} (2000), 4988-4991.


\bibitem{cf:Primer}
 S.~G.~Krantz and H.~R.~Parks, \emph{A primer of real analytic functions},  Second edition, Birkh\"auser Advanced Texts: : Basel Textbooks. Birkh\"auser Boston, Inc., Boston, MA, 2002.

\bibitem{cf:litan} 
A.~Litan, 
\emph{A Statistical Mechanical Treatment of the Open Ends in a Double-Stranded Polynucleotide Molecule},
Biopolymers {\bf 2} (1964), 279-282.

%\bibitem{cf:MR} S. Mercier and M. Roussignol, \emph{A block replacement policy for a bivariate wear subordinator}, 

\bibitem{cf:NG} R.  A. Neher and U. Gerland,
\emph{Intermediate phase in DNA melting}, Phys. Rev. E {\bf 73} (2006), 030902R.

\bibitem{cf:NP} 
G.~B.~Nguyen and N.~P\'etr\'elis, \emph{A variational formula for the free energy of the partially directed polymer collapse}, J. Stat. Phys. {\bf 151} (2013), 1099-1120.


\bibitem{cf:PSbook} D. Poland and H. A. Scheraga, \emph{Theory of helix-coil transitions in biopolymers;: Statistical mechanical theory of order-disorder transitions in biological macromolecules}, Academic Press, 1970.

%\bibitem{cf:PS1-2} D. Poland and H. A. Scheraga, \emph{Occurrence of a phase transition in nucleic acid models}, J. Chem. Phys. {\bf 45} (1966), 1464-1469.

\bibitem{cf:RG} 
C. Richard and A. J. Guttmann,
\emph{Poland-Scheraga Models and the DNA Denaturation
Transition}, 
J.  Statist. Phys. {\bf 115} (2004), 925-947.

\bibitem{cf:SNR}
K. Sture J. Nordholm and S. A. Rice, \emph{Phase transitions and end effects in models of biopolymers},
J. Chem. Phys. {\bf 59} (1973), 5605-5613.

\bibitem{cf:Yer}
E. Yeramian,
\emph{Complexity and Tractability. Statistical Mechanics of Helix-Coil Transitions in Circular DNA as a Model-Problem},
Europhysics Letterr {\bf 25} (1964), 49-55. 


%\bibitem{cf:OW} E. Omey and E. Willekens, \emph{Abelian and Tauberian theorems for the Laplace transform in several variables}, J. Multivariate Anal. {\bf 30} (1989), 292-306.



\end{thebibliography}
\end{document}